\documentclass[11pt,reqno]{amsart}
\usepackage{hyperref}

\allowdisplaybreaks[4]
\usepackage{amsmath}
\usepackage{amssymb}
\usepackage{mathrsfs}
\usepackage{amsthm}
\usepackage{bm}
\usepackage{booktabs}
\usepackage{float}
\usepackage{graphicx,epstopdf}

\usepackage[caption=false]{subfig}

\graphicspath{{Fig/}}
\usepackage{geometry}

\usepackage{color}
\theoremstyle{plain}
\newtheorem{lemma}{Lemma}[section]
\newtheorem{theorem}[lemma]{Theorem}

\newtheorem{corollary}[lemma]{Corollary}
\newtheorem{remark}[lemma]{Remark}
\newtheorem{prop}[lemma]{Proposition}
\newtheorem{assumption}{Assumption}
\begin{document}
\title[Convergence of adaptive time-stepping schemes]{Strong convergence of adaptive time-stepping schemes for the stochastic Allen--Cahn equation}
\author{Chuchu Chen, Tonghe Dang, Jialin Hong}
\address{LSEC, ICMSEC,  Academy of Mathematics and Systems Science, Chinese Academy of Sciences, Beijing 100190, China,
\and 
School of Mathematical Sciences, University of Chinese Academy of Sciences, Beijing 100049, China}
\email{chenchuchu@lsec.cc.ac.cn; dangth@lsec.cc.ac.cn; hjl@lsec.cc.ac.cn}
\thanks{This work is funded by National Natural Science Foundation of China (No. 11971470,
No. 11871068, No. 12031020, No. 12022118), and by Youth Innovation Promotion Association CAS}
\begin{abstract}
It is known in \cite{beccari} that the standard explicit Euler-type scheme (such as the exponential Euler and the linear-implicit Euler schemes) with a uniform timestep, though computationally efficient, may diverge for the stochastic Allen--Cahn equation. To overcome the divergence, this paper proposes and analyzes adaptive time-stepping schemes, which adapt the timestep at each iteration to control numerical solutions from instability.
The \textit{a priori} estimates  in $\mathcal {C}(\mathcal {O})$-norm and $\dot{H}^{\beta}(\mathcal{O})$-norm of numerical solutions are established provided the adaptive timestep function is suitably bounded, which plays a key role in the convergence analysis. We show that 
the adaptive time-stepping schemes converge strongly with order $\frac{\beta}{2}$ in time and $\frac{\beta}{d}$ in space with 
$d$ ($d=1,2,3$) being the dimension and 
$\beta\in(0,2]$. 
Numerical experiments show that the adaptive time-stepping schemes are simple to implement and 
at a lower computational cost than a scheme with the uniform timestep. 
\end{abstract}
\maketitle
\section{Introduction}
Numerical approximations for stochastic partial differential equations (SPDEs) with globally Lipschitz coefficients have been studied in recent decades (see e.g., the monograph \cite{LordBook}). 
In contrast, numerical analysis of SPDEs with non-globally Lipschitz coefficients, for example the stochastic Allen--Cahn equation, has been considered (see e.g., \cite{CH19, CHS21, LQ20} and references therein) and is still not fully understood. 
It is pointed out in \cite{beccari} that the explicit Euler, the exponential Euler and the linear-implicit Euler schemes with the uniform timestep fail to converge in the strong sense for SPDEs with superlinearly growing coefficients; see also \cite{jentzen}. 
Implicit schemes like fully drift-implicit scheme (see e.g. \cite{KLL18, QW19, LQ20, liuzh, prohl18} and references therein) can be strongly convergent in this setting. It is known that the implementation of the implicit scheme requires solving an algebraic equation at each iteration step, which needs additional computational effort.
These reasons have led to the research on the construction of explicit schemes that can ensure  convergence under the non-globally Lipschitz condition. For instance,  \cite{Split19} proposes the splitting scheme and studies the convergence in strong, weak, and  probability senses. It is shown that the  mean-square convergence  order is almost 1/4, localized on an event of arbitrarily large probability, and that the  convergence order in probability  is almost 1/4 for the space-time white noise case; Authors in \cite{BCH2019} study the  strong convergence order of the explicit temporal splitting numerical scheme for the case of  
different noises with varying degrees of smoothness; \cite{caimeng21} proves that the weak convergence order of the tamed exponential Euler scheme is almost $\beta$  for the generalized $Q$-Wiener process case, and \cite{WXJ20} shows that the strong convergence order of the nonlinearity-tamed accelerated exponential Euler  scheme is almost $1/2$ for the space-time white noise case;  \cite{BJ19} studies the strong convergence order of the nonlinearity-truncated Euler-type  schemes, which is almost $1/4$ for the cylindrical Wiener process case. The present work makes further contributions on the numerical study of the explicit, adaptive time-stepping schemes for the stochastic Allen--Cahn equation.

Adaptive time-stepping schemes, which adapt the timestep at each iteration to control the numerical solution from divergence,
have been deeply studied for stochastic ordinary differential equations (SODEs) 
with non-globally Lipschitz drift. As for the selection of adaptive timesteps, we refer to e.g. \cite{Fang20, Lord18, LordK22} for the admissible strategy, \cite{Lamba, Lemaire07} for the strategy based on the local error control, and \cite{Prohl} for the strategy based on the  \textit{a posteriori} weak error estimate. 
Numerically this adaptive scheme is simple to implement and the complexity is similar as that of an Euler scheme, which is a big advantage for high dimensional problems (see \cite{Lamba}). 
 It is also pointed out  in \cite{Fang20, MLMC2012} that such an adaptive scheme can lead to better computational performance for multi-level Monte--Carlo simulations.
 To our knowledge, 
 there are few works on the study of the adaptive time-stepping scheme for SPDEs.
 The first attempt to apply the adaptive time-stepping scheme to the simulation of SPDEs is \cite{lord}, where the strong convergence rate is obtained under the assumption that the Fr\'echet derivative of $F$ is bounded
   polynomially in $L^2(\mathbb{D})$-norm (see \cite[Assumption 2.4]{lord}), where $\mathbb{D}\subset \mathbb{R}^d$.  
   
  Consider another important class of nonlinear SPDEs, including the stochastic Allen--Cahn equation driven by additive noise
\begin{align}
\begin{cases}\label{spde}
\mathrm{d}X(t)+AX(t)\mathrm{d}t=F(X(t))\mathrm{d}t+\mathrm{d}W(t),\;t\in(0,T],\\
X(0)=X_0,
\end{cases}
\end{align}
where $-A:=\Delta:\mathrm{Dom}(A)\subset H\rightarrow H$ is the Laplacian 
operator with homogeneous Dirichlet boundary condition with 
$H:=L^2(\mathcal{O}),\mathcal{O}:=[0,1]^d,d=1,2,3$ endowed with the usual inner product $\langle\cdot,\cdot\rangle$ and the norm $\|\cdot\|,$
  and the stochastic process $\{W(t)\}_{t\in[ 0,T]}$ is a generalized $Q$-Wiener process on a filtered probability space $(\Omega,\mathcal{F},\left\{\mathcal{F}_t\right\}_{0\leq t\leq T},\mathbb{P})$, subject to  $\|A^{\frac{\beta-1}{2}}Q^{\frac12}\|_{\mathcal{L}_2(H)}<\infty,\,\beta\in(0,2]$.
The nonlinear drift $F$ is a Nemytskii operator defined by $F(X)(x)=f(X(x))$ for $X\in H$ and $x\in\mathcal{O}$, where $f$ is a polynomial and satisfies $f(\xi)=\sum_{i=0}^3a_i\xi^i,\; a_3<0,\,\xi\in\mathbb{R}$.
 In this case, the Fr\'echet derivative of $F$ is bounded polynomially in $E$-norm ($E:=\mathcal{C}(\mathcal{O})$).
  The main contribution of this work is to present the \textit{a priori} estimates and rigorous strong convergence analysis
of the adaptive time-stepping schemes for \eqref{spde}.

To be specific, the adaptive time-stepping scheme for \eqref{spde}, whose spatial discretization is using the spectral Galerkin method, and temporal direction is based on the adaptive exponential integrator, is an explicit numerical scheme with adaptive timesteps. The prerequisite of the convergence analysis is the \textit{a priori} estimates in $E$-norm and $\dot{H}^{\beta}$-norm of the numerical solution,  which are derived by a bootstrap argument. We refer to \cite{WXJ20} for the use of this argument for the nonlinearity-tamed scheme with the uniform timestep. 
Based on the above \textit{a priori} estimates,  and combining the smoothing effect of the analytic semigroup and regularity properties of the generalized $Q$-Wiener process, the strong convergence order of this fully discrete scheme for \eqref{spde} is finally carefully analyzed, which is the same as usual, i.e., order $\frac{\beta}{2}$ in time and $\frac{\beta}{d}$ in space. Moreover, we give the numerical analysis of the adaptive time-stepping scheme for the multiplicative noise case, and show that the convergence order is $1/2$ in time and $1$ in space when $d=1$. 

For the feasibility of an adaptive time-stepping scheme, some bounds for the adaptive timestep function are proposed. 
We would like to mention that in practice, instead of  verifying whether an adaptive timestep function satisfies the given lower bound, 
 one generally  introduces   a backstop scheme with a uniform timestep and couple it with the adaptive time-stepping scheme to ensure that a simulation over the interval $[0,T]$ can be completed in a finite number of timesteps (see \cite{Lord18} for the case of SODEs). More precisely, when the lower bound is invalid for the adaptive timestep function, for example, we perform a single step with the tamed exponential integrator with a uniform timestep instead. It can be shown that the corresponding  coupled scheme is strongly convergent with the order being the same as the adaptive time-stepping scheme. Further  it can be observed from numerical experiments  in Section \ref{sec3} that the coupled schemes are at a  lower computational cost, measured in terms of the CPU time. 
 %These show the  efficiency of the coupled schemes  than a scheme with the uniform timestep.

The outline of this paper is as follows. In the next section, some preliminaries are listed. In Section \ref{sec2}, we propose the adaptive time-stepping schemes, and present the main convergence theorem of this paper. Section \ref{sec2.1} presents 
the \textit{a priori} estimates in $E$-norm and $\dot{H}^{\beta}$-norm of numerical solutions. In Section \ref{secproof}, we give the proof of the main convergence theorem of the schemes. In Section \ref{multi_discuss}, we give the discussion of the numerical analysis for the multiplicative noise case. 
 Section \ref{sec3} is devoted to the numerical experiments, which verify our theoretical results.

\section{Preliminaries}\label{preliminary}

In this section, we give assumptions on $A,\, F,\,W(t),$ and the initial datum, as well as the well-posedness of \eqref{spde}, see e.g. \cite{C01, CHS21} for details. Throughout this paper, $C$  is a constant which may change from one line to another, and sometimes we write  $C(a,b,c,\ldots)$ to emphasize the dependence on the parameters $a,b,c,\ldots$

Let $H^s:=H^s(\mathcal{O})$ be the usual Sobolev space.
Then the domain of the operator $A$ is ${\rm Dom}(A):=H^2\cap H^1_0$, and there is a sequence of real numbers $\lambda_i \sim i^{\frac 2 d}, i\in\mathbb{N}_+$ (see \cite[Section 1]{Dirichlet_laplacian}), and an orthonormal basis $\{e_i(x)\}_{i\in\mathbb{N}_+}$ such that $Ae_i=\lambda_ie_i.$ It is known that $A$ is positive, self-adjoint, and densely defined operator on $H$, and  that $-A$ generates an  analytic semigroup $\{S(t):=e^{-tA},t\ge 0\}$ on $H$. Define the Hilbert space $\dot{H}^{\gamma}:=\mathrm{Dom}(A^{\frac{\gamma}{2}}),$ $\gamma\in\mathbb{R}$, equipped with the inner product $\langle\cdot,\cdot\rangle_{\gamma}:=\langle A^{\frac{\gamma}{2}}\cdot,A^{\frac{\gamma}{2}}\cdot\rangle$ and the norm $\|\cdot\|_{\gamma}:=\langle\cdot,\cdot\rangle^{\frac{1}{2}}_{\gamma}$.  Furthermore, $\mathcal{L}(H,U)$ and $\mathcal{L}_2(H,U)$ denote spaces of the usual bounded linear operators and Hilbert–Schmidt operators from a Hilbert space $H$ to another Hilbert space $U$, respectively. When $H=U,$ we use notations $\mathcal{L}(H)$ and $\mathcal{L}_2(H)$ for simplicity. 

It is well-known (see e.g., \cite[Lemma B.9]{Kruse14}) that there is a positive constant $C$ such that
\begin{align}
\|A^{\gamma}S(t)\|_{\mathcal{L}(H)}&\leq Ct^{-\gamma},\quad t>0,\gamma\ge 0,\label{S(t)1}\\
\|A^{-\gamma}(Id-S(t))\|_{\mathcal{L}(H)}&\leq Ct^{\gamma},\quad t>0,\gamma\in[0,1],\label{S(t)2}\\
\int_s^t\|A^{\frac{\gamma}{2}}S(t-r)u\|^2\,\mathrm{d}r&\leq C(t-s)^{1-\gamma}\|u\|^2,\quad u\in H,\;0\leq s\leq t,\;\gamma\in[0,1],\label{S(t)3}
\end{align}
and
\begin{align}\label{S(t)4}
\Big\|A^{\gamma}\int_s^tS(t-r)u\mathrm{d}r\Big\|\leq C(t-s)^{1-\gamma}\|u\|,\quad u\in H,\;0\leq s\leq t,\;\gamma\in[0,1].
\end{align}
\begin{assumption}\label{assumption1}
Let $F:L^6(\mathcal{O})\rightarrow H$ be the Nemytskii operator defined by
\begin{align*}
F(X)(x)=f(X(x)):=a_3X^3(x)+a_2X^2(x)+a_1X(x)+a_0,\,a_3<0,\,x\in \mathcal{O},\; a.e.
\end{align*}
\end{assumption}

It can be verified that there exist positive constants $L_0$ and $L_1$ such that for $X,Y\in E$,
\begin{align}
\langle X-Y,F(X)-F(Y)\rangle&\leq L_0\|X-Y\|^2,\label{one-sided Lip}\\
\|F(X)-F(Y)\|&\leq L_1(1+\|X\|^2_E+\|Y\|^2_E)\|X-Y\|.\label{polynomial Lip}
\end{align}
And for $X,\psi,\psi_1,\psi_2\in L^6(\mathcal{O}),$
\begin{align}
&\big(DF(X)(\psi)\big)(x)=(3a_3X^2(x)+2a_2X(x)+a_1)\psi(x),\quad x\in\mathcal{O},\label{DF_prop}\\
&\big(D^2F(X)(\psi_1,\psi_2)\big)(x)=(6a_3X(x)+2a_2)\psi_1(x)\psi_2(x),\quad x\in\mathcal{O}.\label{D^2F_prop}
\end{align}
Moreover, there is a positive constant $C$ such that for $X\in\dot{H}^2,$
\begin{align}\|F(X)\|^2_1&=\|A^{\frac{1}{2}}F(X)\|^2
=\int_{\mathcal{O}}\lvert(3a_3X^2(x)+2a_2X(x)+a_1)\nabla X(x)\rvert^2 \mathrm{d}x
\leq C(\|X\|^4_E+1)\|X\|^2_1,\label{F_normH^1}\\
\|F(X)\|^2_2&=\int_{\mathcal{O}}\Big|(6a_3X(x)+2a_2)|\nabla X(x)|^2+(3a_3X^2(x)+2a_2X(x)+a_1)\Delta X(x)\Big|^2\mathrm{d}x\notag\\
&\leq C(\|X\|^4_E+1)(\|X\|_2^4+1),\label{F_normH^2}
\end{align}
where we have used the Gagliardo--Nirenberg inequality  $\|u\|_{L^4(\mathcal{O})}\leq C\|\nabla u\|^{\theta}\|u\|^{1-\theta}$ for $u\in H^1$ with $\theta=\frac d4\in(0,1]$, see e.g. \cite[Eq. (5.71)]{rochnerbook}.

We make the following assumption on the stochastic process $\{W(t,\cdot)\}_{t\in [0,T]}.$
\begin{assumption}\label{assumption2}
Let $\{W(t,\cdot)\}_{t\in [0,T]}$ be a generalized $Q$-Wiener process on a filtered probability space $(\Omega,\mathcal{F},\left\{\mathcal{F}_t\right\}_{0\leq t\leq T},\mathbb{P})$, which can be represented as
$
W(t,x):=\sum_{n=1}^{\infty}Q^{\frac{1}{2}}e_n(x)\beta_n(t),
$
where $\{\beta_n(t)\}_{n\in\mathbb{N}_+}$ is a sequence of independent real valued standard Brownian motions. Assume that for some $\beta\in(0,2]$,
\begin{align}\label{noise}
\|A^{\frac{\beta-1}{2}}Q^{\frac{1}{2}}\|_{\mathcal{L}_2(H)}<\infty.
\end{align}
 In the case $\beta\leq \frac{d}{2},$ we in addition assume that $Q$ commutes with $A$.
\end{assumption}
\noindent
There are two important cases included in \eqref{noise}: the trace-class noise case (i.e. $tr(Q)<\infty$) for $\beta\in[1,2]$, and the space-time white noise case (i.e. $Q=\mathrm{Id}$) for $d=1$ and $\beta<\frac{1}{2}$. 
We remark that the condition that $Q$ commutes with $A$ for $\beta\leq \frac d2$ is used to ensure 
\begin{align}\label{sup_noise}\sup_{t\in[0,T]}\mathbb{E}\Big[\Big\|\int_0^tS(t-s)\mathrm{d}W(s)\Big\|^p_E\Big]<\infty
\end{align} for $p>\frac{2d}{\beta\wedge 1},$ while for the case of $\beta>\frac d2$, \eqref{sup_noise} 
can be  proved directly by the Sobolev embedding $\dot{H}^{\beta}\hookrightarrow E$; see also \cite[Lemma 2]{CHS21}.
\begin{assumption}\label{assumption3}
The initial datum satisfies $\mathbb{E}\big[\exp\{\|X_0\|_{\dot{H}^{\beta}\cap E}\}\big]<\infty$,
%$X_0\in L^{s_0}(\Omega;\dot{H}^{\beta}\cap E)$ for some $s_0\ge 2$ sufficiently large,  
where $\beta$ is given in Assumption \ref{assumption2}. In addition, $X_0$ is $\mathcal{F}_0/\mathcal{B}(H\cap E)$-measurable.
\end{assumption}

With the above assumptions, we can get the existence, uniqueness, and regularity estimates of the mild solution of \eqref{spde}. The proofs can be found in e.g. \cite[Proposition 6.2.2]{C01}, \cite[Theorem 2.1]{caimeng21}, and  \cite[Theorem 2.6]{WXJ20}.

\begin{theorem}\label{eau}
Under Assumptions \ref{assumption1}-\ref{assumption3}, the stochastic Allen--Cahn equation \eqref{spde} has a unique mild solution given by
\begin{align*}
X(t)=S(t)X_0+\int_0^tS(t-s)F(X(s))\,\mathrm{d}s+\int_0^tS(t-s)\,\mathrm{d}W(s)\quad a.s.
\end{align*}
For $p\ge 2$, we have \begin{align}
&\sup_{t\in[0,T]}\|X(t)\|_{L^p(\Omega;E)}+\sup_{t\in[0,T]}\|X(t)\|_{L^p(\Omega;\dot{H}^{\beta})}\leq C_1,\label{supest}\\
&\|X(t)-X(s)\|_{L^p(\Omega;H)}\leq C_2(t-s)^{\frac{\beta\wedge 1}{2}},\quad  0\leq s<t\leq T,\label{Holderexact}
\end{align}
where constants $C_1,C_2>0$ depend on $X_0,p,T,\|A^{\frac{\beta-1}{2}}Q^{\frac12}\|_{\mathcal{L}_2(H)}$.
\end{theorem}

\section{Adaptive schemes}\label{sec2}
In this section, we first introduce the adaptive time-stepping scheme, and assumptions to ensure that the final time $T$ can be attained in finite many steps. Then we present the coupled scheme for the practical use.
 Finally, we show the strong convergence orders of these numerical  schemes.

\subsection{Schemes}

Introduce the adaptive timestep function $\tau:H\to\mathbb{R}_+.$
We consider the following adaptive time-stepping scheme, whose spatial discretization is based on the spectral Galerkin method, and temporal direction is the adaptive exponential integrator:
\begin{align}\tag{AE}\label{scheme}
X^N_{t_{m+1}}=S^N(\tau_m)X^N_{t_m}+S^N(\tau_m)F^N(X^N_{t_m})\tau_m+S^N(\tau_m)P^N\Delta W_m,\quad X^N_0=P^NX_0,
\end{align}
where $\tau_m:=\tau(X^N_{t_m}),\,t_0=0,\,t_{m+1}=t_m+\tau_m,\,S^N(t):=P^NS(t)=e^{-tA^N}$ with $P^N:H\to H_N\,(:=span\{e_1,\ldots,e_N\})$ being the spectral projection operator and  $A^N:=P^NA$,  $F^N:=P^NF,$ and the increment 
$
\Delta W_m:=W(t_{m+1})-W(t_m).
$
If the existing time span is longer than $T$ after adding the last timestep, then we take a smaller timestep such that the existing time span just attains $T$ after adding it. Namely, letting $M_{T}$ be the number of timesteps for a given timestep function
$\tau,$ if $t_{M_T-1}+\tau_{M_T-1}>T,$ then we enforce the last timestep $\tau_{M_T-1}:=T-t_{M_T-1}.$ 
In the sequel, we will give some assumptions on the timestep function so that the numerical solution can attain $T$ with finite many timesteps.

The continuous version of \eqref{scheme} is given by
\begin{align}\label{continuous version}
X^N_t=S^N(t)P^NX_0+\int_0^tS^N(t-t_{m_s})F^N(X^N_{t_{m_s}})\,\mathrm{d}s+\int_0^tS^N(t-t_{m_s})P^N\,\mathrm{d}W(s),
\end{align}
where 
%$t_{m_s}:=\max\{t_m:t_m\leq s\}$, i.e., the maximum grid point before time $s$. Denote 
$m_s:=\max\{m:t_m\leq s\}.$ %Then $t_{m_s}=t_{m_s}$. 

In order to bound the number of timesteps, we give the following assumption on the adaptive timestep function with the uniform lower bound. 
\begin{assumption}\label{assumption 4}
The adaptive timestep function $\tau:H\to\mathbb{R}_+$ is continuous and satisfies that for $X(\omega)\in H,$
\begin{align}
\tau(X(\omega))^{1-\frac{d}{4}}\|F(X(\omega))\|&\leq L_2,\quad a.s.,\label{adaptive1}\\
\tau(X(\omega))&\ge \tau_{min},\quad a.s.\label{adaptive2}
\end{align}
with positive constants $L_2$ and $\tau_{min}$ independent of $\omega$. 
\end{assumption} 

Under the assumption \eqref{adaptive2},  we have $M_T\leq \frac{T}{\tau_{min}}<\infty,$ a.s. That is to say, $T$ is a.s. attainable in finite many timesteps. The power $1-\frac{d}{4}$ in \eqref{adaptive1} is for technical reason to get the \textit{a priori} estimate  in $E$-norm of the solution of \eqref{scheme}. 
Examples for adaptive timestep functions that satisty Assumption \ref{assumption 4} are given in Section \ref{sec3}.

\begin{remark}\label{rmkbeta=1}
If the expected supremum of the $p$th moment of the numerical solution is finite, i.e., $\mathbb{E}\Big[\sup_{0\leq t\leq T}\|X^N_t\|^p\Big]<\infty$ for some large $p\ge2,$ then the bounds of adaptive timestep function in Assumption \ref{assumption 4} can be weaken to the adaptive ones:
\begin{align}
\tau(X(\omega))^{1-\frac{d}{4}}\|F(X(\omega))\|&\leq L_2(\omega),\quad a.s.,\label{adaptive4}\\
\tau(X(\omega))&\ge (\zeta_1\|X(\omega)\|^{q_0}+\zeta_2)^{-1},\quad a.s.\label{adaptive5}
\end{align}
with positive constants $\zeta_1,\zeta_2$, and $q_0$ independent of $\omega$, and $q_0\le p,$ $\|L_2(\cdot)\|_{L^{2p}(\Omega)}<\infty$.
In this setting, $T$ is still a.s. attainable, i.e., 
$$
\mathbb{E}[M_{T}]\leq T\mathbb{E}\Big[\sup_{0\leq t_m\leq T}\frac{1}{\tau(X^N_{t_m})}\Big]
\leq T\mathbb{E}\big[\sup_{0\leq t_m\leq T}(\zeta_1\|X^N_{t_m}\|^{q_0}+\zeta_2)\big]<\infty.
$$

Notice that for the trace-class noise case (i.e. $tr(Q)<\infty$ with $\beta\in[1,2]$ in Assumption \ref{assumption2}), 
 under the assumption
\begin{align}
\langle X(\omega),F(X(\omega))\rangle+\frac{1}{2}\tau(X(\omega))\|F(X(\omega))\|^2&\leq \bar{L}_2\|X(\omega)\|^2+L_3,\quad a.s.\label{adaptive3}
\end{align}
with positive constants $\bar{L}_2$ and $L_3$ independent of $\omega,$ following the approach of \cite[Theorem $1$]{Fang20} and combining with the contractivity of the semigroup $\{S(t),t\ge0\}$ in $H$,
 we can get the finiteness of the expected supremum of the $p$th moment of the numerical solution. However, for  $\beta\in(0,1),$ we have not gotten  the finiteness of the expected supremum of the $p$th moment of the numerical solution. Hence, the main result in this paper is still hold under assumptions \eqref{adaptive4}-\eqref{adaptive3} when $\beta\in[1,2]$. 
\end{remark}

In practice, instead of verifying  whether a timestep function satisfies the lower bound in \eqref{adaptive2} or \eqref{adaptive5}, people usually introduce a backstop scheme with a uniform timestep and couple it with \eqref{scheme} to ensure that a simulation over the interval $[0,T]$ can be completed in a finite number of timesteps (see \cite{Lord18} for the case of SODEs).
More precisely,  if $\tau_m<\tau_{min}$ 
at time $t_m$, then we apply a single step of some convergent scheme $\Psi:H_N\times \mathbb{R}_+\times H\to H_N,$ which is called the backstop scheme over a timestep of length $\tau_{min}$ instead, i.e.,  
\begin{align}\label{expression_cou}
X^{N,C}_{t_{m+1}}=\Phi(X^{N,C}_{t_{m}},\tau_m,\Delta W_m)\chi_{\{\tau_m\ge \tau_{min}\}}+\Psi(X^{N,C}_{t_{m}},\tau_{min},\Delta W_m)\chi_{\{\tau_m<\tau_{min}\}},
\end{align}
 where the map $\Phi:H_N\times \mathbb{R}_+\times H\to H_N,\; (x,h,y)\mapsto S^N(h)(x+F^N(x)h+y)$ denotes the scheme \eqref{scheme}.
 
The backstop scheme is usually chosen to be an  explicit and  convergent scheme, for instance, the tamed exponential integrator (see \cite{WXJ20}), the nonlinearity-truncated exponential integrator and the linear-implicit nonlinearity-truncated scheme (see \cite{BJ19}). In the following, we take the backstop scheme $\Psi$  as the tamed exponential integrator: 
\begin{align}\label{Psi1}
\Psi(x,h,y)=S^N(h)\Big(x+\frac{F^N(x)h}{1+\|F^N(x)\|h}+y\Big),
\end{align}
and give estimates of the corresponding  coupled scheme  
\begin{align}\tag{CAU 1}\label{ATEU}
X^{N,(1)}_{t_{m+1}}=&\;S^N(\tau_m)\Big[X^{N,(1)}_{t_m}+F^N(X^{N,(1)}_{t_m})\tau_m+\Delta W_m\Big]\chi_{\{\tau_m\ge \tau_{min}\}}\nonumber\\
&+S^N(\tau_{min})\Bigg[X^{N,(1)}_{t_m}+\frac{F^N(X^{N,(1)}_{t_m})\tau_{min}}{1+\|F^N(X^{N,(1)}_{t_m})\|\tau_{min}}+\Delta W_m\Bigg]\chi_{\{\tau_m<\tau_{min}\}}.\notag
\end{align}
Moreover, the adaptive timestep function that satisfies \eqref{adaptive1} can be chosen as, e.g. 
$\tau_m=\tau(X^N_{t_m})=\big(\frac{L_2}{\|F(X^N_{t_m})\|+c}\big)^{\frac{4}{4-d}}$ with some $c> 0$.
We can write the continuous versions of \eqref{ATEU} into the compact integral form by defining a new timestep function denoted by $\tau_{new}$ with
\begin{equation*}
\tau_{new}:=\left\{
\begin{array}{ll}
\tau_{m},\quad &\text{for  }\tau_m\ge \tau_{min},\\
\tau_{min},\quad &\text{for  }\tau_m<\tau_{min}.
\end{array}
\right.
\end{equation*}
Then $t_{m+1}=t_m+\tau_{new},$ and the continuous versions of \eqref{ATEU} is
\begin{align}
X^{N,(1)}_t&=S^N(t)P^NX_0+\int_0^tS^N(t-t_{m_s})\Big(F^N(X^{N,(1)}_{t_{m_s}})\chi_{\{\tau_{m_s}\ge \tau_{min}\}}\notag\\
&\quad+\frac{F^N(X^{N,(1)}_{t_{m_s}})}{1+\|F^N(X^{N,(1)}_{t_{m_s}})\|\tau_{min}}\chi_{\{\tau_{m_s}<\tau_{min}\}}\Big)\mathrm{d}s
 +\int_0^tS^N(t-t_{m_s})P^N\mathrm{d}W(s).\label{conti_ATEU1}
\end{align}

\subsection{Main result}

In this subsection, we give strong convergence orders for the scheme \eqref{scheme} as well as the coupled scheme \eqref{ATEU}, 
whose proofs are postponed to Section \ref{secproof} and Appendix \ref{appendix2}, respectively.
 
Since the timestep function $\tau$ is determined by the numerical solution, we need to make a modification when considering the convergence of  adaptive time-stepping schemes (see \cite{Fang20}).
Namely, for a given timestep function $\tau$ which satisfies Assumption \ref{assumption 4}, we introduce the refined timestep function $\tau^{\delta}$ controlled by a scalar parameter $\delta\in(0,1)$  and consider the convergence when $\delta\to0$ as well as the order with respect to $\delta.$
\begin{assumption}\label{assumption6}
The refined timestep function $\tau^{\delta}$ satisfies that for $X(\omega)\in H,$
\begin{align*}
\delta\min\{T,\tau(X(\omega))\}\leq \tau^{\delta}(X(\omega))\leq \min\{T\delta,\tau(X(\omega))\},\quad a.s.
\end{align*}
\end{assumption}

Examples of $\tau$ and $\tau^{\delta}$ are given in Section \ref{sec3}. We also remark that in this setting, the lower bound in \eqref{adaptive2} is defined as $\tau^{\delta}_{min}:=\delta\tau_{min}.$ With this assumption in hand, 
in the following, we present strong convergence orders of schemes \eqref{scheme} and  \eqref{ATEU} with the timestep function $\tau^{\delta}.$ %whose proofs are given in Section \ref{secproof}. 
Before that, we put an additional  assumption on the initial datum, which is used to get the \textit{a priori} estimates of  numerical solutions in $E$-norm, see Section \ref{sec2.1} for details.
\begin{assumption}\label{assumpP^NX_0}
The initial datum of \eqref{spde} satisfies that $\sup_{N\in \mathbb{N}_+}\mathbb{E}\big[\exp\{\|P^NX_0\|_E\}\big]<\infty$.
%for sufficiently large $p_0\ge 2$, $\sup_{N\in \mathbb{N}_+}\|P^NX_0\|_{L^{p_0}(\Omega;E)}<\infty.$
\end{assumption}

Based on the Sobolev embedding theorem, Assumption \ref{assumpP^NX_0} is fulfilled if 
$\sup_{N\in\mathbb{N}_+}\mathbb{E}\big[\exp\{\|X_0\|_{\dot{H}^{\beta}}\}\big]<\infty$ for $\beta>\frac{d}{2}.$ 
\begin{theorem}\label{thm3}
Under Assumptions 
\ref{assumption1}-\ref{assumpP^NX_0}, for $p\ge 2,$
\begin{align*}
\sup_{0\leq t\leq T}\|X(t)-\mathbb{X}^N_t\|_{L^p(\Omega;H)}&\leq C({\lambda_N}^{-\frac{\beta}{2}}+\delta^{\frac{\beta}{2}}),
\end{align*}
where $\mathbb{X}^N_t$ is the numerical solution $X^N_t$ (or $X^{N,(1)}_t$) of the scheme \eqref{scheme} (or \eqref{ATEU}), and
$C>0$ depends on $p,T,\tau_{min},L_0,L_1,L_2,$ $X_0$, and  $\|A^{\frac{\beta-1}{2}}Q^{\frac{1}{2}}\|_{\mathcal{L}_2(H)}$.
\end{theorem}

Noting that $\lambda_N\sim N^{-\frac 2d},$ the convergence result in Theorem \ref{thm3} can be rewritten as
$
\sup_{0\leq t\leq T}\|X(t)-\mathbb{X}^N_t\|_{L^p(\Omega;H)}\leq C({N}^{-\frac{\beta}{d}}+\delta^{\frac{\beta}{2}}).
$
Then we say that the spatial convergence order is $\frac{\beta}{d}$ with respect to $N.$

\begin{remark}\label{remark_couple}
As stated in Remark \ref{rmkbeta=1}, for the trace-class noise case (i.e., $\beta\in[1,2]$), 
the bounds of the adaptive timestep function in Assumption \ref{assumption 4} can be weaken to 
\eqref{adaptive4} and \eqref{adaptive5}. Similarly to the definition of \eqref{expression_cou}, when the critical parameter for the adaptive timestep size is $\big(\zeta_1\big\|\widetilde{X^{N,C}_{t_m}}\big\|^{q_0}+\zeta_2\big)^{-1}$, the coupled scheme can also be defined as
\begin{align}\label{expression_cou2}
\widetilde{X^{N,C}_{t_{m+1}}}=\Phi(\widetilde{X^{N,C}_{t_m}},\tau_m,\Delta W_m)\chi_{\big\{\tau_m\ge \big(\zeta_1\big\|\widetilde{X^{N,C}_{t_m}}\big\|^{q_0}+\zeta_2\big)^{-1}\big\}}+\Psi(\widetilde{X^{N,C}_{t_m}},\tau_m,\Delta W_m)\chi_{\big\{\tau_m< \big(\zeta_1\big\|\widetilde{X^{N,C}_{t_m}}\big\|^{q_0}+\zeta_2\big)^{-1}\big\}}.
\end{align} 
If we still choose $\Psi$ to be \eqref{Psi1}, and denote the solution of the continuous version of the corresponding coupled scheme  
 by $\widetilde{X^{N,(1)}_{t}}$,   then the similar proof as that of Theorem \ref{thm3} yields that: under Assumptions \ref{assumption1}-\ref{assumption3}, Assumptions \ref{assumption6}-\ref{assumpP^NX_0} and Eq. \eqref{adaptive4}-\eqref{adaptive3}, for $p\ge 2,$
\begin{align*}
\sup_{0\leq t\leq T}\|X(t)-\widetilde{X^{N,(1)}_t}\|_{L^p(\Omega;H)}\leq C(\lambda_N^{-\frac{\beta}{2}}+\delta^{\frac{\beta}{2}})\leq C(\lambda_N^{-\frac{\beta}{2}}+(\mathbb{E}[M_{T,\delta}])^{-\frac{\beta}{2}}),
\end{align*}
where $M_{T,\delta}$ is the number of timesteps for the given timestep function
$\tau^{\delta},$ and
$C>0$ depends on $p,T,\zeta_1,\zeta_2,L_i(i=0,\ldots,3),\bar{L}_2,$ $X_0,$ and $\|A^{\frac{\beta-1}{2}}Q^{\frac{1}{2}}\|_{\mathcal{L}_2(H)}$.
\end{remark}
\begin{remark}
We remark that it is interesting to investigate if one can obtain the strong convergence order for the coupled scheme directly 
from some    error estimates  of schemes $\Phi$ and $\Psi$.
For this problem, a fundamental  convergence theorem that characterizes the relation between the local error and the global  error might be helpful.
 For the study of such a theorem, we refer to \cite{Milstein} for the case of SODEs with either the globally or locally Lipschitz  drift, and to \cite{CCC} for the case of SPDEs with the globally Lipschitz drift.  
 However,  for SPDEs with  the non-Lipschitz drift, there has been no work on such theorem. 
 We leave this as the future work. 
\end{remark}

\section{Estimates of numerical solutions}\label{sec2.1}
In this section, we analyze the \textit{a priori} estimates of  numerical solutions \eqref{scheme} and \eqref{ATEU} in $E$-norm and $\dot{H}^{\beta}$-norm, respectively. Proofs of all the results in this section are given in Appendix \ref{appendix1} for readers' convenience. 

We give the following lemma on the properties of the semigroup $\{S(t),t\ge 0\}$, which is important in the \textit{a priori} estimates of numerical solutions. 

\begin{lemma}\label{lemma1.4} 
We have\\
$(\romannumeral1)$ $\;\sup_{N\in\mathbb{N}_+}\|A^{\rho}P^NS(t)u\|_E\leq C(\rho,d)(t^{-\rho}+t^{-\rho-\frac{d}{4}})\|u\|,$\quad for $\rho\ge0,\;t>0,\;u\in H;$\\
$(\romannumeral2)$
$\sup_{N\in\mathbb{N}_+}\|P^NS(t)u\|_E\leq C(\rho,d)t^{\frac{2\rho-d}{4}}\|u\|_{\rho},$\quad for $\rho\in[0,\frac{d}{2}),\; t>0,\;u\in\dot{H}^{\rho}.$
\end{lemma}

Let $v^N,z^N:[0,T]\to H_N$ satisfy the perturbed differential equation \begin{align}
\begin{cases}\label{perturbed} 
\mathrm{d}v^N_t=(-A^Nv^N_t+F^N(v^N_t+z^N_t))\mathrm{d}t,\quad t\in(0,T],\\
v^N_0=0.
\end{cases}
\end{align}
In what follows, we aim
to show that the $E$-norm of the solution $v^N_t$ of \eqref{perturbed} can be controlled by the $L^{\eta}([0,t];E)$-norm of the perturbation $z^N_{\cdot}$ with some $\eta>0$, which plays a crucial role in deriving  the moment bounds in $E$-norm for the scheme \eqref{scheme}. 
\begin{lemma}\label{purt_lemma}
The solution of \eqref{perturbed} satisfies 
\begin{align*}
\|v^N_t\|_E\leq C\Big[1+\int_0^t\Big((t-s)^{-\frac d4}\|z^N_s\|^3_E+\|z^N_s\|^{72}_E\Big)\mathrm{d}s\Big],\quad t\in(0,T]
\end{align*}
with $C>0.$
\end{lemma}

With these preparations, we can establish the \textit{a priori} estimates of the numerical solution of \eqref{scheme} in $E$-norm  by the standard bootstrap argument; see \cite{WXJ20} for the description of this approach for the nonlinearity-tamed scheme with the uniform timestep. %Recall that parameters $s_0$ and $p_0$ are given in Assumptions \ref{assumption3} and \ref{assumpP^NX_0}, respectively.

\begin{prop}\label{thm2}
Under conditions in Theorem \ref{thm3}, 
we have for $p\ge 2$, 
\begin{align}\label{stability}
\sup_{N\in\mathbb{N}_+}\sup_{0\leq t\leq T}\|X^N_{t}\|_{L^p(\Omega;E)}\leq C,
\end{align}
where $X^N_{t}$ is the numerical solution of \eqref{scheme} with timestep function $\tau^{\delta}$ for some $\delta\in(0,1)$, and $C>0$ depends on $p,T,L_1,L_2,\tau_{min},$ $X_0,$ and $\|A^{\frac{\beta-1}{2}}Q^{\frac{1}{2}}\|_{\mathcal{L}_2(H)}$.
\end{prop}

With the \textit{a priori} estimate of the numerical solution of \eqref{scheme} in $E$-norm in hand, we can get the following 
\textit{a priori} estimate in $\dot{H}^{\beta}$-norm directly by means of the mild form of the solution. A  standard argument gives the H\"older continuity of the numerical solution; see also \cite[Lemma 4.3]{BJ19}. 

\begin{prop}\label{coro}
Under conditions in Proposition \ref{thm2},  for $p\ge 2$, $\beta\in(0,2],$ and $\gamma\in(0,\beta],$ we have \begin{align}&\sup_{N\in \mathbb{N}_+}\sup_{0\leq t\leq T}\|X^N_t\|_{L^p(\Omega;\dot{H}^{\beta})}\leq C_1,\label{coroeq1}\\
&\sup_{N\in \mathbb{N}_+}\|X^{N}_t-X^N_s\|_{L^p(\Omega;\dot{H}^{\gamma})}\leq C_2(t-s)^{\frac{(\beta-\gamma)\wedge 1}{2}},\quad 0\leq s<t\leq T,\label{coroeq11}
\end{align}
where constants $C_1,C_2>0$ depend on $p, T, L_1,L_2, \tau_{min},X_0,$ and $\|A^{\frac{\beta-1}{2}}Q^{\frac{1}{2}}\|_{\mathcal{L}_2(H)}.$
%\|X_0\|_{L^{q}(\Omega;\dot{H}^{\beta})},$ $\sup\limits_{N\in \mathbb{N}_+}\|P^NX_0\|_{L^{q}(\Omega;E)}$ with $q$ being defined in Proposition \ref{thm2}.
\end{prop}

With the above two propositions, combining regularity estimates of the tamed exponential integrator, which can be proved similarly as Propositions \ref{thm2}-\ref{coro} and  \cite{WXJ20},  we get the following the  \textit{a priori} estimates in $\dot{H}^{\beta}$ and the H\"older continuity for the numerical solution of \eqref{ATEU}.
\begin{corollary}\label{coro2}
Under conditions in Proposition \ref{thm2}, for $p\ge 2$, $\beta\in(0,2],$ and $\gamma\in(0,\beta],$ we have
\begin{align*}
&\sup_{N\in \mathbb{N}_+}\sup_{0\leq t\leq T}\|X^{N,(1)}_t\|_{L^p(\Omega;E\cap\dot{H}^{\beta})}\leq C_1,\\
&\sup_{N\in \mathbb{N}_+}\|X^{N,(1)}_t-X^{N,(1)}_s\|_{L^p(\Omega;\dot{H}^{\gamma})}\leq C_2(t-s)^{\frac{(\beta-\gamma)\wedge 1}{2}},\quad 0\leq s<t\leq T,
\end{align*}
where constants $C_1,C_2>0$ depend on $p,T,L_1,L_2,\tau_{min}, X_0, $ and $\|A^{\frac{\beta-1}{2}}Q^{\frac{1}{2}}\|_{\mathcal{L}_2(H)}$.
%\|X_0\|_{L^{q}(\Omega;\dot{H}^{\beta})}$, $\sup\limits_{N\in \mathbb{N}_+}\|P^NX_0\|_{L^{q}(\Omega;E)}$ with $q$ being defined in Proposition \ref{thm2}.
\end{corollary}

\section{Proof of the main result}\label{secproof}
In this section, 
based on the \textit{a priori} estimates of numerical solutions presented in Section \ref{sec2.1}, we show the proof of the convergence order of the scheme \eqref{scheme} in Theorem \ref{thm3}, 
and leave that of the coupled scheme \eqref{ATEU} to Appendix \ref{appendix2}.
\begin{proof}[Proof of Theorem \ref{thm3}] 
%\textbf{Proof of the convergence of scheme \eqref{scheme}.}
By introducing an auxiliary  process,
\begin{align*}
Y^N_t=S^N(t)P^NX_0+\int_0^tS^N(t-s)F^N(X(s))\mathrm{d}s+\int_0^tS^N(t-t_{m_s})P^N\,\mathrm{d}W(s),
\end{align*}
the error can be divided into the following terms,
\begin{align*}
\|X(t)-X^{N}_t\|_{L^p(\Omega;H)}
\leq \|X(t)-P^NX(t)\|_{L^p(\Omega;H)}+\|P^NX(t)-Y^N_t\|_{L^p(\Omega;H)}+\|Y^N_t-X^{N}_t\|_{L^p(\Omega;H)}.
\end{align*}
The term $\|X(t)-P^NX(t)\|_{L^p(\Omega;H)}$ can be estimated as 
\begin{align}\label{X-P^NX}
\|X(t)-P^NX(t)\|_{L^p(\Omega;H)}=\|A^{-\frac{\beta}{2}}(\mathrm{Id}-P^N)A^{\frac{\beta}{2}}X(t)\|_{L^p(\Omega;H)}\leq C\lambda_N^{-\frac{\beta}{2}}\|X(t)\|_{L^p(\Omega;\dot{H}^{\beta})}.
\end{align}

For the term $\|P^NX(t)-Y^N_t\|_{L^p(\Omega;H)},$  using the Burkholder--Davis--Gundy  inequality (see e.g. \cite[Proposition 2.12]{Kruse14}), \eqref{S(t)2}-\eqref{S(t)3},  the H\"older inequality, and Assumption \ref{assumption6} gives that
\begin{align}\label{eq4}
&\|P^NX(t)-Y^N_t\|^p_{L^p(\Omega;H)}=\Big\|\int_0^tS^N(t-s)(\mathrm{Id}-S^N(s-t_{m_s}))\mathrm{d}W(s)\Big\|^p_{L^p(\Omega;H)}\nonumber\\
&\leq C(p)\mathbb{E}\Big[\Big(\int_0^t\|S^N(t-s)(\mathrm{Id}-S^N(s-t_{m_s}))Q^{\frac{1}{2}}\|^2_{\mathcal{L}_2(H)}\mathrm{d}s\Big)^{\frac{p}{2}}\Big]\nonumber\\
&= C(p)\mathbb{E}\Big[\Big(\int_0^t\|A^{-\frac{\beta}{2}}(\mathrm{Id}-S^N(s-t_{m_s}))A^{\frac{1}{2}}S^N(t-s)A^{\frac{\beta-1}{2}}Q^{\frac{1}{2}}\|^2_{\mathcal{L}_2(H)}\mathrm{d}s\Big)^{\frac{p}{2}}\Big]\nonumber\\
&\leq C(p)\mathbb{E}\Big[\Big(\int_0^t\sum_{j=1}^{\infty}(s-t_{m_s})^{\frac{\beta}{2}}\|A^{\frac{1}{2}}S^N(t-s)A^{\frac{\beta-1}{2}}Q^{\frac{1}{2}}e_j\|^2\mathrm{d}s\Big)^{\frac p2}\Big]\nonumber\\
&\leq C(p,T)\|A^{\frac{\beta-1}{2}}Q^{\frac{1}{2}}\|^p_{\mathcal{L}_2(H)}\delta^{\frac{\beta p}{2}}\leq C(p,T,\|A^{\frac{\beta-1}{2}}Q^{\frac{1}{2}}\|_{\mathcal{L}_2(H)})\delta^{\frac{\beta p}{2}}.
\end{align}

For the estimate of the term $\|Y^N_t-X^{N}_t\|_{L^p(\Omega;H)},$ combining the differential form
\begin{align*}
\mathrm{d}(X^{N}_t-Y^N_t)&=-A^N(X^{N}_t-Y^N_t)\mathrm{d}t+\Big(S^N(t-t_{m_t})F^N(X^{N}_{t_{m_t}})-F^N(X(t))\Big)\mathrm{d}t,
\end{align*}
and applying the Taylor formula yield 
\begin{align}\label{add_X-Y}
\|X^{N}_t-Y^N_t\|^2=&\;2\int_0^t\langle X^{N}_s-Y^N_s,-A^N(X^{N}_s-Y^N_s)\rangle\mathrm{d}s\notag\\
&+2\int_0^t\Big\langle X^{N}_s-Y^N_s,S^N(s-t_{m_s})F^N(X^{N}_{t_{m_s}})-F^N(X(s))\Big\rangle\mathrm{d}s\notag\\
=&\;2\int_0^t\langle X^{N}_s-Y^N_s,-A^N(X^{N}_s-Y^N_s)\rangle\mathrm{d}s\notag\\
& +2\int_0^t\langle X^{N}_s-Y^N_s,S^N(s-t_{m_s})F^N(X^{N}_{t_{m_s}})-F^N(X^N_s) \rangle\mathrm{d}s\notag\\
&+2\int_0^t\langle X^{N}_s-Y^N_s,F^N(X^{N}_{s})-F^N(Y^N_s)\rangle \mathrm{d}s
+2\int_0^t\langle X^{N}_s-Y^N_s,F^N(Y^N_s)-F^N(X(s))\rangle \mathrm{d}s\notag\\
\leq&\, (2L_0+1)\int_0^t\|X^{N}_s-Y^N_s\|^2\mathrm{d}s+\int_0^t\|F^N(Y^N_s)-F^N(X(s))\|^2\mathrm{d}s+2J_{1,1}(t)+2J_{1,2}(t),
\end{align}
where we have used the condition \eqref{one-sided Lip} and the Young inequality, 
and $$J_{1,1}(t):=\int_0^t\langle X^{N}_s-Y^N_s,(S^N(s-t_{m_s})-\mathrm{Id})F^N(X^{N}_{t_{m_s}})\rangle\mathrm{d}s, \quad J_{1,2}(t):=\int_0^t\langle X^N_s-Y^N_s,F^N(X^N_{t_{m_s}})-F^N(X^N_s)\rangle\mathrm{d}s.$$

\textit{Estimate of $J_{1,1}.$}
For the case $\beta=1$ or $\beta=2$, it follows from \eqref{S(t)2}, \eqref{F_normH^1}-\eqref{F_normH^2} and the Young inequality that
\begin{align*}
J_{1,1}(t)&\leq \frac12\int_0^t\|X^{N}_s-Y^N_s\|^2\mathrm{d}s+\frac12\int_0^t\|A^{-\frac{\beta}{2}}(S^N(s-t_{m_s})-\mathrm{Id})A^{\frac{\beta}{2}}F(X^{N}_{t_{m_s}})\|^2\mathrm{d}s\\
&\leq \frac12\int_0^t\|X^{N}_s-Y^N_s\|^2\mathrm{d}s+C\int_0^t(s-t_{m_s})^{\beta}(1+\|X^{N}_{t_{m_s}}\|^4_E)(\|X^N_{t_{m_s}}\|_{\beta}^4+1).
\end{align*}
This, combining Assumption \ref{assumption6},  Propositions \ref{thm2} and \ref{coro} leads to 
\begin{align*}
\|J_{1,1}(t)\|_{L^p(\Omega)}\leq \frac12\int_0^t\|X^{N}_s-Y^N_s\|^2_{L^{2p}(\Omega;H)}\mathrm{d}s+
C\delta^{\beta},\quad t\in[0,T].
\end{align*}

For the case $\beta\in(0,1)\cup(1,2),$ we need to further split the term $J_{1,1}$ into three parts, which are denoted by $J_{1,1}^i,i=1,2,3,$ based on 
\begin{align}\label{splittingX-Y}
X^{N}_s-Y^N_s&=\int_0^s\big(S^N(s-t_{m_r})-S^N(s-r)\big)F^N(X^{N}_{t_{m_r}})\mathrm{d}r+\int_0^sS^N(s-r)\big(F^N(X^N_{t_{m_r}})-F^N(X^N_r)\big)\mathrm{d}r\notag\\
&\quad+\int_0^sS^N(s-r)\big(F^N(X^N_r)-F^N(X(r))\big)\mathrm{d}r.
\end{align}
Namely,
\begin{align*}
J^1_{1,1}(t)&
:= \int_0^t\Big\langle \int_0^s S^N(s-r)(S^N(r-t_{m_r})-\mathrm{Id})F^N(X^{N}_{t_{m_r}})\mathrm{d}r,\;(S^N(s-t_{m_s})-\mathrm{Id})F^N(X^{N}_{t_{m_s}})\Big\rangle \mathrm{d}s,\\
J^2_{1,1}(t)&:=\int_0^t\Big\langle \int_0^s S^N(s-r)\big(F^N(X^{N}_{t_{m_r}})-F^N(X^{N}_r)\big)\mathrm{d}r,\;(S^N(s-t_{m_s})-\mathrm{Id})F^N(X^{N}_{t_{m_s}})\Big\rangle \mathrm{d}s,\\
J^3_{1,1}(t)&:=\int_0^t\Big\langle \int_0^s S^N(s-r)\big(F^N(X^{N}_{r})-F^N(X(r))\big)\mathrm{d}r,\;(S^N(s-t_{m_s})-\mathrm{Id})F^N(X^{N}_{t_{m_s}})\Big\rangle \mathrm{d}s.
\end{align*}
For the term $J^1_{1,1},$ by the Young inequality, we obtain that for $\beta\in(0,1)$,
\begin{align*}
J_{1,1}^1(t)&=\int_0^t\Big\langle \int_0^s A^{\beta}S^N(s-r)A^{-\frac{\beta}{2}}(S^N(r-t_{m_r})-\mathrm{Id})F^N(X^{N}_{t_{m_r}})\mathrm{d}r, \; A^{-\frac{\beta}{2}}(S^N(s-t_{m_s})-\mathrm{Id})F^N(X^{N}_{t_{m_s}})\Big\rangle \mathrm{d}s\\
&\leq C\int_0^t\int_0^s(s-r)^{-\beta}\Big[(r-t_{m_r})^{\beta}\|F(X^{N}_{t_{m_r}})\|^2+(s-t_{m_s})^{\beta}\|F(X^{N}_{t_{m_s}})\|^2\Big]\mathrm{d}r\mathrm{d}s.
\end{align*}
When $\beta\in(1,2),$ applying \eqref{S(t)1}-\eqref{S(t)2} and  \eqref{F_normH^1}, the term $J_{1,1}^1$ can be estimated as
\begin{align*}
J_{1,1}^1(t)&=\int_0^t\Big\langle \int_0^s A^{\beta-1}S^N(s-r)A^{-\frac{\beta}{2}}(S^N(r-t_{m_r})-\mathrm{Id})A^{\frac{1}{2}}F^N(X^{N}_{t_{m_r}})\mathrm{d}r, \\
&~~\quad A^{-\frac{\beta}{2}}(S^N(s-t_{m_s})-\mathrm{Id})A^{\frac{1}{2}}F^N(X^{N}_{t_{m_s}})\Big\rangle \mathrm{d}s\\
&\leq C\int_0^t\int_0^s(s-r)^{1-\beta}\Big[(r-t_{m_r})^{\beta}(1+\|X^{N}_{t_{m_r}}\|^4_E)\|X^{N}_{t_{m_r}}\|^2_1+(s-t_{m_s})^{\beta}(1+\|X^{N}_{t_{m_s}}\|^4_E)\|X^{N}_{t_{m_s}}\|^2_1\Big]\mathrm{d}r\mathrm{d}s.
\end{align*}
Hence, we derive from Propositions \ref{thm2} and \ref{coro} that for $\beta\in(0,1)\cup(1,2),$
\begin{align*}
\|J^1_{1,1}(t)\|_{L^p(\Omega)}\leq C\delta^{\beta},\quad t\in[0,T].
\end{align*}
For the term $J^2_{1,1},$ when $\beta\in(0,1),$
\begin{align*}
J_{1,1}^2(t)&=\int_0^t\Big\langle \int_0^s A^{\frac{\beta}{2}}S^N(s-r)(F^N(X^{N}_{t_{m_r}})-F^N(X^{N}_r))\mathrm{d}r,  \;A^{-\frac{\beta}{2}}(S^N(s-t_{m_s})-\mathrm{Id})F^N(X^{N}_{t_{m_s}})\Big\rangle \mathrm{d}s\\
&\leq C\int_0^t\int_0^s(s-r)^{-\frac{\beta}{2}}\Big[(1+\|X^{N}_{t_{m_r}}\|^4_E+\|X^{N}_{r}\|^4_E)\|X^{N}_{r}-X^{N}_{t_{m_r}}\|^2+(s-t_{m_s})^{\beta}\|F(X^{N}_{t_{m_s}})\|^2\Big]\mathrm{d}r\mathrm{d}s,
\end{align*}
which gives $\|J^2_{1,1}(t)\|_{L^p(\Omega)}\leq C\delta^{\beta}$ based on Proposition \ref{coro}.
And when $\beta\in(1,2)$,
since the order of the H\"older continuity of $X^N_r$ is  $\frac12$ in $H$-norm, hence the term $J^2_{1,1}$ need to be further split, based on the Taylor formula to $F$ 
and the fact that
\begin{align*}
X^{N}_r-X^{N}_{t_{m_r}}&=(S^N(r-t_{m_r})-\mathrm{Id})X^{N}_{t_{m_r}}+\int_{t_{m_r}}^rS^N(r-t_{m_r})F^N(X^{N}_{t_{m_r}})\mathrm{d}u+\int_{t_{m_r}}^rS^N(r-t_{m_r})\mathrm{d}W(u).
\end{align*}
Namely, we arrive at
\begin{align*}
&J^2_{1,1}(t)
= -\int_0^t\Big\langle \int_0^sA^{\frac{\beta-1}{2}}S^N(s-r)DF(X^{N}_{t_{m_r}})(S^N(r-t_{m_r})-\mathrm{Id})X^{N}_{t_{m_r}}\mathrm{d}r,\\
& \quad A^{-\frac{\beta}{2}}(S^N(s-t_{m_s})-\mathrm{Id})A^{\frac{1}{2}}F^N(X^N_{t_{m_s}})\Big\rangle\mathrm{d}s\\
&\quad -\int_0^t\Big\langle\int_0^sA^{\frac{\beta-1}{2}}S^N(s-r)DF(X^{N}_{t_{m_r}})\int_{t_{m_r}}^rS^N(r-t_{m_r})F^N(X^{N}_{t_{m_r}})\mathrm{d}u\mathrm{d}r,\\
&\quad A^{-\frac{\beta}{2}}(S^N(s-t_{m_s})-\mathrm{Id})A^{\frac{1}{2}}F^N(X^N_{t_{m_s}})\Big\rangle\mathrm{d}s\\
&\quad -\int_0^t\Big\langle\int_0^sA^{\frac{\beta-1}{2}}S^N(s-r)DF(X^{N}_{t_{m_r}})\int_{t_{m_r}}^rS^N(r-t_{m_r})\mathrm{d}W(u)\mathrm{d}r,\\
&\quad A^{-\frac{\beta}{2}}(S^N(s-t_{m_s})-\mathrm{Id})A^{\frac{1}{2}}F^N(X^N_{t_{m_s}})\Big\rangle\mathrm{d}s\\
&\quad -\int_0^t\Big\langle \int_0^sS^N(s-r)R_F(X^{N}_{t_{m_r}},X^{N}_r)\mathrm{d}r,\;(S^N(s-t_{m_s})-\mathrm{Id})F^N(X^N_{t_{m_s}})\Big\rangle\mathrm{d}s\\
&=:II_1(t)+II_2(t)+II_3(t)+II_4(t)
\end{align*}
with the remainder $$R_F(X^{N}_{t_{m_r}},X^{N}_r):=\int_0^1D^2F\big(X^{N}_{t_{m_r}}+(1-\theta)(X^{N}_r-X^{N}_{t_{m_r}})\big)\big((X^{N}_r-X^{N}_{t_{m_r}}),(X^{N}_r-X^{N}_{t_{m_r}})\big)(1-\theta)\mathrm{d}\theta.$$
Then we treat the four terms $II_i,i=1,2,3,4$ one by one.
It follows from \eqref{S(t)1}-\eqref{S(t)2} and  \eqref{DF_prop} that
\begin{align*}
&II_1(t)\leq C\int_0^t\int_0^s(s-r)^{-\frac{\beta-1}{2}}(1+\|X^{N}_{t_{m_r}}\|^2_E)(r-t_{m_r})^{\frac{\beta}{2}}\|X^{N}_{t_{m_r}}\|_{\beta}\mathrm{d}r(s-t_{m_s})^{\frac{\beta}{2}}\|F(X^N_{t_{m_s}})\|_1\mathrm{d}s,\\
&II_2(t)\leq C\int_0^t\int_0^s(s-r)^{-\frac{\beta-1}{2}}(1+\|X^{N}_{t_{m_r}}\|^2_E)(r-t_{m_r})\|F(X^{N}_{t_{m_r}})\|\mathrm{d}r(s-t_{m_s})^{\frac{\beta}{2}}\|F(X^N_{t_{m_s}})\|_1\mathrm{d}s,
\end{align*} 
which combining  Propositions \ref{thm2} and  \ref{coro} gives $\|II_1(t)\|_{L^p(\Omega)}+\|II_2(t)\|_{L^p(\Omega)}\leq C\delta^{\beta}, \; t\in[0,T].$
For the term $II_3,$ applying the stochastic Fubini theorem (see e.g. \cite[Theorem 4.18]{Kruse14}) and the Burkholder–Davis–Gundy inequality (see e.g. \cite[Proposition 2.12]{Kruse14}), and combining Propositions \ref{thm2} and \ref{coro} yield 
\begin{align}\label{II_3}
&\quad \|II_3(t)\|_{L^{p}(\Omega)}\notag\\&\leq \int_0^t\Big\|\sum_{i=0}^{m_s}\int_{t_i}^{t_{i+1}}\int_{t_i}^{t_{i+1}}\chi_{[t_i,r)}(u)\chi_{[t_i,s)}(r)A^{\frac{\beta-1}{2}}S^N(s-r)DF(X^{N}_{t_i})S^N(r-t_i)\mathrm{d}W(u)\mathrm{d}r\Big\|_{L^{2p}(\Omega;H)}\times \notag\\
&\qquad \Big\|(s-t_{m_s})^{\frac{\beta}{2}}\|F(X^N_{t_{m_s}})\|_1\Big\|_{L^{2p}(\Omega)}\mathrm{d}s\notag\\
&\leq C\delta^{\frac{\beta}{2}}\int_0^t\Big\|\sum_{i=0}^{m_s}\int_{t_i}^{t_{i+1}}\int_{t_i}^{t_{i+1}}\chi_{[t_i,r)}(u)\chi_{[t_i,s)}(r)A^{\frac{\beta-1}{2}}S^N(s-r)DF(X^{N}_{t_i})S^N(r-t_i)\mathrm{d}r\mathrm{d}W(u)\Big\|_{L^{2p}(\Omega;H)}\mathrm{d}s\notag\\
&\leq C\delta^{\frac{\beta}{2}}\int_0^t\Big(\int_0^{s}\Big\|\int_{t_{m_u}}^{t_{m_u}+\tau^{\delta}_{m_u}}\chi_{[t_{m_u},r)}(u)\chi_{[t_{m_u},s)}(r)A^{\frac{\beta-1}{2}}S^N(s-r)DF(X^{N}_{t_{m_u}})\notag\\
&\quad S^N(r-t_{m_u})Q^{\frac12}\mathrm{d}r\Big\|^2_{L^{2p}(\Omega;\mathcal{L}_2(H))}\mathrm{d}u\Big)^{\frac12}\mathrm{d}s\notag\\
&\leq C\delta^{\frac{1+\beta}{2}}\int_0^t\Big[\Big(\int_0^{(s-T\delta)\vee0}\Big\|\sum_{j=1}^{\infty}\int_{t_{m_u}}^{t_{m_u}+\tau^{\delta}_{m_u}}\|A^{\frac{\beta-1}{2}}S^N(s-r)DF(X^{N}_{t_{m_u}})S^N(r-t_{m_u})Q^{\frac{1}{2}}e_j\|^2\mathrm{d}r\Big\|_{L^{p}(\Omega)}\mathrm{d}u\Big)^{\frac12}\notag\\
&\quad +\Big(\int_{(s-T\delta)\vee0}^s\Big\|\sum_{j=1}^{\infty}\int_{(s-2T\delta)\vee0}^{s}\|A^{\frac{\beta-1}{2}}S^N(s-r)DF(X^{N}_{t_{m_u}})S^N(r-t_{m_u})Q^{\frac{1}{2}}e_j\|^2\mathrm{d}r\Big\|_{L^{p}(\Omega)}\mathrm{d}u\Big)^{\frac12}\Big]
\mathrm{d}s\notag\\
&\leq C\delta^{\frac{1+\beta}{2}}\int_0^t\Big[\delta^{\frac12}\Big(\int_0^{(s-T\delta)\vee0}(s-u-T\delta)^{-\beta+1}\mathrm{d}u\Big)^{\frac12}+\Big(\int_{(s-T\delta)\vee0}^s\int_{(s-2T\delta)\vee0}^s(s-r)^{-\beta+1}\mathrm{d}r\mathrm{d}u\Big)^{\frac12}\Big]\mathrm{d}s\times \notag\\
&\quad\sup_{s\in[0,T]}\|DF(X^N_s)A^{-\frac{\beta-1}{2}}\|_{L^{2p}(\Omega;\mathcal{L}(H))}\|A^{\frac{\beta-1}{2}}Q^{\frac12}\|_{\mathcal{L}_2(H)}.
\end{align}
It can be calculated that
\begin{align*}
&\quad\int_0^t\Big[\delta^{\frac12}\Big(\int_0^{(s-T\delta)\vee0}(s-u-T\delta)^{-\beta+1}\mathrm{d}u\Big)^{\frac12}+\Big(\int_{(s-T\delta)\vee0}^s\int_{(s-2T\delta)\vee0}^s(s-r)^{-\beta+1}\mathrm{d}r\mathrm{d}u\Big)^{\frac12}\Big]\mathrm{d}s\\
&\leq \int_{T\delta}^t\delta^{\frac12}\Big(\int_0^{s-T\delta}(s-u-T\delta)^{-\beta+1}\mathrm{d}u\Big)^{\frac12}+\int_0^{T\delta} \Big(\int_0^s\int_0^s(s-r)^{-\beta+1}\mathrm{d}s\mathrm{d}u\Big)^{\frac12}\mathrm{d}s\\
&\quad +\int_{T\delta}^{2T\delta}\Big(\int_{s-T\delta}^s\int_0^s(s-r)^{-\beta+1}\mathrm{d}r\mathrm{d}u\Big)^{\frac12}\mathrm{d}s+\int_{2T\delta}^t\Big(\int_{s-T\delta}^s\int_{s-2T\delta}^s(s-r)^{-\beta+1}\mathrm{d}r\mathrm{d}u\Big)^{\frac12}\mathrm{d}s\\
&\leq C(\delta^{\frac12}+\delta^{\frac{5-\beta}{2}}+\delta^{\frac{3-\beta}{2}})\leq C\delta^{\frac12},
\end{align*}
which leads to $\|II_3(t)\|_{L^p(\Omega)}\leq C\delta^{1+\frac{\beta}{2}},\;t\in[0,T].$
The term $II_4$ is treated separately for cases $d=1,2$ and $d=3.$ When $d=1,2,$ by the Sobolev embedding $L^1(\mathcal{O})\hookrightarrow \dot{H}^{-\frac{d+\epsilon_0}{2}}$ with $\epsilon_0=\frac{2-\beta}{2}\, (\beta\in(1,2))$,  
we deduce from  \eqref{S(t)1}-\eqref{S(t)2} and  \eqref{D^2F_prop}  that 
\begin{align*}
\|II_4(t)\|_{L^{p}(\Omega)}&\leq \int_0^t\int_0^s\|A^{\frac{\beta-1}{2}+\frac{d+\epsilon_0}{4}}S^N(s-r)\|_{\mathcal{L}(H)}\|A^{-\frac{d+\epsilon_0}{4}}R_F(X^{N}_{t_{m_r}},X^{N}_{r})\|_{L^{2p}(\Omega;H)}\mathrm{d}r\times\\
&\quad  \|A^{-\frac{\beta}{2}}(S^N(s-t_{m_s})-\mathrm{Id})A^{\frac12}F^N(X^N_{t_{m_s}})\|_{L^{2p}(\Omega;H)}\mathrm{d}s\\
&\leq C\delta^{\frac{\beta}{2}}\int_0^t\int_0^s(s-r)^{-\frac{\beta-1}{2}-\frac{d+\epsilon_0}{4}}\|R_F(X^{N}_{t_{m_r}},X^{N}_{r})\|_{L^{2p}(\Omega;L^1(\mathcal{O}))}\mathrm{d}r\|F(X^N_{t_{m_s}})\|_{L^{2p}(\Omega;\dot{H}^1)}\mathrm{d}s\\
&\leq C\delta^{\frac{\beta}{2}}\int_0^t\int_0^s(s-r)^{-\frac{\beta-1}{2}-\frac{d+\epsilon_0}{4}}\big(1+\sup_{0\leq r\leq T}\|X^{N}_r\|_{L^{4p}(\Omega;E)}\big)\big\|X^{N}_r-X^{N}_{t_{m_r}}\big\|^2_{L^{8p}(\Omega;H)}\mathrm{d}r\times\\
&\quad \|F(X^N_{t_{m_s}})\|_{L^{2p}(\Omega;\dot{H}^1)}\mathrm{d}s\leq C\delta^{1+\frac{\beta}{2}},
\end{align*}
where in the last step we use the regularity and the H\"older continuity of $X^N_t$ (see Propositions \ref{thm2} and \ref{coro}). 
And when $d=3,$ applying  \eqref{S(t)1}, \eqref{D^2F_prop}, the Gagliardo--Nirenberg  inequality $\|u\|_{L^4(\mathcal{O})}\leq C\|\nabla u\|^{\frac34}\|u\|^{\frac14},\,u\in H^1$, Propositions \ref{thm2} and \ref{coro} yields 
\begin{align*}
&\quad \|II_4(t)\|_{L^{p}(\Omega)}\\
&=\Big\|\int_0^t\Big\langle \int_0^sA^{\frac12}S^N(s-r)R_F(X^{N}_{t_{m_r}},X^{N}_r)\mathrm{d}r,\;A^{-1}(S^N(s-t_{m_s})-\mathrm{Id})A^{\frac12}F^N(X^N_{t_{m_s}})\Big\rangle\mathrm{d}s\Big\|_{L^p(\Omega)}\\
&\leq C\int_0^t\int_0^s(s-r)^{-\frac{1}{2}}\big(1+\sup_{0\leq r\leq T}\|X^{N}_{r}\|_{L^{4p}(\Omega;E)}\big)\big\|X^{N}_{r}-X^{N}_{t_{m_r}}\big\|^2_{L^{8p}(\Omega;L^4(\mathcal{O}))}\mathrm{d}r\Big\|(s-t_{m_s})\|F(X^N_{t_{m_s}})\|_1\Big\|_{L^{2p}(\Omega)}\mathrm{d}s\\
&\leq C\delta\int_0^t\int_0^s(s-r)^{-\frac{1}{2}}\big(1+\sup_{0\leq r\leq T}\|X^{N}_{r}\|_{L^{4p}(\Omega;E)}\big)\Big\|\|X^{N}_{r}-X^{N}_{t_{m_r}}\|^{\frac14}\|X^{N}_{r}-X^{N}_{t_{m_r}}\|_1^{\frac34}\Big\|_{L^{8p}(\Omega)}^2\mathrm{d}r\mathrm{d}s\\
&\leq C\delta\int_0^t\int_0^s(s-r)^{-\frac{1}{2}}\Big\|X^{N}_{r}-X^{N}_{t_{m_r}}\Big\|_{L^{4p}(\Omega;H)}^{\frac12}\Big\|X^{N}_{r}-X^{N}_{t_{m_r}}\Big\|^{\frac32}_{L^{12p}(\Omega;\dot{H}^1)}\mathrm{d}r\mathrm{d}s\\
&\leq C\delta\int_0^t\int_0^s(s-r)^{-\frac{1}{2}}(r-t_{m_r})^{\frac14+\frac{3(\beta-1)}{4}}\mathrm{d}r\mathrm{d}s
\leq C\delta^{\beta},
\end{align*}
where in the last step we have used the fact that $\frac14+\frac{3(\beta-1)}{4}> \beta-1$ for $\beta\in(1,2).$

Hence, we get that for $\beta\in(0,1)\cup(1,2),$
\begin{align*}
\|J^2_{1,1}(t)\|_{L^p(\Omega)}\leq C\delta^{\beta},\quad t\in[0,T].
\end{align*}

For the term $J^3_{1,1},$ when $\beta\in(0,2),$
\begin{align*}
J_{1,1}^3(t)&=\int_0^t\Big\langle \int_0^s A^{\frac{\beta}{2}}S^N(s-r)\big(F^N(X^{N}_{r})-F^N(X(r))\big)\mathrm{d}r, \;A^{-\frac{\beta}{2}}(S^N(s-t_{m_s})-\mathrm{Id})F^N(X^{N}_{t_{m_s}})\Big\rangle \mathrm{d}s\\
&\leq C\int_0^t\int_0^s(s-r)^{-\frac{\beta}{2}}\Big[\|X(r)-X^{N}_r\|^2+(s-t_{m_s})^{\beta}(1+\|X(r)\|^4_E+\|X^{N}_r\|^4_E)\|F(X^{N}_{t_{m_s}}\|^2)\Big]\mathrm{d}r\mathrm{d}s\\
&\leq C\int_0^t\|X^{N}_s-Y^N_s\|^2\mathrm{d}s+C\int_0^t\|Y^N_s-X(s)\|^2\mathrm{d}s\\
&\quad +C\int_0^t\int_0^s(s-r)^{-\frac{\beta}{2}}(s-t_{m_s})^{\beta}(1+\|X(r)\|^4_E+\|X^{N}_r\|^4_E)\|F(X^{N}_{t_{m_s}}\|^2)\mathrm{d}r\mathrm{d}s,
 \end{align*}
 where in the last step we transform the integral domain and use the Minkowski inequality. This, together with \eqref{X-P^NX}-\eqref{eq4}, implies that
 for $\beta\in(0,2),$
\begin{align*}
\|J^3_{1,1}(t)\|_{L^p(\Omega)}\leq C\int_0^t\|X^{N}_s-Y^N_s\|_{L^{2p}(\Omega;H)}^2\mathrm{d}s+C\lambda_N^{-\beta}+C\delta^{\beta}.
\end{align*}

Altogether, we obtain that for $\beta\in(0,2],$
\begin{align}\label{estJ_11}
\|J_{1,1}(t)\|_{L^p(\Omega)}\leq C\int_0^t\|X^{N}_s-Y^N_s\|_{L^{2p}(\Omega;H)}^2\mathrm{d}s+C\lambda_N^{-\beta}+C\delta^{\beta},\quad t\in[0,T].
\end{align}

\textit{Estimate of $J_{1,2}.$} When $\beta\in(0,1],$ the H\"older continuity of $X^N_t$ and Assumption \ref{assumption6} give 
\begin{align*}
\|J_{1,2}(t)\|_{L^p(\Omega)}&\leq \frac12\int_0^t\big\|\|X^N_s-Y^N_s\|^2\big\|_{L^p(\Omega)}\mathrm{d}s+C\int_0^t\|X^N_s-X^N_{t_{m_s}}\|^2_{L^{4p}(\Omega;H)}(1+\sup_{s\in[0,T]}\|X^N_s\|_{L^{8p}(\Omega;E)}^4)\mathrm{d}s\\
&\leq \frac12\int_0^t\big\|X^N_s-Y^N_s\big\|^2_{L^{2p}(\Omega;H)}\mathrm{d}s+C\delta^{\beta}.
\end{align*}
When $\beta\in(1,2],$ 
by \eqref{splittingX-Y}, we split $J_{1,2}$ further as $J_{1,2}=J^1_{1,2}+J^2_{1,2}+J^3_{1,2}$ with 
\begin{align*}
&J^1_{1,2}(t):=\int_0^t\Big\langle \int_0^sS^N(s-r)(S^N(r-t_{m_r})-\mathrm{Id})F^N(X^N_{t_{m_r}})\mathrm{d}r,\;F^N(X^N_{t_{m_s}})-F^N(X^N_s)\Big\rangle\mathrm{d}s,\\
&J^2_{1,2}(t):=\int_0^t\Big\langle \int_0^sS^N(s-r)(F^N(X^N_{t_{m_r}})-F^N(X^N_r))\mathrm{d}r,\;F^N(X^N_{t_{m_s}})-F^N(X^N_s)\Big\rangle\mathrm{d}s,\\
&J^3_{1,2}(t):=\int_0^t\Big\langle \int_0^sS^N(s-r)(F^N(X^N_{r})-F^N(X(r)))\mathrm{d}r,\;F^N(X^N_{t_{m_s}})-F^N(X^N_s)\Big\rangle\mathrm{d}s.
\end{align*}
The estimate of the term $J^1_{1,2}$ is  similar as that of $J^2_{1,1},$ whose proof is based on the Taylor formula to $F.$  And one can get $\|J^{1}_{1,2}(t)\|_{L^p(\Omega)}\leq C\delta^{\beta}.$ 
Now we show the estimate of the term $J^2_{1,2}=\sum_{i=1}^3J^{2,i}_{1,2},$ where
\begin{align*}
J^{2,1}_{1,2}(t)&:=\int_0^t\Big\langle \int_0^sS^N(s-r)\mathcal{A}_r\mathrm{d}r,\; \mathcal{A}_s\Big\rangle\mathrm{d}s,\\
J^{2,2}_{1,2}(t)&:=\int_0^t\Big\langle \int_0^sS^N(s-r)\mathcal{A}_r\mathrm{d}r,\; DF(X^N_{t_{m_s}})\int_{t_{m_s}}^s S^N(s-t_{m_s})\mathrm{d}W(u)\Big\rangle\mathrm{d}s\\
&\quad +\int_0^t\Big\langle \int_0^sS^N(s-r)DF(X^N_{t_{m_r}})\int_{t_{m_r}}^r S^N(r-t_{m_r})\mathrm{d}W(u)\mathrm{d}r,\;  \mathcal{A}_s\Big\rangle\mathrm{d}s,\\
J^{2,3}_{1,2}(t)&:=\int_0^t\Big\langle \int_0^sS^N(s-r)DF(X^N_{t_{m_r}})\int_{t_{m_r}}^rS^N(r-t_{m_r})\mathrm{d}W(u_2)\mathrm{d}r, \\
&~~\quad DF(X^N_{t_{m_s}})\int_{t_{m_s}}^sS^N(s-t_{m_s})\mathrm{d}W(u_1)\Big\rangle\mathrm{d}s
\end{align*}
with $\mathcal{A}_r:=DF(X^N_{t_{m_r}})\Big((S^N(r-t_{m_r})-\mathrm{Id})X^N_{t_{m_r}}+\int_{t_{m_r}}^rS^N(r-t_{m_r})F^N(X^N_{t_{m_r}})\mathrm{d}u\Big)+R_F(X^N_{t_{m_r}},X^N_r).$
For the term $J^{2,1}_{1,2},$
noting that for $d=1,$
\begin{align*}
&\quad \Big\|\int_0^t\Big\langle\int_0^sS^N(s-r)R_F(X^N_{t_{m_r}},X^N_r)\mathrm{d}r,R_F(X^N_{t_{m_s}},X^N_s)\Big\rangle\mathrm{d}s\Big\|_{L^p(\Omega)}\\
&\leq C\int_0^t\int_0^s(s-r)^{-\frac{d+\epsilon_0}{2}}\Big\|\|R_F(X^N_{t_{m_r}},X^N_r)\|_{-\frac{d+\epsilon_0}{2}}\|R_F(X^N_{t_{m_s}},X^N_s)\|_{-\frac{d+\epsilon_0}{2}}\Big\|_{L^p(\Omega)} \mathrm{d}r\mathrm{d}s\leq C\delta^2,
\end{align*}
and for $d=2,3,$
\begin{align*}
&\quad \Big\|\int_0^t\Big\langle\int_0^sS^N(s-r)R_F(X^N_{t_{m_r}},X^N_r)\mathrm{d}r,R_F(X^N_{t_{m_s}},X^N_s)\Big\rangle\mathrm{d}s\Big\|_{L^p(\Omega)}\\
&\leq C\int_0^t\int_0^s(s-r)^{-\frac{d+\epsilon_0}{4}}\Big\|\|R_F(X^N_{t_{m_r}},X^N_r)\|_{-\frac{d+\epsilon_0}{2}}\|R_F(X^N_{t_{m_s}},X^N_s)\|\Big\|_{L^p(\Omega)} \mathrm{d}r\mathrm{d}s\\
&\leq C\int_0^t\int_0^s(s-r)^{-\frac{d+\epsilon_0}{4}}\|R_F(X^N_{t_{m_r}},X^N_r)\|_{L^{2p}(\Omega;\dot{H}^{-\frac{d+\epsilon_0}{2}})}\sup_{s\in[0,T]}\|X^N_s\|_{L^{4p}(\Omega;E)}\|X^N_s-X^N_{t_{m_s}}\|^2_{L^{8p}(\Omega;L^4(\mathcal{O}))}\mathrm{d}r\mathrm{d}s\\
&\leq C\delta \int_0^t\int_0^s(s-r)^{-\frac{d+\epsilon_0}{4}}\big\|\|X^N_s-X^N_{t_{m_s}}\|^{\frac d4}_1\|X^N_s-X^N_{t_{m_s}}\|^{1-\frac{d}{4}}\big\|_{L^{8p}(\Omega)}^2\mathrm{d}r\mathrm{d}s\leq C\delta^{\beta},
\end{align*}
we get
\begin{align*}
\|J^{2,1}_{1,2}(t)\|_{L^p(\Omega)}&\leq \Big\|\int_0^t\Big\langle\int_0^sS^N(s-r)R_F(X^N_{t_{m_r}},X^N_r)\mathrm{d}r,R_F(X^N_{t_{m_s}},X^N_s)\Big\rangle\mathrm{d}s\Big\|_{L^p(\Omega)}\\
&\quad +C\Big[\delta^{\beta}\sup_{r\in[0,T]}\big\|\|DF(X^N_{t_{m_r}})\|_{\mathcal{L}(H)}\|X^N_{t_{m_r}}\|_{\beta}\big\|^2_{L^{2p}(\Omega)}\\
&\quad +\delta^{2}\sup_{r\in[0,T]}\big\|\|DF(X^N_{t_{m_r}})\|_{\mathcal{L}(H)}\|F(X^N_{t_{m_r}})\|\big\|^2_{L^{2p}(\Omega)}\\
&\quad +\delta^{\frac{\beta}{2}}\int_0^t\int_0^s(s-r)^{-\frac {d+\epsilon_0}{4}}\Big\|\|R_F(X^N_{t_{m_r}},X^N_r)\|_{-\frac {d+\epsilon_0}{2}}\|DF(X^N_{t_{m_s}})\|_{\mathcal{L}(H)}\\
&\quad\big(\|X^N_{t_{m_s}}\|_{\beta}+\|F(X^N_{t_{m_s}})\|\big)\Big\|_{L^p(\Omega)} \mathrm{d}r\mathrm{d}s\Big]
\leq C\delta^{\beta}.
\end{align*}
For the term $J^{2,2}_{1,2},$ applying the stochastic Fubini theorem yields
\begin{align*}
 &\quad\|J^{2,2}_{1,2}(t)\|_{L^p(\Omega)}\\
&\leq \Big\|\sum_{i=0}^{m_t}\int_{t_i}^{t_{i+1}}\chi_{[t_i,s)}(u)\chi_{[t_i,t}(s)\Big\langle \int_0^sS^N(s-r)\mathcal{A}_r\mathrm{d}r,\;DF(X^N_{t_i})S^N(s-t_i)\mathrm{d}s\mathrm{d}W(u)\Big\rangle\Big\|_{L^p(\Omega)}\\
&\quad +\int_0^t\Big\|\sum_{i=0}^{m_s}\int_{t_i}^{t_{i+1}}\int_{t_i}^{t_{i+1}}\chi_{[t_i,r)}(u)\chi_{[t_i,s)}(r)\Big\langle \mathcal{A}_s,\;S^N(s-r)DF(X^N_{t_i})S^N(r-t_i)\mathrm{d}r\mathrm{d}W(u)\Big\rangle
\Big\|_{L^p(\Omega)}\mathrm{d}s.
\end{align*}
And by a similar proof as that of \eqref{II_3}, one can show that $\|J^{2,2}_{1,2}(t)\|_{L^p(\Omega)}\leq C\delta ^{\frac{\beta}{2}+1}.$
For the term $J^{2,3}_{1,2},$
by using the stochastic Fubini theorem twice and letting $G(s):=\int_0^sS^N(s-r)DF(X^N_{t_{m_r}})\int_{t_{m_r}}^rS^N(r-t_{m_r})\mathrm{d}W(u_2)\mathrm{d}r$, we derive
\begin{align*}
&\quad\|J^{2,3}_{1,2}(t)\|_{L^{p}(\Omega)}\\
&= \Big\|\sum_{i=0}^{m_t}\int_{t_i}^{t_{i+1}}\int_{t_i}^{t_{i+1}}\Big\langle G(s),\;DF(X^N_{t_i})\chi_{[t_i,s)}(u_1)\chi_{[t_i,t)}(s)S^N(s-t_{i})\mathrm{d}s\mathrm{d}W(u_1)
\Big\rangle\Big\|_{L^{p}(\Omega)}\\
&\leq C\Big(\int_0^t\Big\|\int_{t_{m_{u_1}}}^{t_{m_{u_1}}+\tau^{\delta}_{m_{u_1}}}Q^{\frac12}S^N(s-t_{m_s})\chi_{[t_{m_{u_1}},s)}(u_1)\chi_{[t_{m_{u_1}},t)}(s)DF(X^N_{t_{m_s}})G(s)\mathrm{d}s\Big\|^2_{L^{p}(\Omega;H)}\mathrm{d}u_1\Big)^{\frac12}\\
&\leq C\delta ^{\frac12}\Big(\int_0^t\int_{(u_1-T\delta)\vee0}^{u_1+T\delta}\sup_{s\in[0,T]}\Big\|Q^{\frac12}S^N(s-t_{m_s})\chi_{[t_{m_{u_1}},s)}(u_1)\chi_{[t_{m_{u_1}},t)}(s)DF(X^N_{t_{m_s}})G(s)\Big\|^2_{L^{p}(\Omega;H)}\mathrm{d}s\mathrm{d}u_1\Big)^{\frac12}\\
&\leq C\delta\sup_{s\in[0,T]}\Big\|\sum_{j=0}^{m_s}\int_{t_j}^{t_{j+1}}\int_{t_j}^{t_{j+1}}Q^{\frac12}S^N(s-t_{m_s})DF(X^N_{t_{m_s}})S^N(s-r)DF(X^N_{t_j})\chi_{[t_j,r)}(u_2)\chi_{[t_j,s)}(r) \\
&\quad S^N(r-t_j)\mathrm{d}r\mathrm{d}W(u_2)
\Big\|_{L^{p}(\Omega;H)}\\
&\leq C\delta \sup_{s\in[0,T]}\Big(\int_0^s\Big\|\int_{t_{m_{u_2}}}^{t_{m_{u_2}}+\tau^{\delta}_{m_{u_2}}}Q^{\frac12}S^N(s-t_{m_s})DF(X^N_{t_{m_s}})S^N(s-r)DF(X^N_{t_{m_{u_2}}})\chi_{[t_{m_{u_2}},r)}(u_2)\chi_{[t_{m_{u_2}},s)}(r)\\
&\quad S^N(r-t_{m_{u_ 2}})Q^{\frac12}\mathrm{d}r\Big\|^2_{L^{p}(\Omega;\mathcal{L}_2(H))}\mathrm{d}u_2\Big)^{\frac12}\leq C\delta^2.
\end{align*}
For the term $J^3_{1,2},$ we split it as $J^3_{1,2}=J^{3,1}_{1,2}+J^{3,2}_{1,2}$ with
\begin{align*}
J^{3,1}_{1,2}(t)&:=\int_0^t\Big\langle \int_0^sS^N(s-r)(F^N(X^N_r)-F^N(X(r)))\mathrm{d}r,\; \mathcal{A}_s\Big\rangle\mathrm{d}s,\\
J^{3,2}_{1,2}(t)&:=\int_0^t\Big\langle \int_0^sS^N(s-r)(F^N(X^N_r)-F^N(X(r)))\mathrm{d}r,\;DF(X^N_{t_{m_s}})\int_{t_{m_s}}^sS^N(s-t_{m_s})\mathrm{d}W(u)\Big\rangle\mathrm{d}s.
\end{align*}
Combining  \eqref{X-P^NX}-\eqref{eq4}, the term $J^{3,1}_{1,2}$ can be estimated as
\begin{align*}
\|J^{3,1}_{1,2}(t)\|_{L^p(\Omega)}&\leq C\int_0^t\int_0^s(s-r)^{-\frac{d+\epsilon_0}{4}}\Big[\|X^N_r-X(r)\|_{L^{2p}(\Omega;H)}^2+\big(1+\|X^N_r\|^4_{L^{8p}(\Omega;E)}+\|X(r)\|^4_{L^{8p}(\Omega;E)}\big)\times \\
&\quad \Big(\delta^{\beta}\Big\|\|DF(X^N_{t_{m_s}})\|_{\mathcal{L}(H)}\big(\|X^N_{t_{m_s}}\|_{\beta}+\|F(X^N_{t_{m_s}})\|\big)\Big\|^2_{L^{4p}(\Omega)}\\
&\quad +\Big\|\|R_F(X^N_{t_{m_s}},X^N_s)\|_{-\frac{d+\epsilon_0}{2}}\Big\|^2_{L^{4p}(\Omega)}\Big)\Big]\mathrm{d}r\mathrm{d}s\\
&\leq C\int_0^t\|X^N_s-Y^N_s\|^2_{L^{2p}(\Omega;H)}\mathrm{d}s+C\delta^{\beta}+C\lambda_N^{-\beta}.
\end{align*}
By the stochastic Fubini theorem and the Young inequality, the term $J^{3,2}_{1,2}$ can be estimated as
\begin{align*}
\|J^{3,2}_{1,2}(t)\|_{L^p(\Omega)}
&=\Big\|\sum_{i=0}^{m_t}\int_{t_i}^{t_{i+1}}\int_{t_i}^{t_{i+1}}\Big\langle \int_0^s S^N(s-r)(F^N(X^N_r)-F^N(X(r)))\mathrm{d}r,\\
&\qquad DF(X^N_{t_i})S^N(s-t_i)\chi_{[t_i,s)}(u)\chi_{[t_i,t)}(s)\mathrm{d}s\mathrm{d}W(u)\Big\rangle\Big\|_{L^p(\Omega)}\\
&\leq C\Big(\int_0^t\Big\|\int_{t_{m_u}}^{t_{m_u}+\tau^{\delta}_{m_u}}Q^{\frac12}S^N(t-t_{m_u})\chi_{[t_{m_u},s)}(u)\chi_{[t_{m_u},t)}(s)DF(X^N_{t_{m_u}})\\
&\quad \int_0^sS^N(s-r)(F^N(X^N_r)-F^N(X(r)))\mathrm{d}r\mathrm{d}s\Big\|^2_{L^p(\Omega;H)}\mathrm{d}u\Big)^{\frac12}\\
&\leq C\Big(\int_0^t\Big\|\int_{t_{m_u}}^{t_{m_u}+\tau^{\delta}_{m_u}}\|Q^{\frac12}S^N(t-t_{m_u})\chi_{[t_{m_u},s)}(u)\chi_{[t_{m_u},t)}(s)DF(X^N_{t_{m_u}})\|_{\mathcal{L}(H)}\times\\
&\quad \int_0^t(1+\|X^N_r\|^2_E+\|X(r)\|^2_E)\|X^N_r-X(r)\|\mathrm{d}r\mathrm{d}s\Big\|_{L^p(\Omega)}^2\mathrm{d}u\Big)^{\frac12}\\
&\leq C\Big(\int_0^t\Big[\int_0^t\Big\|\int_{t_{m_u}}^{t_{m_u}+\tau^{\delta}_{m_u}}\|Q^{\frac12}S^N(t-t_{m_u})\chi_{[t_{m_u},s)}(u)\chi_{[t_{m_u},t)}(s)DF(X^N_{t_{m_u}})\|_{\mathcal{L}(H)}\times\\
&\quad (1+\|X^N_r\|^2_E+\|X(r)\|^2_E)\mathrm{d}s\Big\|_{L^{2p}(\Omega)}\|X^N_r-X(r)\|_{L^{2p}(\Omega;H)}\mathrm{d}r\Big]^2\mathrm{d}u\Big)^{\frac12}\\
&\leq C\int_0^t\|X^N_r-X(r)\|^2_{L^{2p}(\Omega;H)}\mathrm{d}r+C\Big(\int_0^t\Big[\int_0^t\Big\|\int_{t_{m_u}}^{t_{m_u}+\tau^{\delta}_{m_u}}\chi_{[t_{m_u},s)}(u)\chi_{[t_{m_u},t)}(s)\\
&\quad \|Q^{\frac12}S^N(t-t_{m_u})DF(X^N_{t_{m_u}})\|_{\mathcal{L}(H)}(1+\|X^N_r\|^2_E+\|X(r)\|^2_E)\mathrm{d}s\Big\|^2_{L^{2p}(\Omega)}\mathrm{d}r\Big]^2\mathrm{d}u\Big)^{\frac12}\\
&\leq C\int_0^t\|X^N_r-Y^N_r\|^2_{L^{2p}(\Omega;H)}\mathrm{d}r+C\delta^{\beta}+C\lambda_N^{-\beta},
\end{align*}
where in the third step we exchange the order of the integral and use the H\"older inequality, and in the last step we use \eqref{X-P^NX}-\eqref{eq4}.

Hence, we derive that for $\beta\in(0,2],$ \begin{align}\label{estJ_12}
\|J_{1,2}(t)\|_{L^p(\Omega)}\leq C\int_0^t\|X^N_r-Y^N_r\|^2_{L^{2p}(\Omega;H)}\mathrm{d}r+C\delta^{\beta}+C\lambda_N^{-\beta}\quad t\in[0,T].
\end{align}

Therefore, combining estimates of terms $J_{1,1}$ and $J_{1,2}$ (i.e., \eqref{estJ_11} and \eqref{estJ_12}), we get 
$$
\left\|X^{N}_t-Y^N_t\right\|^2_{L^{2p}(\Omega;H)}\leq C\Big(\int_0^t\left\|X^{N}_s-Y^N_s\right\|^2_{L^{2p}(\Omega;H)}\mathrm{d}s+\lambda_N^{-\beta}+\delta^{\beta}\Big),
$$
where $p\ge 1,$ and the constant $C>0$ depends on $p,T,\tau_{min},L_0,L_1,L_2, X_0$, and $\|A^{\frac{\beta-1}{2}}Q^{\frac{1}{2}}\|_{\mathcal{L}_2(H)}.$ 
%Recalling that $\lambda_N\sim N^{-\frac 2d},$ and 
Applying the Gr\"onwall inequality, we obtain
$$
\sup_{0\leq t\leq T}\|X^N_t-Y^N_t\|_{L^{2p}{(\Omega;H)}}\leq C(\lambda_N^{-\frac{\beta}{2}}+\delta^{\frac{\beta}{2}}),\;p\ge 1,
$$
which together with \eqref{X-P^NX}-\eqref{eq4} finishes the proof of the scheme \eqref{scheme}.
\end{proof}

\section{Discussions on the multiplicative noise case}\label{multi_discuss}

For the stochastic Allen--Cahn equation  driven by a multiplicative noise, there have been some works on the numerical study, see  e.g. \cite{jentzen, liuzh, prohl18, fengxb21} and references therein. 
To be specific, for the stochastic Allen--Cahn equation driven by the Brownian motion, authors in \cite{fengxb21} propose the fully discrete finite element method and prove the strong convergence with nearly optimal rates; authors in \cite{prohl18} give variational error analysis for the structure preserving finite element based space-time discretization of  the  strong variational solution. For the stochastic Allen--Cahn equation driven by the $Q$-Wiener process, 
\cite{jentzen} derives the strong convergence rate for the nonlinearity-truncated exponential Euler approximation scheme; \cite{liuzh} proves strong convergence rates for both the drift-implicit Euler--Galerkin finite element scheme  and the Milstein--Galerkin finite element scheme.  

In this section, we present the numerical analysis of the adaptive time-stepping scheme for the multiplicative noise case, i.e., the stochastic Allen--Cahn equation with the multiplicative noise
\begin{equation}
\begin{cases}\label{multi_spde}
\mathrm{d} X(t)+AX(t)\mathrm{d}t=F(X(t))\mathrm{d}t+G(X(t))\mathrm{d}W(t),\quad t\in(0,T],\\
X(0)=X_0,
\end{cases}
\end{equation}  
where $A$ and $F$ are defined as in \eqref{spde}. In this section,  $\{W(t)\}_{t\in [0,T]}$ is a generalized $Q$-Wiener process valued on another separable Hilbert space $U$.
% $U_1\supset U$ with $U$  being a separable Hilbert space. 
Here,  $G:H\to\mathcal{L}^0_2:=\mathcal{L}_2(U_0,H)$  is Lipschitz continuous with $U_0:=Q^{\frac12}(U)$, i.e.,
\begin{align*}
\|G(X)-G(Y)\|_{\mathcal{L}^0_2}\leq C\|X-Y\|,\quad X,Y\in H
\end{align*}
for some constant $C>0.$
The corresponding adaptive time-stepping scheme  for \eqref{multi_spde} is
\begin{align}\tag{AE-M}\label{scheme_multi}
X^N_{t_{m+1}}=S^N(\tau_m)\big(X^N_{t_m}+F^N(X^N_{t_m})\tau_m+P^NG(X^N_{t_m})\Delta W_m\big),\quad X^N_0=P^NX_0,
\end{align}
where $\tau_m$ is defined similarly as in \eqref{scheme}.
The differential form of the  continuous version of \eqref{scheme_multi} is given by
\begin{align}
\label{conti_multi}
\mathrm{d}X^N_t=-A^NX^N_t\mathrm{d}t+S^N(t-t_{m_t})F^N(X^N_{t_{m_t}})\mathrm{d}t+S^N(t-t_{m_t})P^NG(X^N_{t_{m_t}})\mathrm{d}W(t). 
\end{align}
We remark that in this section we still use  notations $X(t),X^N_{t_{m}},$ and $X^N_t$ in the multiplicative noise case for the simplicity of symbols.  
And we will emphasize on them when these symbols  are referred to solutions in the additive noise case.

%Recall the whole numerical analysis of \eqref{scheme}, the prerequisite of the convergence analysis is the \textit{a priori} estimates in $E$-norm and $\dot{H}^{\beta}$-norm, while the latter one can be obtained by the use of the former one. The same idea can be used in the multiplicative noise case, namely, 

%We first give the \textit{a priori} estimate  in $E$-norm of the solution of \eqref{scheme_multi}. 

%As stated in Remark \ref{rmkbeta=1}, for the adaptive time-stepping scheme  \eqref{scheme} of \eqref{spde}, one can obtain $\mathbb{E}\Big[\sup_{0\leq t\leq T}\|X^N_t\|^p\Big]<\infty$ for some $p>2$ under the assumption \eqref{adaptive3}for the trace-class noise case. This can be generalized to the multiplicative case \eqref{multi_spde}, namely, we have $\mathbb{E}\Big[\sup_{0\leq t\leq T}\|X^N_t\|^p\Big]<\infty$ for $p\ge 1$ under the assumption \eqref{adaptive3} with $X^N_t$ being the solution of \eqref{multi_spde}.
The following lemma shows the boundedness of the  expected supremum of the $p$th moment for the solution of \eqref{scheme_multi}.
\begin{lemma}
Under Assumptions \ref{assumption1},\ref{assumption3}, and \eqref{adaptive5}-\eqref{adaptive3}, we have $\mathbb{E}\big[\sup_{0\leq t\leq T}\|X^N_t\|^p\big]\le C$ for  $p\ge 2$ with some constant $C:=C(p,T,X_0)>0$.
\end{lemma}
\begin{proof}
We only consider the case of $p\ge 4$ since the case of $2\leq p<4$ can be obtained by the use of the  H\"older inequality. 
The proof is based on  the truncation technique. 
For $K>\|P^NX_0\|,$ define the truncated function $\Theta_K\in\mathcal{C}^{\infty}_0:H\to[0,1]$ satisfying $\Theta_K(x)=1$ for $\|x\|\leq K$ and $\Theta_K(x)=0$ for $\|x\|\ge 2K.$
%operator $P_K$ by $P_Ku:=\min\{1,\frac{K}{\|u\|}\}u$ for some constant $K>\|X_0\|$ and any $u\in H.$ 
Introduce the truncated stochastic process of the numerical solution $\{X^N_{t_m}\}_{m\ge 0}$ as 
\begin{align*}
X^{N,K}_{t_{m+1}}=S^N(\tau_m)\big(X^{N,K}_{t_m}+\Theta_K(X^{N,K}_{t_m})F^N(X^{N,K}_{t_m})\tau_m+P^NG(X^{N,K}_{t_m})\Delta W_m\big).
\end{align*} 
The corresponding continuous version is defined as
\begin{align*}
X^{N,K}_t=S^N(t)P^NX_0+\int_0^tS^N(t-t_{m_s})\Theta_K(X^{N,K}_{t_{m_s}})F^N(X^{N,K}_{t_{m_s}})\mathrm ds+\int_0^tS^N(t-t_{m_s})P^NG(X^{N,K}_{t_{m_s}})\mathrm{d}W(s).
\end{align*}
%Noticing that the adaptive timestep function $\tau(\cdot)$ is continuous, hence $\min_{\|x\|\leq K}\tau(x)>0,$ which implies that the time $T$ is always attainable  for the truncated numerical solution $X^{N,K}_{t_m}.$
Then the proof is separated into two steps.

\textit{Step 1: Expected supremum of the $p$th moment of the truncated numerical solution.}

By the contraction property of $P^N$ and $S(t)$ in $H$, i.e., $\|P^Nu\|\leq \|u\|,$ $\|S(t)u\|\leq \|u\|$, and using \eqref{adaptive3}, we have
\begin{align*}
&\|X^{N,K}_{t_{m+1}}\|^2
\leq \|X^{N,K}_{t_m}+\Theta_K(X^{N,K}_{t_m})F^N(X^{N,K}_{t_m})\tau_m+P^NG(X^{N,K}_{t_m})\Delta W_m\|^2\\
&= \|X^{N,K}_{t_m}\|^2+2\tau_m\Theta_K(X^{N,K}_{t_m})\Big(\langle X^{N,K}_{t_m},F^N(X^{N,K}_{t_m})\rangle+\frac{1}{2}\tau_m\Theta_K(X^{N,K}_{t_m})\|F^N(X^{N,K}_{t_m})\|^2\Big)\\
&\quad +2\Big\langle X^{N,K}_{t_m}+\Theta_K(X^{N,K}_{t_m})F^N(X^{N,K}_{t_m})\tau_m,\,P^N G(X^{N,K}_{t_m})\Delta W_m\Big\rangle+\|P^NG(X^{N,K}_{t_m})\Delta W_m\|^2\\
&\leq \|X^{N,K}_{t_m}\|^2+2\tau_m(\bar{L}_2\|X^{N,K}_{t_m}\|^2+L_3)+2\Big\langle X^{N,K}_{t_m}+\Theta_K(X^{N,K}_{t_m})F^N(X^{N,K}_{t_m})\tau_m,G(X^{N,K}_{t_m})\Delta W_m\Big\rangle\\
&\quad +\|G(X^{N,K}_{t_m})\Delta W_m\|^2.
\end{align*}
Similarly, for $t\in(t_m,t_{m+1}),$
\begin{align*}
&\|X^{N,K}_t\|^2\leq \|X^{N,K}_{t_m}+\Theta_K(X^{N,K}_{t_m})F^N(X^{N,K}_{t_m})(t-t_m)+P^NG(X^{N,K}_{t_m})(W(t)-W(t_m))\|^2\\
&\leq \|X^{N,K}_{t_m}\|^2+2(t-t_m)(\bar{L}_2\|X^{N,K}_{t_m}\|^2+L_3)\\
&\quad +2\Big\langle X^{N,K}_{t_m}+\Theta_K(X^{N,K}_{t_m})F^N(X^{N,K}_{t_m})(t-t_m),G(X^{N,K}_{t_m})(W(t)-W(t_m))\Big\rangle\\
&\quad +\|G(X^{N,K}_{t_m})(W(t)-W(t_m))\|^2.
\end{align*}
By induction, we have
\begin{align*}
\|X^{N,K}_{t}\|^2
&\leq \|P^NX_0\|^2+2\bar{L}_2\int_0^{t}\|X^{N,K}_{t_{m_s}}\|^2\mathrm{d}s+2L_3t\\
&\quad +2\int_0^{t}\Big\langle X^{N,K}_{t_{m_s} }+\Theta_K(X^{N,K}_{t_{m_s}})F^N(X^{N,K}_{t_{m_s} })((t_{m_s}+\tau_{m_s})\wedge t-t_{m_s}),\\
&\quad \,G(X^{N,K}_{t_{m_s}})\mathrm{d} W(s)\Big\rangle+\sum_{i=0}^{m_t-1}\|G(X^{N,K}_{t_i})\Delta W_i\|^2 +\|G(X^{N,K}_{t_{m_t}})(W(t)-W(t_{m_t}))\|^2.
\end{align*}
Taking the $\frac{p}{2}$th moment with $p\ge 4$ and applying the Burkholder--Davis--Gundy inequality  yield
\begin{align*}
\mathbb{E}[\|X^{N,K}_{t}\|^p]&
\leq C(p,T,\bar{L}_2)\Big(\mathbb{E}[\|X_0\|^p]+T^{\frac{p}{2}}+\int_0^{t}\mathbb{E}[\|X^{N,K}_{t_{m_s}}\|^p]\mathrm{d}s\\
&\quad +\mathbb{E}\Big[\Big(\int_0^{t}\Big\|X^{N,K}_{t_{m_s}}+\Theta_K(X^{N,K}_{t_{m_s}})F^N(X^{N,K}_{t_{m_s}})((t_{m_s}+\tau_{m_s})\wedge t-t_{m_s})\Big\|^{2}\|G(X^{N,K}_{t_{m_s}})\|^{2}_{\mathcal{L}^0_2}\mathrm{d}s\Big)^{\frac{p}{4}}\Big]
\\
&\quad +\mathbb{E}\Big[\Big(\sum_{i=0}^{m_t-1}\|G(X^{N,K}_{t_i})\Delta W_i\|^2+\|G(X^{N,K}_{t_{m_t}})(W(t)-W(t_{m_t}))\|^2\Big)^{\frac{p}{2}}\Big]\Big).
\end{align*}
Notice that
\begin{align*}
&\quad \|X^{N,K}_{t_{m_s}}+\Theta_K(X^{N,K}_{t_{m_s}})F^N(X^{N,K}_{t_{m_s}})((t_{m_s}+\tau_{m_s})\wedge t-t_{m_s})\|^2\\
&=\|X^{N,K}_{t_{m_s}}\|^2+2((t_{m_s}+\tau_{m_s})\wedge t-t_{m_s})\Theta_K(X^{N,K}_{t_{m_s}})\times\\
&\quad \Big(\langle X^{N,K}_{t_{m_s}},F(X^{N,K}_{t_{m_s}})\rangle+\frac{1}{2}((t_{m_s}+\tau_{m_s})\wedge t-t_{m_s})\Theta_K(X^{N,K}_{t_{m_s}})\|F(X^{N,K}_{t_{m_s}})\|^2\Big)\\
&\leq \|X^{N,K}_{t_{m_s}}\|^2+2((t_{m_s}+\tau_{m_s})\wedge t-t_{m_s})(\bar{L}_2\|X^{N,K}_{t_{m_s}}\|^2+L_3)
\end{align*}
and
\begin{align*}
&\mathbb{E}\Big[\Big(\sum_{i=0}^{m_t-1}\|G(X^{N,K}_{t_i})\Delta W_i\|^2\Big)^{\frac{p}{2}}\Big]=
\mathbb{E}\Big[\Big(\sum_{i=0}^{m_t-1}\tau_i\frac{\|G(X^{N,K}_{t_i})\Delta W_i\|^2}{\tau_i}\Big)^{\frac{p}{2}}\Big]\\
&\leq \mathbb{E}\left[\Big(\sum_{i=0}^{m_t-1}\tau_i\Big)^{\frac{p}{2}-1}\sum_{i=0}^{m_t-1}\tau_i\frac{\|G(X^{N,K}_{t_i})\Delta W_i\|^p}{\tau^{\frac{p}{2}}_i}\right]\\
&\leq T^{\frac{p}{2}-1}\mathbb{E}\left[\sum_{i=0}^{m_t-1}\tau_i\frac{\|G(X^{N,K}_{t_i})\Delta W_i\|^p}{\tau^{\frac{p}{2}}_i}\right]\leq T^{\frac{p}{2}-1}\mathbb{E}\left[\int_0^t\frac{\|G(X^{N,K}_{t_{m_s}})\Delta W_{m_s}\|^p}{\tau^{\frac{p}{2}}_{m_s}}\mathrm{d}s\right]\\
&= T^{\frac{p}{2}-1}\int_0^t\mathbb{E}\Big[\mathbb{E}\big[\tau^{-\frac{p}{2}}_{m_s}\|G(X^{N,K}_{t_{m_s}})\Delta W_{m_s}\|^p\big|\mathcal{F}_{t_{m_s}}\big]\Big]\mathrm{d}s=T^{\frac{p}{2}-1}\int_0^t\mathbb{E}\Big[\tau^{-\frac{p}{2}}_{m_s}\mathbb{E}\big[\|G(X^{N,K}_{t_{m_s}})\Delta W_{m_s}\|^p\big|\mathcal{F}_{t_{m_s}}\big]\Big]\,\mathrm{d}s\\
&\leq C\int_0^t\mathbb{E}[1+\|X^{N,K}_{t_{m_s}}\|^p]\mathrm{d}s,
\end{align*}
where we have used the inequality $\big(\sum_{i}\tau_iu_i\big)^p\leq \big(\sum_i\tau_i\big)^{p-1}\sum_i\tau_iu^p_i$, the property of the conditional expectation, the Burkholder--Davis--Gundy inequality, and the Lipschitz continuity of $G$.
Therefore, combining the H{\"o}lder inequality, we get for $p\ge 4,$
\begin{align*}
&\mathbb{E}[\|X^{N,K}_{t}\|^p]\leq C(p,T,\|Q\|_{\mathcal{L}_2(H)},\bar{L}_2)\left(1+\mathbb{E}[\|X_0\|^p]+\int_0^{t}\mathbb{E}[\|X^{N,K}_{t_{m_s}}\|^p]\mathrm{d}s\right)\\
&\leq C(p,T,\|Q\|_{\mathcal{L}_2(H)},\bar{L}_2)\left(1+\mathbb{E}[\|X_0\|^p]+\int_0^t\sup_{0\leq u\leq s}\mathbb{E}[\|X^{N,K}_u\|^p]\mathrm{d}s\right),
\end{align*}
which implies
\begin{align*}
\sup_{0\leq r\leq t}\mathbb{E}[\|X^{N,K}_r\|^p]\leq C(p,T,\|Q\|_{\mathcal{L}_2(H)},\bar{L}_2)\left(1+\mathbb{E}[\|X_0\|^p]+\int_0^t\sup_{0\leq u\leq s}\mathbb{E}[\|X^{N,K}_u\|^p]\mathrm{d}s\right).
\end{align*}
Applying the Gr\"onwall inequality leads to that for $p\ge 4,$
\begin{align}\label{gronwall}
\sup_{0\leq t\leq T}\mathbb{E}[\|X^{N,K}_{t}\|^p]\leq C(p,T,\|Q\|_{\mathcal{L}_2(H)},\bar{L}_2)(1+\mathbb{E}[\|X_0\|^p]).
\end{align}
%When $2\leq p<4,$ the H\"older's inequality gives that 
%$$\sup_{0\leq t\leq T}\mathbb{E}\|X^{N,K}_t\|^p\leq \sup_{0\leq t\leq T}(\mathbb{E}\|X^{N,K}_t\|^4)^{\frac{p}{4}}\leq C(p,T,\|Q\|_{\mathcal{L}_2(H)},\bar{L}_2)(1+\mathbb{E}\|X_0\|^4).$$
To obtain the expected supremum of the $p$th moment of $X^{N,K}_{t},$ we need to apply the following Burkholder--Davis--Gundy inequality
 \begin{align*}
&\mathbb{E}\Big[\sup_{0\leq t\leq T}\Big|\int_0^{t}\langle X^{N,K}_{t_{m_s}}+\Theta_K(X^{N,K}_{t_{m_s}})F(X^{N,K}_{t_{m_s}})\tau_{m_s},G(X^{N,K}_{t_{m_s}})\mathrm{d}W(s)\rangle\Big|^p\Big]\\
&\leq C\,\mathbb{E}\Big[\Big(\int_0^T\|X^{N,K}_{t_{m_s}}+\Theta_K(X^{N,K}_{t_{m_s}})F^N(X^{N,K}_{t_{m_s}})\tau_{m_s}\|^2\|G(X^{N,K}_{t_{m_s}})\|^2_{\mathcal{L}^0_2}\,\mathrm{d}s\Big)^{\frac{p}{2}}\Big]\\
&\leq C\int_0^T\sup_{s\in[0,T]}\mathbb{E}\big[\|X^{N,K}_{t_{m_s}}\|^p\big]\mathrm{d}s\leq C,
\end{align*}
where we use \eqref{gronwall}.
%Based on the above inequality and \eqref{gronwall}, 
And the remaining proof of $\mathbb{E}\big[\sup_{0\leq t\leq T}\|X^{N,K}_{t}\|^p\big]\leq C(1+\mathbb{E}[\|X_0\|^p])$
 can be given similarly. 
 Moreover, the assumption \eqref{adaptive5} gives that $\mathbb{E}[M_{T}]\leq T\mathbb{E}\Big[\sup_{0\leq t\leq T}\tau(X^{N,K}_{t})^{-1}\Big]\leq T\mathbb{E}[\sup_{0\leq t\leq T}(\zeta_1\|X^{N,K}_t\|^{q_0}+\zeta_2)]\leq C$ with $C$ independent of $K,$ which implies  that $T$ is a.s. attainable.

 \textit{Step 2:} Expected supremum of the $p$th moment of $X^N_t.$
 
 Define the stopping time by $\tilde \tau_K=\inf \{t\in[0,T]:\|X^{N,K}_t\|\ge K\}.$
 Then for $t\in [0,\tilde \tau_{K_1}\wedge \tilde \tau_{K_2}),$ we have $\Theta_{K_1}(X^{N,K_1}_t)=\Theta_{K_2}(X^{N,K_2}_t)=1$ and thus $X^{N,K_1}_t=X^{N,K_2}_t$ a.s. due to the existence and uniqueness of the solution. Hence, $\tilde \tau_K$ is nondecreasing with respect to $K.$ Denote $\lim_{K\to\infty}\tilde \tau_K=T_*$ a.s.
 The Chebyshev inequality gives 
 \begin{align*}
 &\mathbb{P}\{\tilde \tau_K=T\}=\mathbb{P}\Big\{\sup_{0\leq t\leq T}\|X^{N,K}_t\|\leq K\Big\}\\
 &\ge 1-K^{-4}\mathbb{E}\Big[\sup_{0\leq t\leq T}\|X^{N,K}_t\|^4\Big]\rightarrow 1\;\;\text{ as }K\to \infty,
 \end{align*}
  which implies $T=T_*$. Define $X^N_t$ on $[0,T]$ by $X^N_t=X^{N,K}_t$ on $[0,\tilde \tau_K).$ Then $X^N_t$ is the solution of \eqref{conti_multi} with $ \lim_{K\to\infty}\sup_{t\in[0,\tilde\tau_K)}\|X^{N,K}_t\|=\sup_{t\in [0,T]}\|X^N_t\|$ a.s. 
 Hence, applying the Fatou lemma and \textit{Step 1} leads to
 \begin{align*}
 \mathbb{E}\Big[\sup_{0\leq t\leq T}\|X^N_t\|^p\Big]\leq \lim_{K\to\infty}\mathbb{E}\Big[\sup_{0\leq t\leq T}\|X^{N,K}_t\|^p\Big]\leq C(1+\mathbb{E}[\|X_0\|^p]).
 \end{align*}
 The proof is finished.
\end{proof}

Below, we give the \textit{a priori} estimate in $E$-norm of the solution of \eqref{scheme_multi} for the case of $d=1$ based on the Sobolev embedding $H^{\frac12+}\hookrightarrow E.$
\begin{prop}\label{prop6}
Let $d=1.$ Under Assumptions \ref{assumption1}, \ref{assumption3}, \ref{assumption6}-\ref{assumpP^NX_0}, and \eqref{adaptive4}-\eqref{adaptive3}, we have for $p\ge 2,$
\begin{align*}
\sup_{N\in\mathbb{N}_+}\sup_{0\leq t\leq T}\|X^N_t\|_{L^p(\Omega;E)}\leq C
\end{align*}
for some $C:=C(p,T,X_0)>0,$ where $X^N_t$ is the solution of \eqref{scheme_multi} with timestep function $\tau^{\delta}$ for $\delta\in(0,1)$. 
\end{prop}
\begin{proof}
We split the numerical solution as $X^N_t=Y^N_{1,t}+Y^N_{2,t}+\mathbb{Z}^N_t,$ where $Y^N_{1,t},Y^N_{2,t},$ and $\mathbb{Z}^N_t$ are solutions of
\begin{align}\label{split_Y}
\mathrm{d}Y^N_{1,t}=-A^NY^N_{1,t}\mathrm{d}t+S^N(t-t_{m_t})F^N(X^N_{t_{m_t}})\mathrm{d}t-F^N(X^N_t)\mathrm{d}t,\quad Y^N_{1,0}=P^NX_0,
\end{align}
\begin{align*}
\mathrm{d}Y^N_{2,t}=F^N(X^N_t)\mathrm{d}t,\quad Y^N_{2,0}=0,
\end{align*}
and 
\begin{align}\label{split_Z}
\mathrm{d}\mathbb{Z}^N_t=-A^N\mathbb{Z}^N_t\mathrm{d}t+S^N(t-t_{m_t})P^NG(X^N_{t_{m_t}})\mathrm{d}W_t,\quad \mathbb{Z}^N_0=0,
\end{align}
respectively. 

%Introduce  processes $Y^N_{1,t}$ and $Y^N_{2,t}$ satisfying 
%\begin{align*}
%Y^N_{1,t}:=S^N(t)P^NX_0+\int_0^tS^N(t-t_{m_s})F^N(X^N_{t_{m_s}})\mathrm{d}s-\int_0^tS^N(t-s)F^N(X^N_s)\mathrm{d}s
%\end{align*}
%and 
%\begin{align*}
%Y^N_{2,t}=\int_0^tS^N(t-s)F^N(X^N_s)\mathrm{d}s,
%\end{align*}
%respectively. 
%Then $X^N_t=Y^N_{1,t}+Y^N_{2,t}+\mathbb{Z}^N_t.$ 
Applying Lemma \ref{purt_lemma} to  $Y^N_{2,t}$ yields that
\begin{align*}
\|Y^N_{2,t}\|_E\leq C\big(1+\int_0^t(t-s)^{-\frac d4}\|Y^N_{1,s}+\mathbb{Z}^N_s\|^3_E+\|Y^N_{1,s}+\mathbb{Z}^N_s\|^{72}_E\mathrm{d}s\big).
\end{align*}
Recall the definition $Z^N_t$ (see \eqref{def_Z}) in the proof of Proposition \ref{thm2}. It can be observed that $Y^N_{1,t}+\mathbb{Z}^N_t$ has a similar expression to that of $Z^N_t,$ in which the counterpart of $W^N_A$ is $\mathbb{Z}^N_t$ in this setting. 
%Note that $\mathbb{E}\big[\sup_{0\leq t\leq T}\|W^N_A(t)\|_E^q\big]<\infty$ for some large $q>2$ plays an important role in the proof of Proposition \ref{thm2}, and in fact $\sup_{0\leq t\leq T}\mathbb{E}\big[\|W^N_A(t)\|_E^q\big]<\infty$ is enough for the proof. 
 Thus it suffices to prove that 
$\sup_{0\leq t\leq T}\mathbb{E}\big[\|\mathbb{Z}^N_t\|_E^q\big]<\infty
$ for $q\ge 2$ in the multiplicative noise case.
In fact,  the Sobolev embedding $H^{\frac12+\epsilon}\hookrightarrow E$ with $\epsilon>0$,  the Burkholder--Davis--Gundy inequality and the Lipschitz continuity of $G$ give that
\begin{align*}
\sup_{0\leq t\leq T}\mathbb{E}\Big[\|\mathbb{Z}^N_t\|^q_E\Big]&\leq C\sup_{0\leq t\leq T}\mathbb{E}\Big[\Big\|\int_0^tS^N(t-t_{m_s})P^NG(X^N_{t_{m_s}})\mathrm{d}W(s)\Big\|_{\frac12+\epsilon}^q\Big]\\
&\leq C\sup_{0\leq t\leq T}\mathbb{E}\Big[\Big(\int_0^t\|A^{\frac{1+2\epsilon}{4}}S^N(t-t_{m_s})P^NG(X^N_{t_{m_s}})\|^2_{\mathcal{L}^0_2}\mathrm{d}s\Big)^{\frac q2}\Big]\\
&\leq C\sup_{0\leq t\leq T}\mathbb{E}\Big[\Big(\int_0^t(t-t_{m_s})^{-\frac{1+2\epsilon}{2}}\|P^NG(X^N_{t_{m_s}})\|^2_{\mathcal{L}^0_2}\mathrm{d}s\Big)^{\frac q2}\Big]\\
&\leq C\sup_{0\leq t\leq T}\mathbb{E}\Big[\Big(\int_0^t(t-s)^{-\frac{1+2\epsilon}{2}}(1+\|X^N_{t_{m_s}}\|^2)\mathrm{d}s\Big)^{\frac q2}\Big]\\
&\leq C\sup_{0\leq t\leq T}\Big(\int_0^t(t-s)^{-\frac{1+2\epsilon}{2}}\big(1+\sup_{0\leq s\leq T}\|X^N_{t_{m_s}}\|^2_{L^q(\Omega;H)}\big)\mathrm{d}s\Big)^{\frac{q}{2}}\leq C.
\end{align*} 
The remaining proof is similar to that of Proposition \ref{thm2} under conditions \eqref{adaptive4}-\eqref{adaptive5} and hence is omitted.  
\end{proof}
With the above regularity estimates, one can derive the \textit{a priori} estimate in $\dot{H}^1$-norm of the solution of \eqref{scheme_multi}  by a similar approach to Proposition \ref{coro}. 
Then similar to the proof of Theorem \ref{thm3}, one can obtain the following convergence order of the adaptive time-stepping scheme \eqref{scheme_multi}. %whose  proof is omitted.  
\begin{prop}
Let  $d=1$. Under conditions in Proposition \ref{prop6}, for $p\ge 4$, 
\begin{align*}
\sup_{t\in[0,T]}\|X(t)-X^N_t\|_{L^p(\Omega;H)}\leq C(\lambda_N^{-\frac12}+\delta^{\frac 12}),
\end{align*}
where $X^N_t$ is the solution of \eqref{scheme_multi} and $C:=C(p,T,X_0)>0.$
\end{prop}

\begin{proof}
Similar to the proof of Theorem \ref{thm3}, we introduce the auxiliary process $\mathbb Y^N_t,$
\begin{align*}
\mathbb Y^N_t=S^N(t)P^NX_0+\int_0^tS^N(t-s)F^N(X(s))\mathrm ds+\int_0^tS^N(t-t_{m_s})P^NG(\mathbb Y^N_{s})\mathrm dW(s).
\end{align*}
It can be shown that $\|\mathbb Y^N_t\|_{L^p(\Omega;H)}\leq C$ and $\|\mathbb Y^N_t-\mathbb Y^N_s\|_{L^p(\Omega;H)}\leq C(t-s)^{\frac12}.$
The error can be divided into 
\begin{align*}
\|X(t)-X^N_t\|_{L^p(\Omega;H)}\leq \|X(t)-P^NX(t)\|_{L^p(\Omega;H)}+\|P^NX(t)-\mathbb Y^N_t\|_{L^p(\Omega;H)}+\|\mathbb Y^N_t-X^N_t\|_{L^p(\Omega;H)}.
\end{align*}
With the regularity of the solution $X$, we have   $\|X(t)-P^NX(t)\|_{L^p(\Omega;H)}\leq C\lambda_N^{-\frac{1}{2}}\|X(t)\|_{L^p(\Omega;\dot{H}^{1})}$. 

For the term $\|P^NX(t)-\mathbb Y^N_t\|_{L^p(\Omega;H)},$ it follows from the Burkholder--Davis--Gundy inequality that for $p\ge 2,$
\begin{align*}
 &\quad \|P^NX(t)-\mathbb Y^N_t\|_{L^p(\Omega;H)}\\
 &\leq \Big\|\int_0^tS^N(t-s)(G(X(s))-G(\mathbb Y^N_s))\mathrm dW(s)\Big\|_{L^p(\Omega;H)}\\
 &\quad +\Big\|\int_0^tS^N(t-s)(\mathrm{Id}-S^N(s-t_{m_s}))G(\mathbb Y^N_s)\mathrm dW(s)\Big\|_{L^p(\Omega;H)}\\
 &\leq C\Big[\Big\|\int_0^t\|G(X(s))-G(P^NX(s))\|^2_{\mathcal L^0_2}\mathrm ds\Big\|_{L^{\frac p2}(\Omega)}^{\frac12}+\Big\|\int_0^t\|G(P^NX(s))-G(\mathbb Y^N_s)\|^2_{\mathcal L^0_2}\mathrm ds\Big\|_{L^{\frac p2}(\Omega)}^{\frac12}\\
 &\quad +\Big\|\int_0^t\|S^N(t-s)(\mathrm{Id}-S^N(s-t_{m_s}))G(\mathbb Y^N_s)\|^2_{\mathcal L^0_2}\mathrm ds\Big\|^{\frac12}_{L^{\frac p2}(\Omega)}\Big]\\
 &\leq C\Big[\lambda_N^{-\frac12}+\Big(\int_0^t\|P^NX(s)-\mathbb Y^N_s\|^2_{L^p(\Omega;H)}\mathrm ds\Big)^{\frac12}\\
 &\quad +\Big\|\int_0^t\|A^{\frac12}S^N(t-s)A^{-\frac12}(\mathrm{Id}-S^N(s-t_{m_s}))G(\mathbb Y^N_t)\|^2_{\mathcal L^0_2}\mathrm ds\Big\|^{\frac12}_{L^{\frac p2}(\Omega)}\\
 &\quad +\Big\|\int_0^t\|A^{\frac12}S^N(t-s)A^{-\frac12}(\mathrm{Id}-S^N(s-t_{m_s}))(G(\mathbb Y^N_t)-G(\mathbb Y^N_s))\|^2_{\mathcal L^0_2}\mathrm ds\Big\|^{\frac12}_{L^{\frac p2}(\Omega)}\Big],
\end{align*} 
which together with \eqref{S(t)1}-\eqref{S(t)3},  the H\"older continuity of $\mathbb Y^N_t$ and the Gr\"onwall inequality yields that $\|P^NX(t)-\mathbb Y^N_t\|_{L^p(\Omega;H)}\leq C(\lambda_N^{-\frac12}+\delta^{\frac12}).$

For the term $\|\mathbb Y^N_t-X^N_t\|_{L^p(\Omega;H)},$ 
applying the It\^o formula gives 
\begin{align}\label{mul_X-Y}
\|X^N_t-\mathbb Y^N_t\|^2&=2\int_0^t\langle X^N_s-\mathbb Y^N_s,-A^N(X^N_s-\mathbb Y^N_s)\rangle \mathrm dt\notag\\
&\quad+2\int_0^t\Big\langle X^N_s-\mathbb Y^N_s,S^N(s-t_{m_s})F^N(X^N_{t_{m_s}})-F^N(X(s))\Big\rangle\mathrm ds\notag\\
& +2\int_0^t\Big\langle X^N_s-\mathbb Y^N_s,S^N(s-t_{m_s})P^N(G(X^N_{t_{m_s}})-G(\mathbb Y^N_s))\mathrm dW(s)\Big\rangle\notag\\
& +\int_0^t\|S^N(s-t_{m_s})P^N(G(X^N_{t_{m_s}})-G(\mathbb Y^N_s))\|^2_{\mathcal L^0_2}\mathrm ds. 
\end{align}
Then compared with the additive noise case (see \eqref{add_X-Y}), the main difference of the estimation of $\|X^N_t-\mathbb Y^N_t\|_{L^p(\Omega;H)}$  lies in the estimation of the last two terms in \eqref{mul_X-Y}, which can be estimated as
 for $p\ge 4,$\begin{align*}
&\quad  \mathbb E \Big[\Big|\int_0^t\Big\langle X^N_s-\mathbb Y^N_s,S^N(s-t_{m_s})P^N(G(X^N_{t_{m_s}})-G(\mathbb Y^N_s))\mathrm dW(s)\Big\rangle\Big|^{\frac p2}\Big]\\
 &\leq C \mathbb E\Big[\Big(\int_0^t\|X^N_s-\mathbb Y^N_s\|^2\|G(X^N_{t_{m_s}})-G(\mathbb Y^N_s)\|^2_{\mathcal L^0_2}\mathrm ds\Big)^{\frac p4}\Big]\\
 &\leq C\int_0^t\mathbb E[\|X^N_s-\mathbb Y^N_s\|^{p}]\mathrm ds+C\delta^{\frac p2},
 \end{align*}
 and
 \begin{align*}
 &\quad \mathbb E\Big[\Big|\int_0^t\|S^N(t-t_{m_s})P^N(G(X^N_{t_{m_s}})-G(\mathbb Y^N_s))\|^2_{\mathcal L^0_2}\mathrm ds\Big|^{\frac p2}\Big]\\
 &\leq C\int_0^t\mathbb E[\|X^N_s-\mathbb Y^N_s\|^{p}]\mathrm ds+C\delta^{\frac p2},
 \end{align*}
 respectively, where the H\"older continuity of $X^N_t$, the H\"older inequality, the Young inequality and the Lipschitz condition of $G$ are used.
 Hence, similar to the proof of \eqref{add_X-Y}, one can obtain that $\|X^N_t-\mathbb Y^N_t\|_{L^p(\Omega;H)}\leq C\lambda_N^{-\frac12}+\delta^{\frac12}.$ This finishes the proof. 
\end{proof}

%With the above regularity estimates, similar to the proof Theorem \ref{thm3}, the convergence order for  \eqref{scheme_multi} can be proved to be $1/2$ in time and $1$ in space for $d=1$. 
Note that the proof of Proposition \ref{prop6}  is not applicable for the case of $d>1$ due to that the Sobolev embedding $H^{\frac d2+}\hookrightarrow E$ requires a higher space regularity. %{\color{blue}This can be solved by letting $\|A^{\frac r2}G(X)\|_{\mathcal L^0_2}\leq C\|X\|$ with some $r>\frac d2-1$ for $d>1.$ Otherwise,} 
Hence other skills  
%such as studying \eqref{multi_spde} and \eqref{scheme_multi} in M-type 2 Banach spaces, 
are needed to obtain the regularity estimate of \eqref{scheme_multi} in $E$-norm. We leave this as the future work and attempt to study it in the future.

\section{Numerical experiments}\label{sec3}
In this section, we present numerical experiments for the trace-class noise case and the space-time white noise case to verify the previous theoretical results, respectively. Meanwhile, some alternative choices of the adaptive timestep function are given.

Consider the following  stochastic Allen--Cahn equation with the generalized Q-Wiener process:
\begin{align}
\begin{cases}
\frac{\partial u}{\partial t}=\frac{\partial^2 u}{\partial x^2}+u-u^3+\dot{W}(t,x),\quad t\in(0,1],\;x\in(0,1),\\
u(0,x)=\sqrt{2}\sin(\pi x),\quad x\in(0,1),\\
u(t,0)=u(t,1)=0,\quad t\in(0,1].
\end{cases}
\end{align} 
For the trace-class noise case, we choose $Q$ such that $Qe_i=\frac{1}{i^2}e_i,\; i\ge 1$, which implies that $tr(Q)<\infty$ and Assumption \ref{assumption2} holds for some $\beta\in[1,2]$. For the space-time white noise case, i.e., $Q=\mathrm{Id}$, Assumption \ref{assumption2} holds for $\beta\in(0,\frac{1}{2})$. 

We are going to compare the performance of adaptive schemes with the scheme with a uniform timestep. The scheme with a uniform timestep is chosen as the tamed exponential integrator (TE), see e.g. \cite{WXJ20}. 
In the following, error bounds are measured in root-mean-square (RMS) sense at the end point $T=1$, caused by spatio-temporal discretization. And the expectations are approximated by computing averages over $1000$ samples. The infinite-dimensional Hilbert space $L^2(0,1)$ is approximated through the finite-dimensional subspace spanned by the first $2^{10}$ eigenfunctions of the Laplacian, i.e., in the spectral Galerkin method, we take $N=2^{10}.$ Since the exact solution cannot be given explicitly, we take the solution generated by the same method with a timestep which is three times smaller as the reference solution. Numerical experiments are tested by Matlab R2017a in MacBook Pro (13-inch, 2019, Two Thunderbolt 3 ports).

Besides \eqref{Psi1}, there are other choices for the backstop scheme $\Psi$, for example,  the nonlinearity-truncated exponential integrator (see \cite[Eq. (6)]{BJ19})
\begin{align}
\Psi(x,h,y)&=S^N(h)\Big(x+hF^N(x)\chi_{\{\|x\|^6_{L^{18}(\mathcal{O})}\leq h^{-1}\}}+y\Big),\label{Psi2}
\end{align}
and the linear-implicit nonlinearity-truncated scheme (see \cite[Eq. (7)]{BJ19})
\begin{align}
\Psi(x,h,y)&=(\mathrm{Id}-A^Nh)^{-1}\Big(x+hF^N(x)\chi_{\{\|x\|^6_{L^{18}(\mathcal{O})}\leq h^{-1}\}}+y\Big).\label{Psi3}
\end{align}
Hence, for the coupled scheme \eqref{expression_cou}, when the backstop scheme $\Psi$ is chosen to be \eqref{Psi2} or \eqref{Psi3}, we get two coupled schemes denoted by (CAU 2) and  (CAU 3), respectively. Similarly, for the coupled scheme \eqref{expression_cou2}, when the backstop scheme $\Psi$ is chosen to be \eqref{Psi1}, \eqref{Psi2} or \eqref{Psi3}, we get three coupled schemes denoted by (CAU a), (CAU b) and (CAU c), respectively.

\subsection{Trace-class noise}
We show the experiment results for the trace-class noise case in this subsection. 
For the trace-class noise case,  the seven methods tested are: the scheme TE which is with the  uniform timestep, schemes \eqref{ATEU}, (CAU 2), (CAU 3), (CAU a), (CAU b) and (CAU c) which are with the adaptive timestep. Recall that $\tau^{\delta}$ is the refined timestep function controlled by the scalar parameter $\delta\in(0,1)$ and satisfies Assumption \ref{assumption6}. 
We choose the following two types of adaptive timestep functions: \\
\textit{type 1:}
$$
 \tau^{\delta}_{11}(X)=\delta\tau_{11}(X) \text{ with }\tau_{11}(X)=\left(\frac{\|X\|+0.5}{2\|F(X)\|+1.8}\right)^{\frac{4}{3}},
 $$
 $$
 \tau^{\delta}_{12}(X)=\delta \tau_{12}(X)\text{ with }\tau_{12}(X)=\min\left\{\frac{2\|X\|^4_{L^4}+0.01}{\|X\|^6_{L^6}+0.01},\left(\frac{\|X\|+0.2}{\|F(X)\|+0.8}\right)^{\frac{4}{3}}\right\},
 $$
 \textit{type 2:}
$$
 \tau^{\delta}_{21}(X)=\delta \tau_{21}(X)\text{ with }\tau_{21}(X)=\left(\frac{1.1}{\|F(X)\|+3.2}\right)^{\frac{4}{3}},
 $$
 $$
 \tau^{\delta}_{22}(X)=\delta\tau_{22}(X)\text{ with }\tau_{22}(X)= \min\left\{\frac{\|X\|^2+0.01}{\|X\|^6_{L^6}+0.04},\left(\frac{1.1}{\|F(X)\|+3.2}\right)^{\frac{4}{3}}\right\},
$$
where $\delta=2^{-l},l=2,\ldots,7$. 
% $\phi_{1,1}=\phi_{1,2}=1,\phi_{2,1}=3,\phi_{2,2}=1$. 
Here, $\tau_{11},\tau_{21}$ are adaptive timestep functions for \eqref{ATEU}, (CAU 2) and (CAU 3), and $\tau_{12}, \tau_{22}$ are adaptive timestep functions for (CAU a), (CAU b) and (CAU c). %{\color{blue}Constants in \textit{type $i$ ($i=1,2$)} timestep functions are flexible.} 
In \textit{type 1}, we let $L_2(\omega)=\|X(\omega)\|+c$ with $X$ being the reference solution or numerical solution, $c=0.5$ in $\tau^{\delta}_{1,1}$ and $c=0.2$ in $\tau^{\delta}_{1,2}$,  and $\bar{L}_2=\frac{3}{2}$. 
In \textit{type 2}, we let $L_2=1.1$
and $\bar{L}_2=2$ 
independent of $\omega.$
One can verify that the \textit{type 1} and \textit{type 2} timestep functions satisfy 
 Assumption \ref{assumption6}.

\begin{figure}[tbhp]
\centering

\subfloat[timestep type 1]{
\begin{minipage}[t]{0.4\linewidth}
\centering
\includegraphics[height=5cm,width=6.8cm]{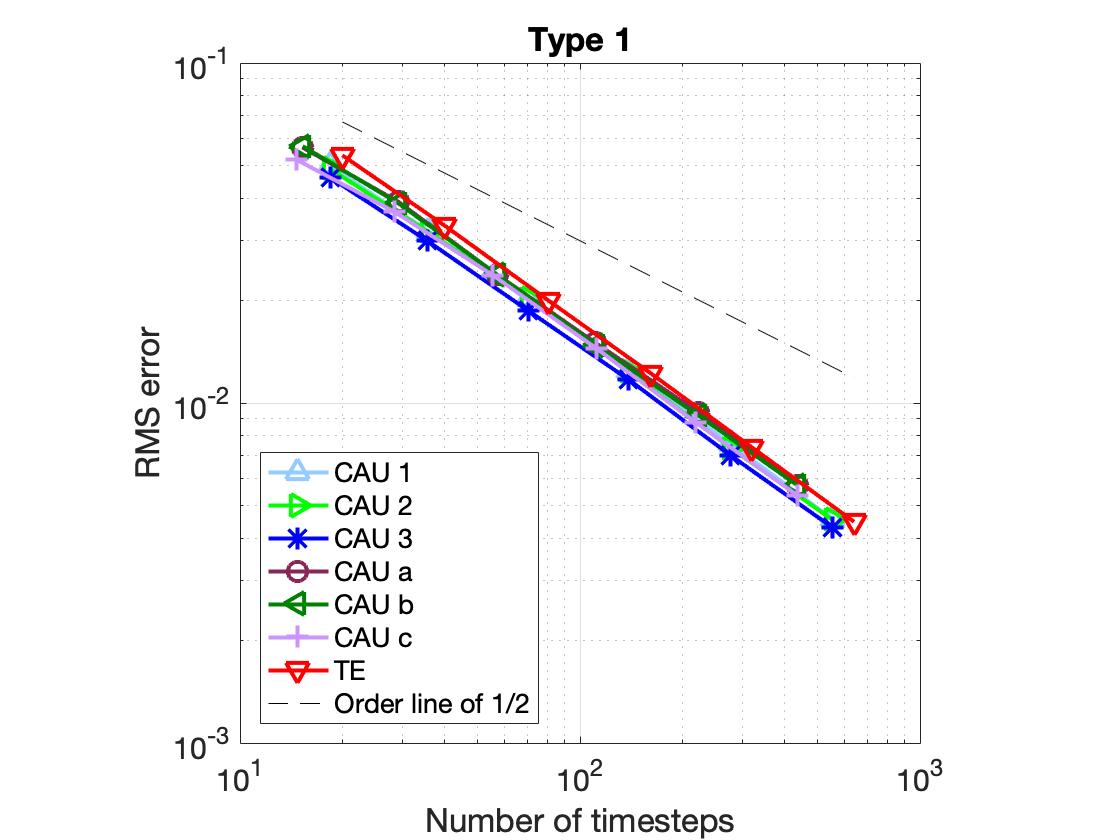}
\end{minipage}
}
\subfloat[timestep type 1]{
\begin{minipage}[t]{0.4\linewidth}
\centering
\includegraphics[height=5cm,width=6.8cm]{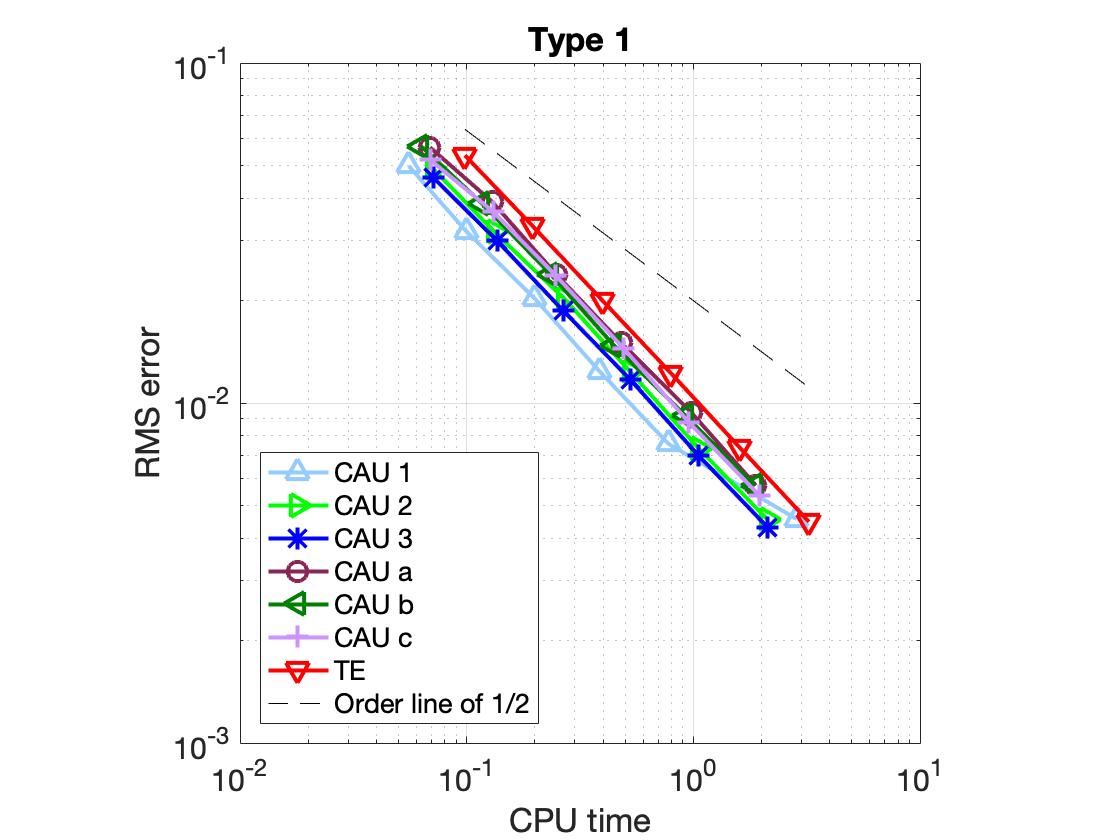}
\end{minipage}
}

\centering
\caption{RMS error for numerical schemes with trace-class noise for type $1$}
\label{RMStrace2}
\end{figure}

The comparisons of the seven schemes with \textit{type 1} and \textit{type 2} timestep functions are presented in Figure \ref{RMStrace2} and Figure \ref{RMStrace3}, respectively, where the left ones are about the RMS error against the average number of timesteps, and the right ones are about the RMS error against the CPU time. 
 Recall  \eqref{adaptive2} and \eqref{adaptive5}, and set $\tau_{min}=0.2\delta,$ $\zeta_1=3,q_0=1$ and $\zeta_2=4$ for \textit{type 1} and \textit{type 2}. The length of a timestep for TE is set by $\tau_{min}.$
From  Figures \ref{RMStrace2} and \ref{RMStrace3}, we observe that the convergence orders in the temporal direction of these seven schemes are slightly bigger than $\frac{1}{2}$. As for the \textit{type 1} timestep function, for a given RMS error, we see from Figure \ref{RMStrace2} $(a)$ that the adaptive schemes cost slightly fewer numbers of timesteps than TE which has the uniform timestep. And we observe from Figure \ref{RMStrace2} $(b)$ that the adaptive schemes cost less CPU time than TE for a given RMS error. A similar phenomenon is observed from Figure \ref{RMStrace3} for the \textit{type 2} timestep function.
These mean that for the \textit{type 1} and \textit{type 2} timestep functions, the adaptive schemes perform slightly better than TE which uses the uniform timestep.

\begin{figure}[tbhp]
\centering

\subfloat[timestep type 2]{
\begin{minipage}[t]{0.4\linewidth}
\centering
\includegraphics[height=5cm,width=6.8cm]{picture/%o_1101_trace22_step.jpg}
step20507}
%step2_0508modi}
\end{minipage}
}
\subfloat[timestep type 2]{
\begin{minipage}[t]{0.4\linewidth}
\centering
\includegraphics[height=5cm,width=6.8cm]{picture/%o_1101_trace22_cpu.jpg}
cpu20507}
%cpu2_0508modi}
\end{minipage}
}

\centering
\caption{RMS error for numerical schemes with trace-class noise for type $2$}
\label{RMStrace3}
\end{figure}

\begin{table}[tbhp]
\begin{tabular}{ccccccccc}
\hline
 & CAU 1 & CAU 2& CAU 3& CAU a &CAU b &CAU c&TE \\
\hline
Type 1 &81.78&82.80& 82.83&77.66&77.00&86.08&93.52\\
\hline
Type 2 &76.84&82.74&79.54&87.29&88.75&86.55&93.52\\
\hline
\end{tabular}
\vspace{0.1in}

\centering
\caption{Time in seconds ($\delta=2^{-6}$)}
\label{table1}
\end{table}

\begin{figure}[tbhp]
\centering
\subfloat{
\begin{minipage}[t]{0.4\linewidth}
\centering
\includegraphics[height=5cm,width=6.8cm]{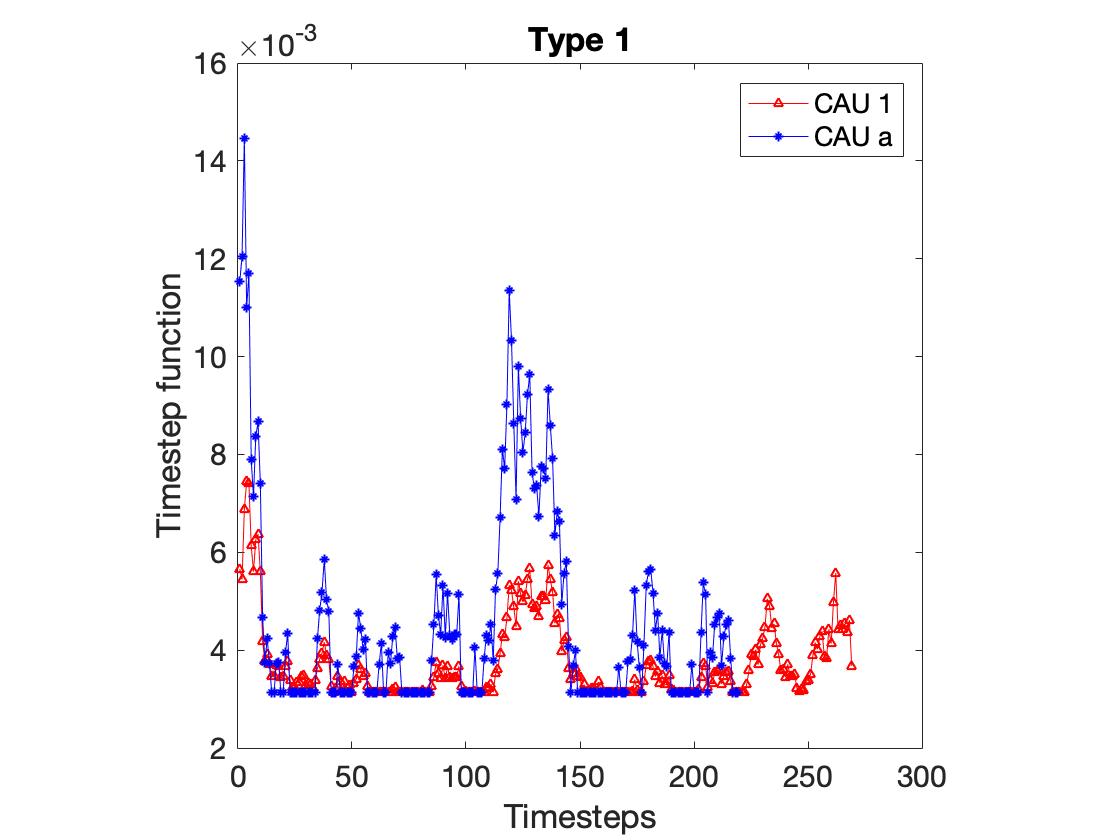}%o_1101_stepfunction1.jpg}
\end{minipage}
}
\subfloat{
\begin{minipage}[t]{0.4\linewidth}
\centering
\includegraphics[height=5cm,width=6.8cm]{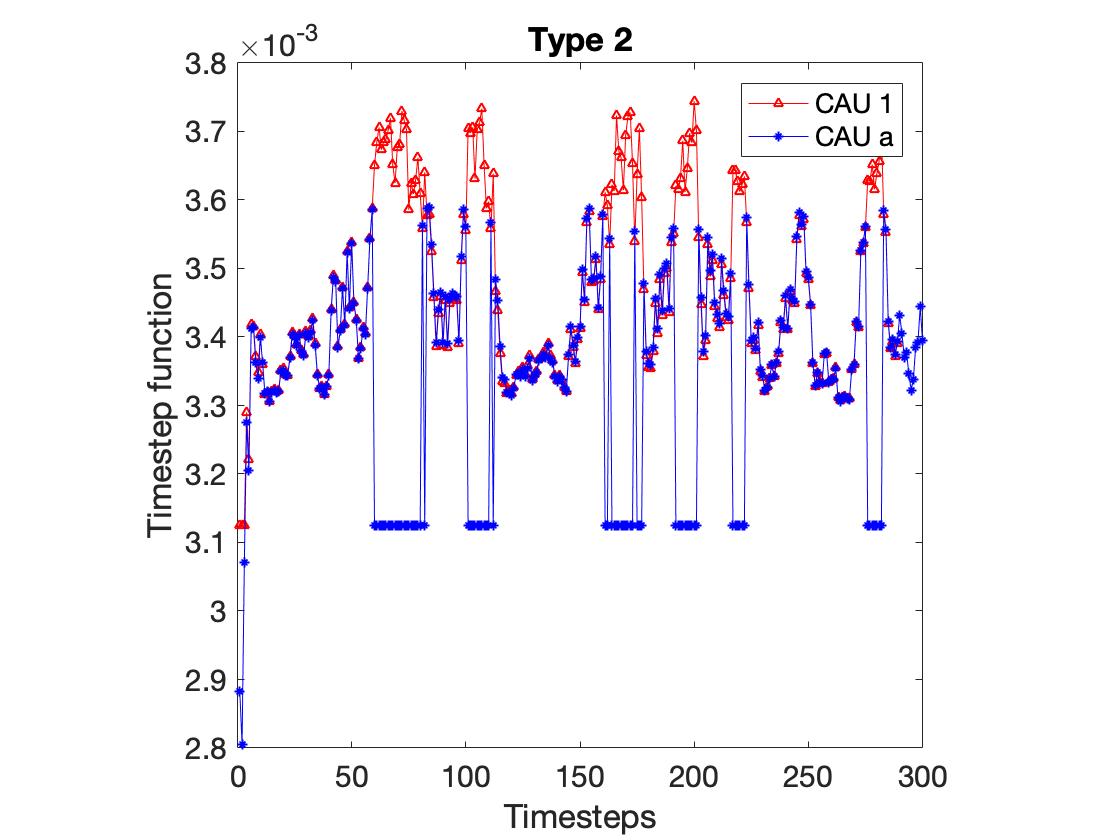}%o_1101_stepfunction2.jpg}
\end{minipage}
}
\caption{Timestep function for numerical schemes with trass-class noise ($\delta=2^{-6}$)}
\label{Timesteptrace1}
\end{figure}

We list the CPU time costed by these schemes with the  \textit{type $i$} ($i=1,2,3$) timestep functions  in Table \ref{table1}.  Here, the uniform timestep is $\tau_{min}=0.2\delta$ with $\delta=2^{-6}$. The critical parameter  for (CAU 1)-(CAU 3) is $\tau_{min},$ and that for (CAU a)-(CAU c) is $(\zeta_1\|X\|^{q_0}+\zeta_2)^{-1}$ with $\zeta_1=3,\zeta_2=4,q_0=1.$ And 1000 realizations are calculated. It can be observed that the coupled schemes are in a lower computational cost in terms of the CPU time compared with  TE that uses the uniform timestep. 
This  can be explained by Figure \ref{Timesteptrace1}, since coupled schemes like (CAU 1) and (CAU c) use many large timesteps compared with the uniform timestep, which may reduce the computational cost of the schemes.

\subsection{Space-time white noise}
We show the experiment results for the space-time white noise case in this subsection. 
For the space-time white noise case, the four schemes tested are: the scheme TE which is with the  uniform timestep, schemes \eqref{ATEU}, (CAU 2) and (CAU 3) which are with the adaptive timestep. We choose the following adaptive timestep function for the considered three adaptive schemes:\\
\textit{type 3:}
$$
 \tau^{\delta}_{31}(X)=\delta \tau_{31}(X)\text{ with }\tau_{31}(X)=\left(\frac{1.1}{\|F(X)\|+3}\right)^{\frac{4}{3}},
$$
where $\delta$ is the same as before.
In this case, the \textit{type 3} timestep function can also be verified to satisfy   Assumption \ref{assumption6}.

\begin{figure}[tbhp]
\centering

\subfloat[timestep type 3]{
\begin{minipage}[t]{0.4\linewidth}
\centering
\includegraphics[height=5cm,width=6.8cm]{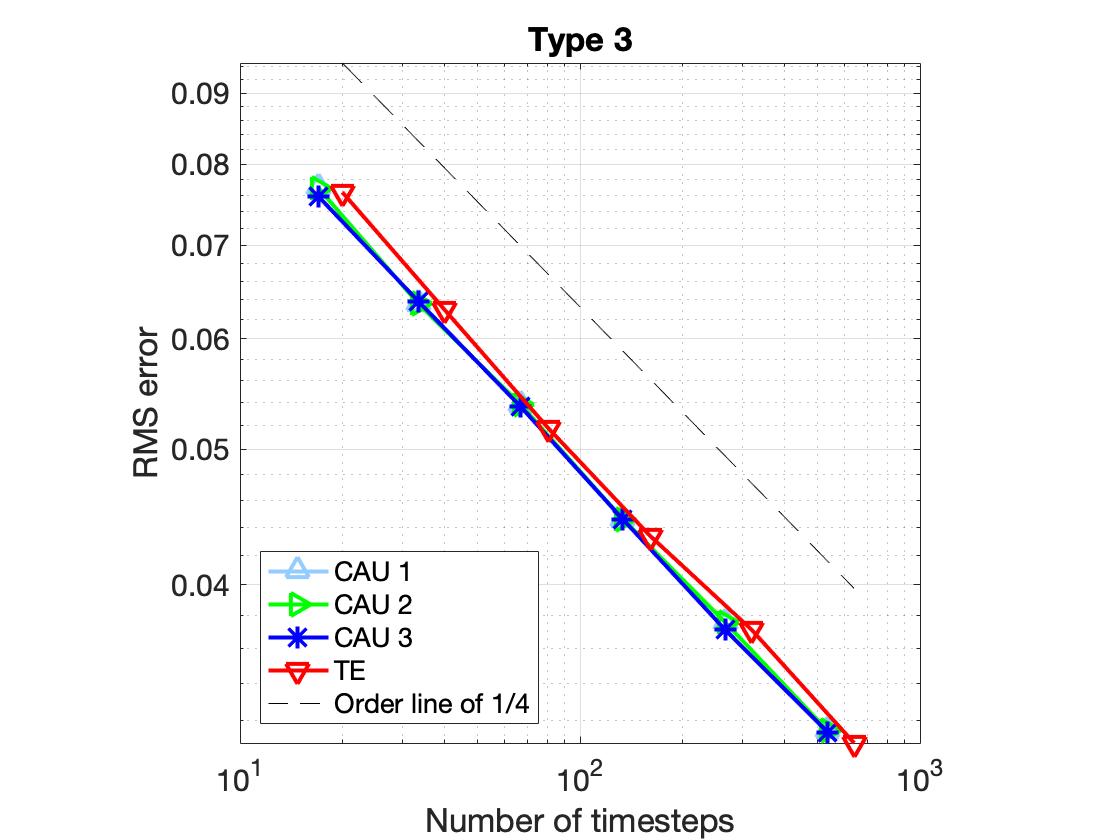}%o_1030_white_step3.jpg}
\end{minipage}
}
\subfloat[timestep type 3]{
\begin{minipage}[t]{0.4\linewidth}
\centering
\includegraphics[height=5cm,width=6.8cm]{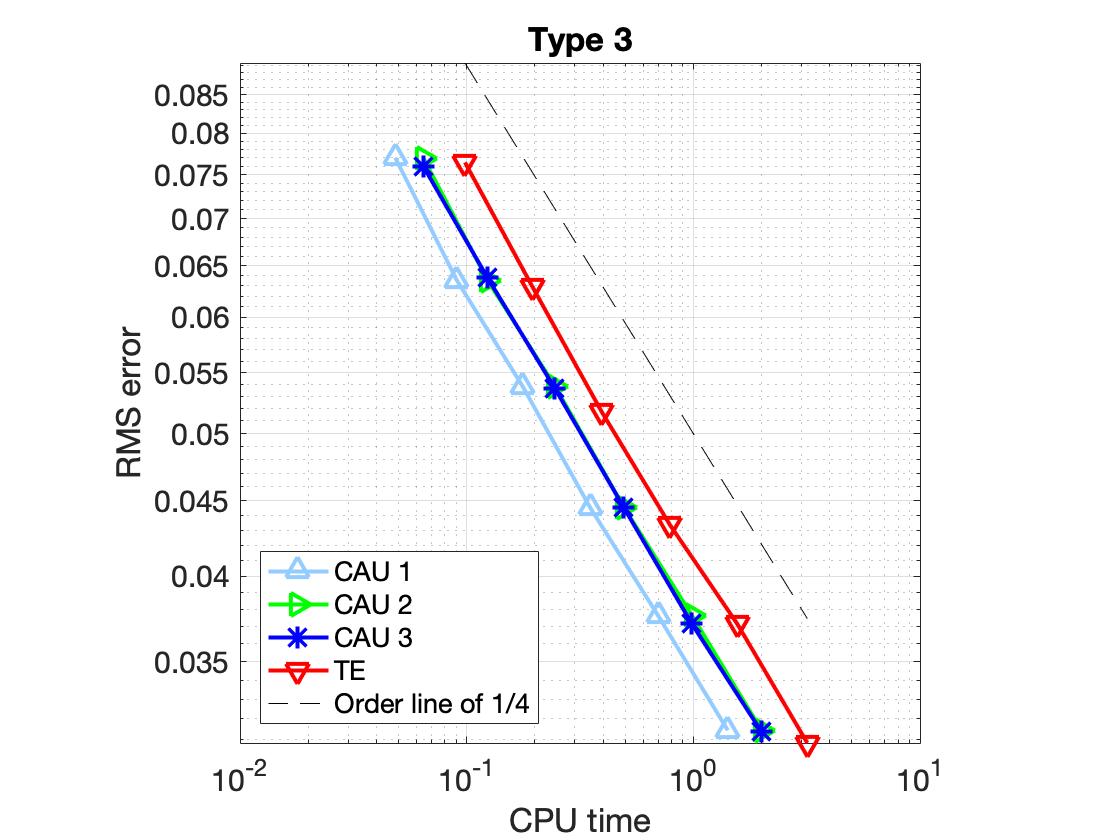}%o_1030_white_cpu3.jpg}
\end{minipage}
}

\centering
\caption{RMS error for numerical schemes with space-time white noise}
\label{RMSwhite}
\end{figure}

For space-time white noise case, Figure \ref{RMSwhite} shows the comparison of RMS error against the number of timesteps and CPU time. 
%The green, blue, red and black solid lines represent \eqref{ATEU}, (CAU 2), (CAU 3) and TE, respectively, and the black dotted line is the reference line of order $\frac{1}{4}$. 
It can be gained from Figure \ref{RMSwhite} that the convergence orders in the temporal direction of the four schemes are $\frac{1}{4}$. 
For a given RMS error, we see from Figure \ref{RMSwhite} $(a)$ that the adaptive schemes cost slightly fewer numbers of timesteps than TE which uses the uniform timestep. And we observe from Figure \ref{RMSwhite} $(b)$ that the adaptive schemes cost less CPU time than TE for a given RMS error. The total time for these schemes with \textit{type 3} timestep function is listed in Table \ref{table2}, from which we can see that the coupled schemes are in a lower computational cost in terms of the CPU time.   
These mean that for the \textit{type 3} timestep function, the adaptive schemes still perform slightly better than TE which uses the uniform timestep. 

\begin{table}[tbhp]
\begin{tabular}{ccccc}
\hline
 & CAU 1 & CAU 2& CAU 3 &TE \\
\hline
Type 3 &82.80&81.98& 77.85&102.62\\
\hline
\end{tabular}
\vspace{0.1in}
\caption{Time in seconds ($\delta=2^{-6}$)}
\label{table2}
\end{table}

\begin{figure}[tbhp]
\centering

\subfloat[timestep type 1]{
\begin{minipage}[t]{0.4\linewidth}
\centering
\includegraphics[height=5cm,width=6.8cm]{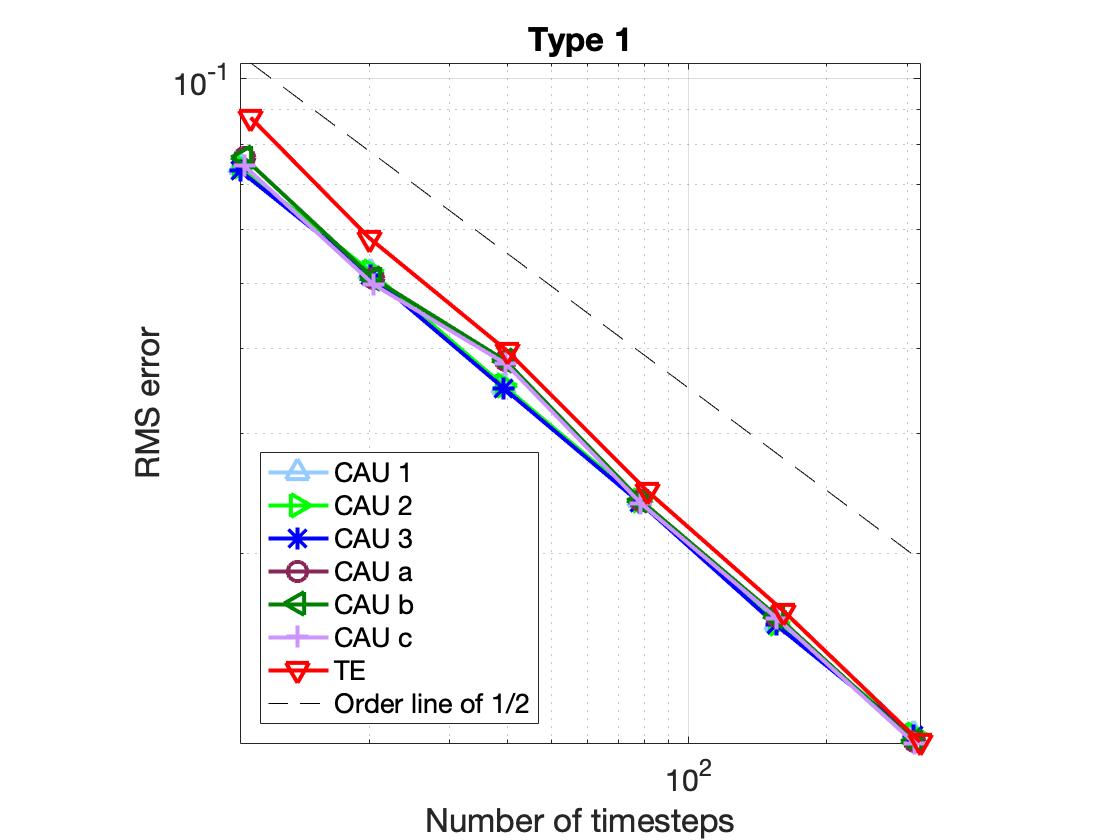}
\end{minipage}
}
\subfloat[timestep type 1]{
\begin{minipage}[t]{0.4\linewidth}
\centering
\includegraphics[height=5cm,width=6.8cm]{picture/%o_1101_trace22_cpu.jpg}
multi_cpu1_7}
\end{minipage}
}

\centering
\caption{RMS error for numerical schemes of multiplicative noise case}
\label{RMSmul1}
\end{figure}

At the end of this section, we present the  experiment result of the multiplicative noise case in one dimensional case with the diffusion coefficient $G(X)(x)=X(x)+1,x\in \mathcal O$ and the trace-class operator $Q$ satisfying  $Qe_i=\frac{1}{i^2}e_i,i \ge 1.$ Here, we only take \textit{type 1} timestep functions as an example, which are taken as
$$
 \tau^{\delta}_{11}(X)=\delta\left(\frac{\|X\|+0.6}{\|F(X)\|+1.4}\right)^{\frac{4}{3}},\quad 
 \tau^{\delta}_{12}(X)=\delta \min\left\{\frac{2\|X\|^4_{L^4}+0.01}{\|X\|^6_{L^6}+0.01},\left(\frac{\|X\|+0.3}{\|F(X)\|+0.8}\right)^{\frac{4}{3}}\right\}.
 $$
 Figure \ref{RMSmul1} shows the comparison of RMS error against the number of timesteps and CPU time, which implies the convergence orders in the temporal direction of schemes are $\frac12.$ For a given RMS error, it can be observed that in Figure \ref{RMSmul1} $(a)$, the adaptive schemes cost slightly fewer number of timesteps than TE with the uniform timestep, and that in Figure \ref{RMSmul1} $(b),$ the adaptive schemes cost less CPU time than TE. The total CPU time is listed in Table \ref{table3}, from which we can see that the coupled schemes are in a lower computational cost in terms of the CPU time. These show the better performance of adaptive schemes for the multiplicative noise case. 

\begin{table}[tbhp]
\begin{tabular}{cccccccc}
\hline
 & CAU 1 & CAU 2& CAU 3 & CAU a & CAU b & CAU c &TE \\
\hline
Type 1 & 48.88 &49.68 & 55.21 & 54.61 & 55.55 & 54.75 &97.26\\
\hline
\end{tabular}
\vspace{0.1in}
\caption{Time in seconds ($\delta=2^{-6}$)}
\label{table3}
\end{table}

\appendix
\section{Appendix}\label{appendix1}

\begin{proof}[Proof of Lemma \ref{lemma1.4}] 
$(\romannumeral1)$ Let $u=\sum_{i=1}^{\infty}\langle u,e_i\rangle e_i.$ Then the H\"older inequality gives that 
\begin{align*}
&\|A^{\rho}P^NS(t)u\|_E=\sup_{x\in\mathcal{O}}\Big|\sum_{i=1}^Ne^{-\lambda_it}\lambda^{\rho}_i\langle u,e_i\rangle e_i(x)\Big|\leq C\sum_{i=1}^Ne^{-\lambda_it}\lambda^{\rho}_i\lvert\langle u,e_i\rangle\rvert\\
&\leq C\Big(\sum_{i=1}^N\lambda^{2\rho}_ie^{-2\lambda_it}\Big)^{\frac{1}{2}}\Big(\sum_{i=1}^N\lvert\langle u,e_i\rangle\rvert^2\Big)^{\frac{1}{2}}\leq C\Big(\sum_{i=1}^N(i^{2/d})^{2\rho}e^{-cti^{2/d}}\Big)^{\frac{1}{2}}\|u\|.
\end{align*}
By the monotonicity of the function $y^{4\rho/d}e^{-cty^{2/d}}$ w.r.t. $y\in(0,\infty),$ which takes the maximum value at the point $y=(\frac{2\rho}{ct})^{d/2}$, and letting $i_0:=\lfloor (\frac{2\rho}{ct})^{d/2}\rfloor$ with $\lfloor \cdot\rfloor$ being the floor function, we have
\begin{align*}
&\sum_{i=1}^{\infty}i^{4\rho/d}e^{-cti^{2/d}}= i_0^{4\rho/d}e^{-cti_0^{2/d}}+\sum_{i=1}^{i_0-1}i^{4\rho/d}e^{-cti^{2/d}}+\sum_{i=i_0+1}^{\infty}i^{4\rho/d}e^{-cti^{2/d}}\\
&\leq \big(\frac{2\rho}{ct}\big)^{2\rho}e^{-2\rho}+\int_0^{\infty}y^{4\rho/d}e^{-cty^{2/d}}\mathrm{d}y\\
&\leq C(\rho,d)\Big(\frac{1}{t^{2\rho}}+\frac{1}{t^{2\rho+\frac d2}}\int_0^{\infty}z^{2\rho+\frac d2-1}e^{-cz}\mathrm{d}z\Big)\leq C(\rho,d)\Big(\frac{1}{t^{2\rho}}+\frac{1}{t^{2\rho+\frac d2}}\Big),
\end{align*}
where in the last step we have used the fact that $\int_0^{\infty}z^{2\rho+\frac d2-1}e^{-cz}\mathrm{d}z=c^{-2\rho-\frac d2}\Gamma(2\rho+\frac d2)$ with $\Gamma$ being the Gamma function.
Applying the inequality $(a+b)^{\frac{1}{2}}\leq a^{\frac{1}{2}}+b^{\frac{1}{2}}$ leads to the desired result.

$(\romannumeral2)$ By the H\"older inequality, we obtain
\begin{align*}
&\|P^NS(t)u\|_E=
\sup_{x\in\mathcal{O}}\Big|\sum_{i=1}^Ne^{-\lambda_i t}\langle u,e_i\rangle e_i(x)\Big|\leq C\sum_{i=1}^Ne^{-\lambda_it}|\langle u,e_i\rangle|\\
&\leq C\Big(\sum_{i=1}^Ne^{-2\lambda_it}\lambda_i^{-\rho}\Big)^{\frac12}\Big(\sum_{i=1}^N\lambda_i^{\rho}|\langle u,e_i\rangle|^2\Big)^{\frac12}\leq C\Big(\int_0^{\infty}y^{-2\rho/d}e^{-cy^{2/d}t}\mathrm{d}y\Big)^{\frac12}\|u\|_{\rho}\\
&\leq C(\rho,d)\Big(t^{\rho-\frac d2}\int_0^{\infty}z^{-\rho+d/2-1}e^{-cz}\mathrm{d}z\Big)^{\frac12}\|u\|_{\rho}\leq C(\rho,d)t^{\frac{\rho}{2}-\frac d4}\|u\|_{\rho},
\end{align*}
where in the last step we have used the fact that $\int_0^{\infty}z^{-\rho+\frac d2-1}e^{-cz}\mathrm{d}z=c^{\rho-\frac d2}\Gamma(-\rho+\frac d2)$ for $\rho\in[0,\frac d2).$
The proof is finished. 
\end{proof}

\begin{proof}[Proof of Lemma \ref{purt_lemma}]
By Lemma \ref{lemma1.4} (\romannumeral2) with $\rho=0$ and the Sobolev embedding $\dot{H}^1\hookrightarrow L^6(\mathcal{O})$, we get
\begin{align}\label{v^N}
\|v^N_t\|_E\leq C\int_0^t(t-s)^{-\frac d4}(1+\|\nabla v^N_s\|^3+\|z^N_s\|^3_{L^6(\mathcal{O})})\mathrm{d}s.
\end{align}
The Taylor formula, the Young inequality, the Gagliardo--Nirenberg inequality $\|u\|_{L^{\eta_1}(\mathcal{O})}\leq C\big(\|\Delta u\|^{\alpha}\|u\|^{1-\alpha}+\|u\|_{L^{\eta_2}(\mathcal{O})}\big)$ for $u\in H^2$ with $\alpha=\frac{(\eta_1-2)d}{4\eta_1}\in(0,1],\,\eta_1>2,\,\eta_2>0$ imply that
\begin{align}
\|\nabla v^N_t\|^2&=\int_0^t \Big(-2\|Av^N_s\|^2+2\langle \nabla v^N_s,\nabla F^N(v^N_s+z^N_s)\rangle\Big)\mathrm{d}s\notag\\
&\leq -(2-\epsilon)\int_0^t\|Av^N_s\|^2\mathrm{d}s\notag\\
&\quad +C(\epsilon)\int_0^t\Big(1+\|v^N_s\|^4_{L^4(\mathcal{O})}(\|z^N_s\|^2_E+1)+\|v^N_s\|^2\|z^N_s\|^4_E+\|z^N_s\|^6_E+\|v^N_s\|^2\Big)\mathrm{d}s\notag\\
&\leq -(2-\epsilon)\int_0^t\|Av^N_s\|^2\mathrm{d}s+C(\epsilon)\int_0^t\Big(1+(\|Av^N_s\|^{\frac d2}\|v^N_s\|^{4-\frac d2}+\|v^N_s\|^4)(\|z^N_s\|^2_E+1)+\|v^N_s\|^2\|z^N_s\|^4_E\notag\\
&\quad +\|z^N_s\|^6_E+\|v^N_s\|^2\Big)\mathrm{d}s\notag\\
&\leq C(\epsilon)\int_0^t\Big(1+\|v^N_s\|^{\frac{2(8-d)}{4-d}}(\|z^N_s\|^{\frac{8}{4-d}}_E+1)+\|v^N_s\|^6+\|z^N_s\|^6_E\Big)\mathrm{d}s,\label{nablav^N}
\end{align}
where in the second step we use the fact that $\langle \nabla v^N_s,3a_3(v^N_s)^2\nabla v^N_s\rangle \leq 0,$ and in the third step we take $\eta_1=4, \eta_2=2.$
Using the Taylor formula again, and combining  \eqref{one-sided Lip} and the Young inequality give that
\begin{align*}
\|v^N_t\|^2&\leq -2\int_0^t\|\nabla v^N_s\|^2\mathrm{d}s+2L_0\int_0^t\|v^N_s\|^2\mathrm{d}s+2\int_0^t\langle v^N_s,F^N(z^N_s)\rangle \mathrm{d}s\\
&\leq 2L_0\int_0^t\|v^N_s\|^2\mathrm{d}s+C\int_0^t\big( 1+\|v^N_s\|^2+\|z^N_s\|^6_E\big)\mathrm{d}s,
\end{align*}
which yields that $\|v^N_t\|^2\leq Ce^{(2L_0+C)T}\int_0^t(1+\|z^N_s\|^6_E)\mathrm{d}s$ due to the Gr\"onwall inequality.
This, together with  \eqref{nablav^N} leads to
\begin{align*}
\sup_{0\leq s\leq t}\|\nabla v^N_s\|^2&\leq C\Big[\sup_{0\leq s\leq t}\|v^N_s\|^{\frac{2(8-d)}{4-d}}\int_0^t(\|z^N_s\|_E^{\frac{8}{4-d}}+1)\mathrm{d}s+\int_0^t(1+\|z^N_s\|^6_E)\mathrm{d}s+\sup_{0\leq s\leq t}\|v^N_s\|^6\Big]\\
&\leq C\Big[1+\Big(\int_0^t(1+\|z^N_s\|^8_E)\mathrm{d}s\Big)^{\frac{8-d}{4-d}+1}\Big].
\end{align*}
Plugging the above inequality into the right-hand side of \eqref{v^N}, one has
\begin{align*}
\|v^N_t\|_E&\leq C\int_0^t(t-s)^{-\frac d4}(1+\|z^N_s\|^3_E)\mathrm{d}s+C\sup_{0\leq s\leq t}\|\nabla v^N_s\|^3\\
&\leq C\int_0^t(t-s)^{-\frac d4}(1+\|z^N_s\|^3_E)\mathrm{d}s+C\Big[1+\Big(\int_0^t(1+\|z^N_s\|^8_E)\mathrm{d}s\Big)^{\frac{8-d}{4-d}+1}\Big]^{\frac32}\\
&\leq C\Big[1+\int_0^t\Big((t-s)^{-\frac d4}\|z^N_s\|^3_E+\|z^N_s\|^{72}_E\Big)\mathrm{d}s\Big],
\end{align*}
where we have used the fact that $\frac32(\frac{8-d}{4-d}+1)\leq 9$ for $d\leq 3.$ The proof is finished. 
\end{proof}

\begin{proof}[Proof of Proposition \ref{thm2}]
First, define a sequence of non-increasing events as follows:
\begin{align*}
\Omega_t:=\Big\{\omega\in\Omega:\sup_{j\in\{0,1,\cdots,m_t(\omega)\}}(\tau^{\delta}_j)^{\alpha}\|X^N_{t_j}(\omega)\|_E\leq 1\Big\},\quad 0\leq t\leq T
\end{align*}
with $\alpha>0$ being determined later.
It is clear that $\chi_{\Omega_t}\in\mathcal{F}_{t}.$ Note that one can always choose $\delta>0$ sufficiently small so that $\mathbb{P}\big(\big\{(\tau^{\delta}_{0})^{\alpha}\|X^N_0\|_E>1\big\}\big)\leq (T\delta)^{\alpha}\mathbb{E}[\|X^N_0\|_E]< 1,$ and hence $\mathbb{P}(\Omega_0)>0$, which implies $\Omega_0\neq \emptyset.$  Intuitively, $\Omega_t$ is the event that the $E$-norm of the numerical solution can be bounded by timestep sizes with a certain order before time $t$ (including $t$).

We claim that 
\begin{align}\label{claim1}
\mathbb{E}\left[\chi_{\Omega_t}\|X^N_t\|^p_E\right]
\leq C\Big(p,T,L_1,L_2,\|A^{\frac{\beta-1}{2}}Q^{\frac{1}{2}}\|_{\mathcal{L}_2(H)},\|X_0\|_{L^{q_1}(\Omega;\dot{H}^{\beta})},\|P^NX_0\|_{L^{q_1}(\Omega;E)}\Big)
\end{align}
and 
\begin{align}\label{claim2}
\mathbb{E}\left[\chi_{\Omega^c_t}\|X^N_t\|^p_E\right]
\leq C\Big(p,T,\tau_{min},L_1,L_2,\|A^{\frac{\beta-1}{2}}Q^{\frac{1}{2}}\|_{\mathcal{L}_2(H)},\|X_0\|_{L^{q_2}(\Omega;\dot{H}^{\beta})},\|P^NX_0\|_{L^{q_2}(\Omega;E)}\Big)
\end{align}
with $q_1,q_2>0$ sufficiently large, whose proofs are given in \textit{Step 1} and \textit{Step 2}, respectively. Once we prove \eqref{claim1} and \eqref{claim2}, the proof of \eqref{stability} is finished due to the Minkowski inequality.

\textit{Step 1: Proof of \eqref{claim1}.}

Let
\begin{align}\label{def_Z}
Z^N_t:=S^N(t)P^NX_0+\int_0^tS^N(t-t_{m_s})F^N(X^N_{t_{m_s}})\mathrm{d}s-\int_0^tS^N(t-s)F^N(X^N_s)\mathrm{d}s+W^N_A(t),
\end{align}
where $W^N_A(t):=P^NW_A(t)$ with $W_A(t)=\int_0^tS(t-t_{m_s})\mathrm{d}W(s).$
Then $X^N_t=\int_0^tS^N(t-s)F^N(X^N_s)\mathrm{d}s+Z^N_t.$
Let $\bar{X}^N_t:=X^N_t-Z^N_t.$
It satisfies
\begin{align*}
\begin{cases}
\mathrm{d}\bar{X}^N_t=(-A^N\bar{X}^N_t+F^N(\bar{X}^N_t+Z^N_t))\mathrm{d}t,\quad t\in(0,T],\\
\bar{X}^N_0=0.
\end{cases}
\end{align*}
Then applying  Lemma \ref{purt_lemma} yields that $\|\bar{X}^N_t\|_E\leq C\int_0^t\Big(1+(t-s)^{-\frac d4}\|Z^N_{s}\|^3_{E}+\|Z^N_s\|^{72}_E\Big)\mathrm{d}s$. This implies
\begin{align}\label{thm1eq2}
\left\|\chi_{\Omega_{t}}\|\bar{X}^N_{t}\|_{E}\right\|_{L^p(\Omega)}\leq C\int_0^t\Big(1+(t-s)^{-\frac d4}\Big\|\chi_{\Omega_{t}}\|Z^N_s\|_E\Big\|^3_{L^{3p}(\Omega)}+\Big\|\chi_{\Omega_t}\|Z^N_s\|_E\Big\|^{72}_{L^{72p}(\Omega)}\Big)\mathrm{d}s.
\end{align}
It follows from the definition of $Z^N_s$ that for $0\leq s\leq t,$
\begin{align}\label{thm1eq3}
\chi_{\Omega_{t}}\|Z^N_s\|_E&\leq \|S^N(s)P^NX_0\|_E+\chi_{\Omega_{t}}\left\|\int_0^sS^N(s-r)(F^N(X^N_{t_{m_r}})-F^N(X^N_r))\mathrm{d}r\right\|_E\nonumber\\
&\quad+\chi_{\Omega_{t}}\left\|\int_0^s(S^N(s-t_{m_r})-S^N(s-r))F^N(X^N_{t_{m_r}})\mathrm{d}r\right\|_E+\|W^N_A(s)\|_E\nonumber\\
&=:\|S^N(s)P^NX_0\|_E+I_0(s)+I_1(s)+\|W^N_A(s)\|_E.
\end{align}

\textit{Estimate of $I_0.$}
Since
\begin{align}\label{since}
X^N_r=S^N(r-t_{m_r})X^N_{t_{m_r}}
+\int_{t_{m_r}}^rS^N(r-t_{m_u})F^N(X^N_{t_{m_u}})\mathrm{d}u+W^N_A(r)-S^N(r-t_{m_r})W^N_A(t_{m_r}),
\end{align}
and by the contractivity of the semigroup in $E$, i.e., $\|S(t)u\|_{E}\leq \|u\|_{E}$,  \eqref{adaptive1}, Assumption \ref{assumption6}, and Lemma \ref{lemma1.4} $(\romannumeral2)$ with $\rho=0$, we obtain
\begin{align}\label{chiX^N_r}
\chi_{\Omega_{t}}\|X^N_r\|_E
&\leq \chi_{\Omega_{t}}\Big(\|X^N_{t_{m_r}}\|_E+C\int_{t_{m_r}}^r(r-t_{m_u})^{-\frac{d}{4}}\|F^N(X^N_{t_{m_u}})\|\mathrm{d}u+\|W^N_A(r)\|_E+\|W^N_A(t_{m_r})\|_E\Big)\notag\\
&\leq \chi_{\Omega_{t}}\left(\|X^N_{t_{m_r}}\|_E+C(\tau^{\delta}_{m_r})^{1-\frac{d}{4}}\|F(X^N_{t_{m_r}})\|+\|W^N_A(r)\|_E+\|W^N_A(t_{m_r})\|_E\right)\notag\\
&\leq \chi_{\Omega_{t}}\left(\|X^N_{t_{m_r}}\|_E+CL_2+\|W^N_A(r)\|_E+\|W^N_A(t_{m_r})\|_E\right).
\end{align}
We derive from \eqref{S(t)1}-\eqref{S(t)2} and  \eqref{continuous version}-\eqref{adaptive1} that for $\beta\in(0,2],\,0\leq r\leq t,$
\begin{align}\label{chiX^N_r-X^N_s}
\chi_{\Omega_{t}}\|X^N_r-X^N_{t_{m_r}}\|
&\leq\|(S(r)-S(t_{m_r}))P^NX_0\|+\chi_{\Omega_{t}}\Big\|\int_0^rS(r-t_{m_u})F^N(X^N_{t_{m_u}})\mathrm{d}u\notag\\
&\quad -\int_0^{t_{m_r}}S(t_{m_r}-t_{m_u})F^N(X^N_{t_{m_u}})\mathrm{d}u\Big\|
+\|P^N(W_A(r)-W_A(t_{m_r}))\|\notag\\
&\leq (\tau^{\delta}_{m_r})^{\frac{\beta}{2}}\|P^NX_0\|_{\beta}+\chi_{\Omega_{t}}\Big\|\int_0^{t_{m_r}}S(t_{m_r}-t_{m_u})(S(r-t_{m_r})-\mathrm{Id})F^N(X^N_{t_{m_u}})\mathrm{d}u\Big\|\notag\\
&\quad +\chi_{\Omega_t}\Big\|\int_{t_{m_r}}^rS(r-t_{m_r})F^N(X^N_{t_{m_r}})\mathrm{d}u\Big\|+\|P^N(W_A(r)-W_A(t_{m_r}))\|\notag\\
&\leq (\tau^{\delta}_{m_r})^{\frac{\beta}{2}}\|P^NX_0\|_{\beta}+C\chi_{\Omega_t}\int_0^{t_{m_r}}(t_{m_r}-t_{m_u})^{-1+\epsilon}(\tau^{\delta}_{m_r})^{1-\epsilon}\|F(X^N_{t_{m_u}})\|\mathrm{d}u\notag\\
&\quad +C\chi_{\Omega_t}\tau_{m_r}^{\delta}\|F(X^N_{t_{m_r}})\|+\|W_A(r)-W_A(t_{m_r})\|\notag\\
&\leq C(L_2)(\tau^{\delta}_{m_r})^{\frac{\beta\wedge1\wedge \frac d2}{2}-\epsilon}\Big(1+\|X_0\|_{\beta}+\frac{\|W_A(r)-W_A(t_{m_r})\|}{(\tau^{\delta}_{m_r})^{\frac{\beta\wedge1}{2}}}\Big),
\end{align}
where $0<\epsilon\ll1$, and in the last step we have used $\frac{\tau^{\delta}_{m_r}}{\tau^{\delta}_{m_u}}\leq \frac{T\delta}{\delta\min\{T,\tau_{m_u}\}}\leq T(T^{-1}+\tau^{-1}_{min})$ due to \eqref{adaptive2} and Assumption \ref{assumption6}.
We deduce from Lemma \ref{lemma1.4} $(\romannumeral2)$ with $\rho=0$, \eqref{polynomial Lip}, and \eqref{chiX^N_r} that
\begin{align*}
I_0(s)&\leq \chi_{\Omega_{t}}\int_0^s(s-r)^{-\frac{d}{4}}\|F(X^N_r)-F(X^N_{t_{m_r}})\|\mathrm{d}r\\
&\leq C(L_1)\chi_{\Omega_{t}}\int_0^s(s-r)^{-\frac{d}{4}}(1+\|X^N_r\|^2_E+\|X^N_{t_{m_r}}\|^2_E)\|X^N_r-X^N_{t_{m_r}}\|\mathrm{d}r\notag\\
&\leq C(L_1,L_2)\chi_{\Omega_{t}}\int_0^s(s-r)^{-\frac{d}{4}}\Big(1+\|X^N_{t_{m_r}}\|^2_E+\|W^N_A(r)\|^2_E+\|W^N_A(t_{m_r})\|^2_E\Big)\|X^N_r-X^N_{t_{m_r}}\|\mathrm{d}r\notag,
\end{align*}
which together with   \eqref{chiX^N_r-X^N_s} leads to
\begin{align*}
I_0(s)&\leq C(L_1,L_2)\chi_{\Omega_{t}}\int_0^s(s-r)^{-\frac{d}{4}}\Big(1+\|X^N_{t_{m_r}}\|^2_E+\|W^N_A(r)\|^2_E+\|W^N_A(t_{m_r})\|^2_E\Big)\times \notag\\
&\quad(\tau^{\delta}_{m_r})^{\frac{\beta\wedge 1\wedge \frac d2}{2}}\Big(1+\|X_0\|_{\beta}+\frac{\|W_A(r)-W_A(t_{m_r})\|}{(\tau^{\delta}_{m_r})^{\frac{\beta\wedge 1}{2}}}\Big)\mathrm{d}r\\
&\leq C(L_1,L_2)\int_0^s(s-r)^{-\frac{d}{4}}(1+\|W^N_A(r)\|^2_E+\|W^N_A(t_{m_r})\|^2_E)\Big(1+\|X_0\|_{\beta}+\frac{\|W_A(r)-W_A(t_{m_r})\|}{(\tau^{\delta}_{m_r})^{\frac{\beta\wedge 1}{2}}}\Big)\mathrm{d}r,
\end{align*}
where we have used the definition of $\Omega_t$ and taken $\alpha\in\big(0, \frac{\beta\wedge 1\wedge \frac d2}{4}\big).$
Taking $L^{\tilde{p}}(\Omega)$-norm ($\tilde{p}\ge 2$) on both sides of the above equation, and applying the Minkowski inequality and the H\"older  inequality, we derive
\begin{align*}
\|I_0(s)\|_{L^{\tilde{p}}(\Omega)}
&\leq C(L_1,L_2)\int_0^s(s-r)^{-\frac{d}{4}}\Big\|1+\|W^N_A(r)\|^2_E+\|W^N_A(t_{m_r})\|^2_E\Big\|_{L^{2\tilde{p}}(\Omega)}\times \\
&\quad\Big\|1+\|X_0\|_{\beta}+\frac{\|W_A(r)-W_A(t_{m_r})\|}{(\tau^{\delta}_{m_r})^{\frac{\beta\wedge1}{2}}}\Big\|_{L^{2\tilde{p}}(\Omega)}\mathrm{d}r.
\end{align*}
Note that
$\sup_{t\in[0,T]}\|A^{\frac{\beta\wedge1}{2}}W_A( t)\|_{L^{2\tilde p}(\Omega;H)}\leq C$ and   \eqref{S(t)2} imply the H\"older continuity of $W_A$ 
\begin{align*}
\|W_A(r)-W_A(t_{m_r})\|_{L^{2\tilde{p}}(\Omega;H)}&\leq \|A^{\frac{\beta\wedge 1}{2}}(S(r-t_{m_r})-\mathrm{Id})\|_{\mathcal{L}(H)}\|W_A(t_{m_r})\|_{L^{2\tilde{p}}(\Omega;\dot{H}^{\beta\wedge1})}\\
&\quad +\Big\|\int_{t_{m_r}}^rS(r-t_{m_r})\mathrm{d}W(u)\Big\|_{L^{2\tilde{p}}(\Omega;H)}\leq C(\tilde{p},\|A^{\frac{\beta-1}{2}}Q^{\frac12}\|_{\mathcal{L}_2(H)})(\tau^{\delta}_{m_r})^{\frac{\beta\wedge1}{2}},
\end{align*}
which combining  properties of the conditional expectation 
gives
\begin{align*}
&\quad\mathbb{E}\Big[\Big(\frac{\|W_A(r)-W_A(t_{m_r})\|}{(\tau^{\delta}_{m_r})^{\frac{\beta\wedge1}{2}}}\Big)^{2\tilde{p}}\Big]=\mathbb{E}\Big[\mathbb{E}\Big(\frac{\|W_A(r)-W_A(t_{m_r})\|^{2\tilde{p}}}{(\tau^{\delta}_{m_r})^{\frac{\beta\wedge1}{2}\times 2\tilde{p}}}\Big|\mathcal{F}_{t_{m_r}}\Big)\Big]\\
&=\mathbb{E}\Big[(\tau^{\delta}_{m_r})^{-\frac{\beta\wedge1}{2}\times 2\tilde{p}}\mathbb{E}\Big(\|W_A(r)-W_A(t_{m_r})\|^{2\tilde{p}}\Big|\mathcal{F}_{t_{m_r}}\Big)\Big]\leq C\big(\tilde{p},\|A^{\frac{\beta-1}{2}}Q^{\frac{1}{2}}\|_{\mathcal{L}_2(H)}\big).
\end{align*}
Similar as the proof of \eqref{sup_noise}, we have  \begin{align}\label{supW_A(t)}
\sup_{0\leq t\leq T}\mathbb{E}\Big[\|W^N_A(t)\|^{4\tilde{p}}_E\Big]<\infty.
\end{align}
Therefore, 
we arrive at
\begin{align*}
\sup_{0\leq s\leq t}\|I_0(s)\|_{L^{\tilde{p}}(\Omega)}\leq C\big(\tilde{p},T,L_1,L_2,\|A^{\frac{\beta-1}{2}}Q^{\frac{1}{2}}\|_{\mathcal{L}_2(H)}\big)\big(1+\|X_0\|_{L^{2\tilde{p}}(\Omega;\dot{H}^{\beta})}\big).
\end{align*}

\textit{Estimate of $I_1.$}
Applying Lemma \ref{lemma1.4} $(\romannumeral1)$ with $\rho=1-\frac{d}{4}-\epsilon=:\rho_0$ and \eqref{S(t)2} yields that for $0\leq s\leq t,$
\begin{align*}
I_1(s)&=\chi_{\Omega_{t}}\left\|\int_0^sA^{\rho_0}S^N(s-r)A^{-\rho_0}(S^N(r-t_{m_r})-\mathrm{Id})F^N(X^N_{t_{m_r}})\mathrm{d}r\right\|_E\\
&\leq C(L_1)\chi_{\Omega_{t}}\int_0^s\Big((s-r)^{-\rho_0}+(s-r)^{-\rho_0-\frac{d}{4}}\Big)(\tau^{\delta}_{m_r})^{\rho_0}(1+\|X^N_{t_{m_r}}\|^3_E)\mathrm{d}r\\
&\leq C(L_1)\int_0^s\Big((s-r)^{-\rho_0}+(s-r)^{-\rho_0-\frac{1}{4}}\Big)\mathrm{d}r,
\end{align*}
where we have used the definition of $\Omega_t$ and taken $\alpha\in(0,\frac{\rho_0}{3}].$
Therefore, we get
$
\sup_{0\leq s\leq t}\|I_1(s)\|_{L^{\tilde{p}}(\Omega)}\leq C(\tilde{p},T,L_1).
$ 

Combining estimates of $I_0$ and $I_1$, and taking parameter $\alpha$ in $\Omega_t$ as $\alpha=\alpha_0:=(\frac{\beta\wedge1\wedge \frac d2}{4}\wedge \frac{4-d}{12})-\epsilon,$
we derive that
\begin{align*}
&\quad\sup_{s\in[0,T]}\Big[\big\|\chi_{\Omega_t}\|Z^N_s\|_E\big\|^3_{L^{3p}(\Omega)}+\big\|\chi_{\Omega_t}\|Z^N_s\|_E\big\|^{72}_{L^{72p}(\Omega)}\Big]\\
&\leq C\big(p,T,L_1,L_2,\|A^{\frac{\beta-1}{2}}Q^{\frac{1}{2}}\|_{\mathcal{L}_2(H)},\|P^NX_0\|_{L^{72p}(\Omega;E)},\|X_0\|_{L^{144p}(\Omega;\dot{H}^{\beta})}\big).
\end{align*}
Then \eqref{thm1eq2} can be estimated as
\begin{align}\label{eqto}
\Big\|\chi_{\Omega_{t}}\|\bar{X}^N_{t}\|_E\Big\|_{L^p(\Omega)}\leq C\big(p,T,L_1,L_2,\|A^{\frac{\beta-1}{2}}Q^{\frac{1}{2}}\|_{\mathcal{L}_2(H)},\|P^NX_0\|_{L^{72p}(\Omega;E)},\|X_0\|_{L^{144p}(\Omega;\dot{H}^{\beta})}\big),
\end{align}
which leads to
\begin{align*}
\Big\|\chi_{\Omega_{t}}\|{X}^N_{t}\|_E\Big\|_{L^p(\Omega)}&\leq \Big\|\chi_{\Omega_{t}}\|\bar{X}^N_{t}\|_E\Big\|_{L^p(\Omega)}+\Big\|\chi_{\Omega_{t}}\|Z^N_{t}\|_E\Big\|_{L^p(\Omega)}\\
&\leq C\big(p,T,L_1,L_2,\|A^{\frac{\beta-1}{2}}Q^{\frac{1}{2}}\|_{\mathcal{L}_2(H)},\|P^NX_0\|_{L^{72p}(\Omega;E)},\|X_0\|_{L^{144p}(\Omega;\dot{H}^{\beta})}\big).
\end{align*} 

\textit{Step 2: Proof of \eqref{claim2}.} 
Recall that in \textit{Step 1} we take $\alpha=\alpha_0.$ Note that for $0\leq t\leq T,$
\begin{align*}
\Omega^c_{t}=\Omega^c_{t_{m_t}}\cup\Big(\Omega_{t_{m_t}}\cap \left\{\omega\in\Omega:(\tau^{\delta}_{m_t})^{\alpha_0}\|X^N_{t_{m_t}}(\omega)\|_E>1\right\}\Big).
\end{align*}
Hence, by iteration,
\begin{align*}
\chi_{\Omega^c_{t}}(\omega)&=\chi_{\Omega^c_{t_{m_t}}}(\omega)+\chi_{\Omega_{t_{m_t}}}(\omega)\cdot \chi_{\{(\tau^{\delta}_{m_t})^{\alpha_0}\|X^N_{t_{m_t}}\|_E>1\}}(\omega)\\
&=\sum_{j=0}^{m_t(\omega)}\chi_{\Omega_{t_{j-1}}}(\omega)\cdot \chi_{\{(\tau^{\delta}_j)^{\alpha_0}\|X^N_{t_j}\|_E>1\}}(\omega),
\end{align*}
where $\chi_{\Omega_{t_{-1}}}:=1.$
Then for $0\leq t\leq T,$ combining the H\"older inequality yields
\begin{align*}
\mathbb{E}\Big[\chi_{\Omega^c_{t}}\|X^N_{t}\|^p_E\Big]&=\mathbb{E}\Big[\sum_{j=0}^{m_t}\|X^N_{t}\|^p_E\cdot\chi_{\Omega_{t_{j-1}}}\cdot\chi_{\{(\tau^{\delta}_j)^{\alpha_0}\|X^N_{t_j}\|_E>1\}}\Big]\\
&=\mathbb{E}\Big[\sum_{j=0}^{m_t}\Big(\chi_{\Omega_{t_{j-1}}}\cdot\chi_{\{(\tau^{\delta}_j)^{\alpha_0}\|X^N_{t_j}\|_E>1\}}\tau^{\delta}_j\Big)\frac{1}{\tau^{\delta}_j}\|X^N_{t}\|^p_E\Big]\\
&\leq \mathbb{E}\Big[\sup_{0\leq t_j\leq T}\big(\frac{1}{\tau^{\delta}_{j}}\big)\int_0^{t} \chi_{\Omega_{t_{m_s-1}}}\cdot\chi_{\{(\tau^{\delta}_{m_s})^{\alpha_0}\|X^N_{t_{m_s}}\|_E>1\}}\mathrm{d}s\|X^N_{t}\|^p_E\Big]\\
&\leq \Big(\mathbb{E}\Big[\sup_{0\leq t_j\leq T}\big(\frac{1}{\tau^{\delta}_j}\big)^3\Big]\Big)^{\frac{1}{3}}\Big(\mathbb{E}\Big[\int_0^{t}\chi_{\Omega_{t_{m_s-1}}}\cdot\chi_{\{(\tau^{\delta}_{m_s})^{\alpha_0}\|X^N_{t_{m_s}}\|_E>1\}}\mathrm{d}s\Big]\Big)^{\frac{1}{3}}\left(\mathbb{E}\big[\|X^N_{t}\|^{3p}_E\big]\right)^{\frac{1}{3}}.
\end{align*}
It can be shown that
\begin{align}\label{eq1}
\Big(\mathbb{E}\Big[\sup_{0\leq t_j\leq T}\big(\frac{1}{\tau^{\delta}_j}\big)^3\Big]\Big)^{\frac{1}{3}}\leq C(T,\tau_{min})\delta^{-1}.
\end{align}
Notice that 
\begin{align*}
\mathbb{E}\Big[\chi_{\Omega_{t_{m_s-1}}}\cdot\chi_{\{(\tau^{\delta}_{m_s})^{\alpha_0}\|X^N_{t_{m_s}}\|_E>1\}}\Big]&=\mathbb{P}\left(\Big\{\omega\in\Omega_{t_{m_s-1}}:(\tau^{\delta}_{m_s})^{\alpha_0}\|X^N_{t_{m_s}}\|_E>1\Big\}\right)\\
&=\mathbb{P}\left(\Big\{\omega\in\Omega:\chi_{\Omega_{t_{m_s-1}}}\cdot(\tau^{\delta}_{m_s})^{\alpha_0}\|X^N_{t_{m_s}}\|_E>1\Big\}\right).
\end{align*}
Combining the Chebyshev inequality, the H\"{o}lder  inequality and Assumption \ref{assumption6} gives 
\begin{align}\label{eq2}
&\quad \mathbb{P}\left(\Big\{\omega\in\Omega:\chi_{\Omega_{t_{m_s-1}}}\cdot(\tau^{\delta}_{m_s})^{\alpha_0}\|X^N_{t_{m_s}}\|_E>1\Big\}\right)\\
&\leq \mathbb{E}\Big[\Big(\chi_{\Omega_{t_{m_s-1}}}\cdot(\tau^{\delta}_{m_s})^{\alpha_0}\|X^N_{t_{m_s}}\|_E\Big)^{\frac{3(p+1)}{\alpha_0}}\Big]\nonumber\\
&\leq \Big(\mathbb{E}\big[(\tau^{\delta}_{m_s})^{6(p+1)}\big]\Big)^{\frac{1}{2}}\Big(\mathbb{E}\Big[\chi_{\Omega_{t_{m_s-1}}}\|X^N_{t_{m_s}}\|^{\frac{6(p+1)}{\alpha_0}}_E\Big]\Big)^{\frac{1}{2}}\nonumber\\
&\leq C\Big(p,T,\tau_{min},L_1,\|A^{\frac{\beta-1}{2}}Q^{\frac{1}{2}}\|_{\mathcal{L}_2(H)},\|P^NX_0\|_{L^q(\Omega;E)},\|X_0\|_{L^{2q}(\Omega;\dot{H}^{\beta})}\Big)\delta^{3(p+1)},\nonumber
\end{align}
where we have utilized the fact that
\begin{align}\label{eqstep1}
\mathbb{E}\Big[\chi_{\Omega_{t_{m_s-1}}}\|X^N_{t_{m_s}}\|^{\frac{6(p+1)}{\alpha_0}}_E\Big]
\leq C\Big(p,T,\tau_{min},L_1,\|A^{\frac{\beta-1}{2}}Q^{\frac{1}{2}}\|_{\mathcal{L}_2(H)},\|P^NX_0\|_{L^q(\Omega;E)},\|X_0\|_{L^{2q}(\Omega;\dot{H}^{\beta})}\Big)
\end{align}
with $q=72\times \frac{6(p+1)}{\alpha_0}$. The proof of \eqref{eqstep1} is a combination of  \eqref{claim1} and
\begin{align*}
\chi_{\Omega_{t_{m_s-1}}}\|X^N_{t_{m_s}}\|_E
&=\chi_{\Omega_{t_{m_s-1}}}\|S^N(\tau_{m_s-1})\big(X^N_{t_{m_s-1}}+F^N(X^N_{t_{m_s-1}})\tau_{m_s-1}+P^N\Delta W_{m_s-1}\big)\|_E\\
&\leq C\chi_{\Omega_{t_{m_s-1}}}\big(\|X^N_{t_{m_s-1}}\|_E+\tau_{m_s-1}^{1-\frac{d}{4}}\|F(X^N_{t_{m_s-1}})\|+\|P^N\Delta W_{m_s-1}\|_E\big)\\
&\leq C\chi_{\Omega_{t_{m_s-1}}}\big(\|X^N_{t_{m_s-1}}\|_E+L_2+\|P^N\Delta W_{m_s-1}\|_E\big).
\end{align*} 
It remains to estimate $\big(\mathbb{E}\big[\|X^N_{t}\|^{3p}_E\big]\big)^{\frac{1}{3}}$.
For $0\leq t_{m+1}\leq T$ and the fixed $\omega\in\Omega,$ it follows from Lemma \ref{lemma1.4} $(\romannumeral2)$ with $\rho=0$ and the contractivity of the semigroup $\{S(t),t\ge0\}$ in $E$ that
\begin{align*}
\|X^N_{t_{m+1}}\|_E&\leq \|S(\tau^{\delta}_m)X^N_{t_m}\|_E+\|S(\tau^{\delta}_m)F^N(X^N_{t_m})\tau^{\delta}_m\|_E+\|S(\tau^{\delta}_m)P^N\Delta W_m\|_E\\
&\leq \|X^N_{t_m}\|_E+C(\tau^{\delta}_m)^{1-\frac{d}{4}}\|F(X^N_{t_m})\|+\|W^N_A(t_{m+1})\|_E+\|W^N_A(t_m)\|_E\\
&\leq \|P^NX_{0}\|_E+C\sum_{i=0}^m\Big((\tau^{\delta}_i)^{1-\frac{d}{4}}\|F(X^N_{t_i})\|+\|W^N_A(t_{i+1})\|_E+\|W^N_A(t_i)\|_E\Big).
\end{align*}
Similarly,
\begin{align*}
\|X^N_t\|_E\leq \|X^N_{t_{m_t}}\|_E+C(t-t_{m_t})^{1-\frac{d}{4}}\|F(X^N_{t_{m_t}})\|+\|W^N_A(t)\|_E+\|W^N_A(t_{m_t})\|_E.
\end{align*}
Hence, by Assumption \ref{assumption 4}, we obtain
\begin{align*}
\|X^N_t\|_E&\leq \|P^NX_0\|_E+C\sum_{i=0}^{m_t}\Big((\tau^{\delta}_i)^{1-\frac{d}{4}}\|F(X^N_{t_i})\|+\|W^N_A(t)\|_E+\|W^N_A(t_i)\|_E\Big).
\end{align*}
This, together with \eqref{supW_A(t)}  implies that for $0\leq t\leq T,$
\begin{align*}%\label{L^2omega}
\|X^N_t\|_{L^p(\Omega;E)}\leq \|P^NX_{0}\|_{L^p(\Omega;E)}+C\int_0^T\|(\tau^{\delta}_{m_s})^{-1}\|_{L^{2p}(\Omega)}(L_2+\|W^N_A(s)\|_{L^{2p}(\Omega;E)}+\|W^N_A(t_{m_s})\|_{L^{2p}(\Omega;E)})\mathrm{d}s.
\end{align*}
It can be deduced from
% \begin{align}\label{M_delta^p}
% M_{T,\delta}^p\leq \Big(\frac{T}{\tau^{\delta}}\Big)^p\leq \delta^{-p}\Big(\frac{T}{\tau_{min}+T}\Big)^p
% \end{align}
$(\tau^{\delta})^{-1}\leq \delta^{-1}(\tau_{min}+T)^{-1}$
and \eqref{supW_A(t)} that
\begin{align}
\label{eq3}
\mathbb{E}\big[\|X^N_t\|^{3p}_E\big]\leq C\big(p,T,\tau_{min},L_2,\|A^{\frac{\beta-1}{2}}Q^{\frac{1}{2}}\|_{\mathcal{L}_2(H)}\big)\big(\|P^NX_0\|^{3p}_{L^{3p}(\Omega;E)}+\delta^{-3p}\big).
\end{align}
It follows from \eqref{eq1}, \eqref{eq2} and \eqref{eq3} that
\begin{align*}
\mathbb{E}\Big[\chi_{\Omega^c_{t}}\|X^N_{t}\|^p_E\Big]&\leq C\delta^{-1}\delta^{p+1}
(1+\|P^NX_0\|^p_{L^{3p}(\Omega;E)}+\delta^{-p})
\leq C\end{align*}
with $C>0$ depending on $p,T, L_1,L_2,\tau_{min},\|A^{\frac{\beta-1}{2}}Q^{\frac{1}{2}}\|_{\mathcal{L}_2(H)},\|P^NX_0\|_{L^q(\Omega;E)},\|X_0\|_{L^{2q}(\Omega;\dot{H}^{\beta})}$.

Therefore, combining \textit{Step 1} and \textit{Step 2},  we get 
\begin{align*}
\sup_{0\leq t\leq T}\mathbb{E}[\|X^N_{t}\|^p_E]
\leq C\Big(p,T,L_1,L_2,\tau_{min},\|A^{\frac{\beta-1}{2}}Q^{\frac{1}{2}}\|_{\mathcal{L}_2(H)},\|X_0\|^q_{L^{2q}(\Omega;\dot{H}^{\beta}),\|P^NX_0\|^q_{L^{q}(\Omega;E)}}\Big)
\end{align*}
for some large $q:=q(p,\beta,d)>0,$
which finishes the proof.
\end{proof}

\begin{proof}[Proof of Proposition \ref{coro}]
 For the proof of the regularity \eqref{coroeq1}, 
when $\beta\in(0,2),$ applying \eqref{S(t)1} and  \eqref{S(t)3} yields that
$$
\sup_{0\leq t\leq T}\|X^N_t\|_{L^p(\Omega;\dot{H}^{\beta})}\leq C(1+\|X_0\|_{L^p(\Omega;\dot{H}^{\beta})}+\sup_{0\leq s\leq T}\|X^N_s\|^3_{L^{3p}(\Omega; E)}+\|A^{\frac{\beta-1}{2}}Q^{\frac{1}{2}}\|_{\mathcal{L}_2(H)}).
$$
And when $\beta=2,$ it follows from \eqref{S(t)1}, \eqref{F_normH^1}, and the H\"older inequality that 
\begin{align*}
\|X^N_t\|_{L^p(\Omega;\dot{H}^2)}&\leq \|X_0\|_{L^p(\Omega;\dot{H}^2)}+C\int_0^t(t-s)^{-\frac{1}{2}}\big\|\|X^N_{t_{m_s}}\|_1(1+\|X^{N}_{t_{m_s}}\|^2_E)\big\|_{L^p(\Omega)}\mathrm{d}s+C\|A^{\frac12}Q^{\frac12}\|_{\mathcal{L}_2(H)}\\
&\leq C(1+\|X_0\|_{L^p(\Omega;\dot{H}^2)}+\sup_{0\leq s\leq T}\|X^N_s\|_{L^{2p}(\Omega;\dot{H}^1)}+\sup_{0\leq s\leq T}\|X^N_s\|^2_{L^{4p}(\Omega;E)}).
\end{align*}
These, combining  Proposition \ref{thm2} finish the proof of \eqref{coroeq1}.

For the proof of  the H\"older continuity \eqref{coroeq11},
when $\gamma\in(0,\beta]$ with $\beta\in(0,2)$, or $\gamma\in(0,\beta)$ with $\beta=2$, we can deduce from \eqref{S(t)1}-\eqref{S(t)3}, Proposition \ref{thm2},  and \eqref{coroeq1} that
\begin{align*}
 \|X^N_t-X^N_s\|_{L^p(\Omega;\dot{H}^{\gamma})}
&\leq \|A^{\frac{\gamma-\beta}{2}}(S^N(t-s)-\mathrm{Id})A^{\frac{\beta}{2}}X^N_s\|_{L^p(\Omega;H)}+\Big\|\int_s^tA^{\frac{\gamma}{2}}S^N(t-t_{m_s})F^N(X^N_{t_{m_s}})\mathrm{d}s\Big\|_{L^p(\Omega;H)}\\
&\quad +\Big\|\int_s^tA^{\frac{1-\beta+\gamma}{2}}S^N(t-t_{m_s})A^{\frac{\beta-1}{2}}\mathrm{d}W(s)\Big\|_{L^p(\Omega;H)}\\
&\leq C(t-s)^{\frac{\beta-\gamma}{2}}\|X^N_s\|_{L^p(\Omega;\dot{H}^{\beta})}+C\int_s^t(t-t_{m_s})^{-\frac{\gamma}{2}}(1+\|X^N_{t_{m_s}}\|^3_{L^{3p}(\Omega;E)})\mathrm{d}s\\
&\quad+C\|A^{\frac{\beta-1}{2}}Q^{\frac12}\|_{\mathcal{L}_2(H)}(t-s)^{\frac{(\beta-\gamma)\wedge 1}{2}}\leq C(t-s)^{\frac{(\beta-\gamma)\wedge 1}{2}}.
\end{align*}
And when $\gamma=\beta=2,$ it follows from \eqref{coroeq1}
that $\|X^N_t-X^N_s\|_{L^p(\Omega;\dot{H}^{2})}\leq C$. The proof is finished.
\end{proof}

\section{Appendix}\label{appendix2}
\begin{proof}
[Proof of the convergence of  \eqref{ATEU}.] 
Based on the auxiliary process $Y^N_t,$
the error can be divided into the following terms,
\begin{align*}
\|X(t)-X^{N,(1)}_t\|_{L^p(\Omega;H)}
\leq \|X(t)-P^NX(t)\|_{L^p(\Omega;H)}+\|P^NX(t)-Y^N_t\|_{L^p(\Omega;H)}+\|Y^N_t-X^{N,(1)}_t\|_{L^p(\Omega;H)},
\end{align*}
where terms $\|X(t)-P^NX(t)\|_{L^p(\Omega;H)}$ and $\|P^NX(t)-Y^N_t\|_{L^p(\Omega;H)}$ have been estimated; see  \eqref{X-P^NX}-\eqref{eq4}.

For the estimate of the term $\|Y^N_t-X^{N,(1)}_t\|_{L^p(\Omega;H)},$ combining the differential form
\begin{align*}
\mathrm{d}(X^{N,(1)}-Y^N_t)&=-A^N(X^{N,(1)}-Y^N_t)\mathrm{d}t+S^N(t-t_{m_t})\Big(F^N(X^{N,(1)}_{t_{m_t}})\chi_{\{\tau^{\delta}_{m_t}\ge \tau^{\delta}_{min}\}}\\
&\quad+\frac{F^N(X^{N,(1)}_{t_{m_t}})}{1+\|F^N(X^{N,(1)}_{t_{m_t}})\|\tau^{\delta}_{min}}\chi_{\{\tau^{\delta}_{m_t}<\tau^{\delta}_{min}\}}\Big)\mathrm{d}t-F^N(X(t))\mathrm{d}t
\end{align*}
with $\tau^{\delta}_{min}=\delta\tau_{min},$ and applying the Taylor formula give 
\begin{align*}
\|X^{N,(1)}_t-Y^N_t\|^2=&\;2\int_0^t\langle X^{N,(1)}_s-Y^N_s,\;-A^N(X^{N,(1)}_s-Y^N_s)\rangle\mathrm{d}s\\
&+2\int_0^t\Big\langle X^{N,(1)}_s-Y^N_s,\;S^N(s-t_{m_s})F^N(X^{N,(1)}_{t_{m_s}})-F^N(X(s))\Big\rangle\chi_{\{\tau^{\delta}_{m_s}\ge \tau^{\delta}_{min}\}}\mathrm{d}s\\
& +2\int_0^t\Big\langle X^{N,(1)}_s-Y^N_s,\;S^N(s-t_{m_s})\frac{F^N(X^{N,(1)}_{t_{m_s}})}{1+\|F^N(X^{N,(1)}_{t_{m_s}})\|\tau^{\delta}_{min}}-F^N(X(s))\Big\rangle\chi_{\{\tau^{\delta}_{m_s}< \tau^{\delta}_{min}\}}\mathrm{d}s\\
=&:\;2\int_0^t\langle X^{N,(1)}_s-Y^N_s,\;-A^N(X^{N,(1)}_s-Y^N_s)\rangle\mathrm{d}s+2K_1(t)+2K_2(t).
\end{align*}
For the estimate of the term $K_1$,
 by \eqref{one-sided Lip} and the Young inequality, we have
\begin{align*}
K_1(t)&=\int_0^t\langle X^{N,(1)}_s-Y^N_s,\;(S^N(s-t_{m_s})-\mathrm{Id})F^N(X^{N,(1)}_{t_{m_s}}) \rangle\chi_{\{\tau^{\delta}_{m_s}\ge \tau^{\delta}_{min}\}}\mathrm{d}s\\
&\quad +\int_0^t\langle X^{N,(1)}_s-Y^N_s,\;F^N(X^{N,(1)}_{t_{m_s}})-F^N(Y^N_s)\rangle \chi_{\{\tau^{\delta}_{m_s}\ge \tau^{\delta}_{min}\}}\mathrm{d}s\\
&\quad+\int_0^t\langle X^{N,(1)}_s-Y^N_s,\;F^N(Y^N_s)-F^N(X(s))\rangle \chi_{\{\tau^{\delta}_{m_s}\ge \tau^{\delta}_{min}\}}\mathrm{d}s\\
&\leq (L_0+\frac{1}{2})\int_0^t\|X^{N,(1)}_s-Y^N_s\|^2\mathrm{d}s+\frac{1}{2}\int_0^t\|F^N(Y^N_s)-F^N(X(s))\|^2\mathrm{d}s+K_{1,1}(t)+K_{1,2}(t),
\end{align*}
where \begin{align*}&K_{1,1}(t):=\int_0^t\langle X^{N,(1)}_s-Y^N_s,\;(S^N(s-t_{m_s})-\mathrm{Id})F^N(X^{N,(1)}_{t_{m_s}}) \rangle\chi_{\{\tau^{\delta}_{m_s}\ge \tau^{\delta}_{min}\}}\mathrm{d}s,\\
&K_{1,2}(t):=\int_0^t\langle X^{N,(1)}_s-Y^N_s,\;F^N(X^{N,(1)}_{t_{m_s}})-F^N(X^{N,(1)}_s) \rangle\chi_{\{\tau^{\delta}_{m_s}\ge \tau^{\delta}_{min}\}}\mathrm{d}s.
\end{align*}

We first show the estimate of the term $K_{1,1}$.
For the case $\beta=1$ or $\beta=2$, similar to the estimate of $J_{1,1}$, the term  $K_{1,1}$ can be estimated as  $\|K_{1,1}(t)\|_{L^p(\Omega)}\leq C\int_0^t\|X^{N,(1)}_r-Y^N_r\|^2_{L^{2p}(\Omega;H)}\mathrm{d}r+C\delta^{\beta},\;t\in[0,T].$ The proof is omitted.

For the case $\beta\in(0,1)\cup(1,2),$ we need to further split the term $K_{1,1}$ into two parts, which are denoted by $K_{1,1,i},i=1,2$, based on 
\begin{align*}
X^{N,(1)}_s-Y^N_s&=\int_0^s\Big(S^N(s-t_{m_r})F^N(X^{N,(1)}_{t_{m_r}})-S^N(s-r)F^N(X(r))\Big)\chi_{\{\tau^{\delta}_{m_r}\ge \tau^{\delta}_{min}\}}\mathrm{d}r\\
&\quad+\int_0^s\Big(S^N(s-t_{m_r})\frac{F^N(X^{N,(1)}_{t_{m_r}})}{1+\|F^N(X^{N,(1)}_{t_{m_r}})\|\tau^{\delta}_{min}}-S^N(s-r)F^N(X(r))\Big)\chi_{\{\tau^{\delta}_{m_r}<\tau^{\delta}_{min}\}}\mathrm{d}r.
\end{align*}
Namely, 
\begin{align*}
K_{1,1,1}(t)&:=\int_0^t\Big\langle\int_0^s\Big(S^N(s-t_{m_r})F^N(X^{N,(1)}_{t_{m_r}})-S^N(s-r)F^N(X(r))\Big)\chi_{\{\tau^{\delta}_{m_r}\ge \tau^{\delta}_{min}\}}\mathrm{d}r,\\
&\quad (S^N(s-t_{m_s})-\mathrm{Id})F^N(X^{N,(1)}_{t_{m_s}})\Big\rangle\chi_{\{\tau^{\delta}_{m_s}\ge\tau^{\delta}_{min}\}}\mathrm{d}s,\\
K_{1,1,2}(t)&:= \int_0^t\Big\langle\int_0^s\Big(S^N(s-t_{m_r})\frac{F^N(X^{N,(1)}_{t_{m_r}})}{1+\|F^N(X^{N,(1)}_{t_{m_r}})\|\tau^{\delta}_{min}}-S^N(s-r)F^N(X(r))\Big)\chi_{\{\tau^{\delta}_{m_r}<\tau^{\delta}_{min}\}}\mathrm{d}r,\\
&\quad (S^N(s-t_{m_s})-\mathrm{Id})F^N(X^{N,(1)}_{t_{m_s}})\Big\rangle\chi_{\{\tau^{\delta}_{m_s}\ge\tau^{\delta}_{min}\}}\mathrm{d}s.
\end{align*}
 For the term $K_{1,1,1},$ it can be further split  as  
\begin{align*}
K_{1,1,1}(t)&
= \int_0^t\Big\langle \int_0^s S^N(s-r)(S^N(r-t_{m_r})-\mathrm{Id})F^N(X^{N,(1)}_{t_{m_r}})\chi_{\{\tau^{\delta}_{m_r}\ge \tau^{\delta}_{min}\}}\mathrm{d}r, \\
&\quad (S^N(s-t_{m_s})-\mathrm{Id})F^N(X^{N,(1)}_{t_{m_s}})\Big\rangle \chi_{\{\tau^{\delta}_{m_s}\ge \tau^{\delta}_{min}\}}\mathrm{d}s\\
&\quad +\int_0^t\Big\langle \int_0^s S^N(s-r)\big(F^N(X^{N,(1)}_{t_{m_r}})-F^N(X^{N,(1)}_r)\big)\chi_{\{\tau^{\delta}_{m_r}\ge \tau^{\delta}_{min}\}}\mathrm{d}r, \\
&\quad (S^N(s-t_{m_s})-\mathrm{Id})F^N(X^{N,(1)}_{t_{m_s}})\Big\rangle \chi_{\{\tau^{\delta}_{m_s}\ge \tau^{\delta}_{min}\}}\mathrm{d}s\\
&\quad +\int_0^t\Big\langle \int_0^s S^N(s-r)\big(F^N(X^{N,(1)}_{r})-F^N(X(r))\big)\chi_{\{\tau^{\delta}_{m_r}\ge\tau^{\delta}_{min}\}}\mathrm{d}r,\\
&\quad (S^N(s-t_{m_s})-\mathrm{Id})F^N(X^{N,(1)}_{t_{m_s}})\Big\rangle \chi_{\{\tau^{\delta}_{m_s}\ge \tau^{\delta}_{min}\}}\mathrm{d}s=:K_{1,1,1}^1(t)+K_{1,1,1}^2(t)+K_{1,1,1}^3(t).
\end{align*}
Similar to estimates of $J^1_{1,1}$ and $J^3_{1,1}$, terms $K_{1,1,1}^1$ and $K_{1,1,1}^3$ can be estimated as
$\|K_{1,1,1}^1(t)\|_{L^p(\Omega)}+\|K_{1,1,1}^3(t)\|_{L^p(\Omega)}\leq C\int_0^t\|X^{N,(1)}_r-Y^N_r\|^2_{L^{2p}(\Omega;H)}\mathrm{d}r+C\lambda_N^{-\beta}+C\delta^{\beta}$ for $\beta\in(0,1)\cup(1,2)$ and $t\in[0,T].$ The proof is omitted.

For the term $K^2_{1,1,1},$ when $\beta\in(0,1),$ the estimate is still similar as $J^2_{1,1}$ and $\|K^2_{1,1,1}(t)\|_{L^p(\Omega)}\leq C\delta^{\beta}.$ 
And when $\beta\in(1,2)$,
by the Taylor formula, and the fact that
\begin{align*}
X^{N,(1)}_r-X^{N,(1)}_{t_{m_r}}&=(S^N(r-t_{m_r})-\mathrm{Id})X^{N,(1)}_{t_{m_r}}+\int_{t_{m_r}}^rS^N(r-t_{m_r})\Big(F^N(X^{N,(1)}_{t_{m_u}})\chi_{\{\tau^{\delta}_{m_u}\ge\tau^{\delta}_{min}\}}\\
&\quad +\frac{F^N(X^{N,(1)}_{t_{m_u}})}{1+\|F^N(X^{N,(1)}_{t_{m_u}})\|\tau^{\delta}_{min}}\chi_{\{\tau^{\delta}_{m_u}<\tau^{\delta}_{min}\}}\Big)\mathrm{d}u+\int_{t_{m_r}}^rS^N(r-t_{m_r})\mathrm{d}W(u),
\end{align*}
we further split $K^2_{1,1,1}$ as 
\begin{align*}
&\quad K_{1,1,1}^2(t)\\
&= -\int_0^t\Big\langle\int_0^sA^{\frac{\beta-1}{2}}S^N(s-r)DF(X^{N,(1)}_{t_{m_r}})(S^N(r-t_{m_r})-\mathrm{Id})X^{N,(1)}_{t_{m_r}}\chi_{\{\tau^{\delta}_{m_r}\ge\tau^{\delta}_{min}\}}\mathrm{d}r,\\
&\qquad A^{-\frac{\beta}{2}}(S^N(s-t_{m_s})-\mathrm{Id})A^{\frac12}F^N(X^{N,(1)}_{t_{m_s}})\Big\rangle \chi_{\{\tau^{\delta}_{m_s}\ge \tau^{\delta}_{min}\}}\mathrm{d}s\\
&\quad -\int_0^t\Big\langle\int_0^sA^{\frac{\beta-1}{2}}S^N(s-r)DF(X^{N,(1)}_{t_{m_r}})\int_{t_{m_r}}^rS^N(r-t_{m_r})\Big(F^N(X^{N,(1)}_{t_{m_u}})\chi_{\{\tau^{\delta}_{m_u}\ge\tau^{\delta}_{min}\}}\\
&\quad +\frac{F^N(X^{N,(1)}_{t_{m_u}})}{1+\|F^N(X^{N,(1)}_{t_{m_u}})\|\tau^{\delta}_{min}}\chi_{\{\tau^{\delta}_{m_u}<\tau^{\delta}_{min}\}}\Big)\mathrm{d}u\chi_{\{\tau^{\delta}_{m_r}\ge\tau^{\delta}_{min}\}}\mathrm{d}r,\\
&\qquad A^{-\frac{\beta}{2}}(S^N(s-t_{m_s})-\mathrm{Id})A^{\frac12}F^N(X^{N,(1)}_{t_{m_s}})\Big\rangle \chi_{\{\tau^{\delta}_{m_s}\ge \tau^{\delta}_{min}\}}\mathrm{d}s\\
&\quad -\int_0^t\Big\langle\int_0^sA^{\frac{\beta-1}{2}}S^N(s-r)DF(X^{N,(1)}_{t_{m_r}})\int_{t_{m_r}}^rS^N(r-t_{m_r})\mathrm{d}W(u)\chi_{\{\tau^{\delta}_{m_r}\ge\tau^{\delta}_{min}\}}\mathrm{d}r,\\
&\qquad A^{-\frac{\beta}{2}}(S^N(s-t_{m_s})-\mathrm{Id})A^{\frac12}F^N(X^{N,(1)}_{t_{m_s}})\Big\rangle \chi_{\{\tau^{\delta}_{m_s}\ge \tau^{\delta}_{min}\}}\mathrm{d}s\\
&\quad -\int_0^t\Big\langle\int_0^sS^N(s-r)R_F(X^{N,(1)}_{t_{m_r}},X^{N,(1)}_r)\chi_{\{\tau^{\delta}_{m_r}\ge\tau^{\delta}_{min}\}}\mathrm{d}r, \;(S^N(s-t_{m_s})-\mathrm{Id})F^N(X^{N,(1)}_{t_{m_s}})\Big\rangle \chi_{\{\tau^{\delta}_{m_s}\ge \tau^{\delta}_{min}\}}\mathrm{d}s\\
&=:\widetilde{II}_1(t)+\widetilde{II}_2(t)+\widetilde{II}_3(t)+\widetilde{II}_4(t).
\end{align*}
 Proofs of terms $\widetilde{II}_1,\widetilde{II}_3,\widetilde{II}_4$ are similar to those of $II_1,II_3,II_4,$ and $\|\widetilde{II}_1(t)\|_{L^p(\Omega)}+\|\widetilde{II}_3(t)\|_{L^p(\Omega)}+\|\widetilde{II}_4(t)\|_{L^p(\Omega)}\leq C\delta^{\beta}$ for $\beta\in(1,2).$ 
For the term $\widetilde{II}_2,$ it follows from \eqref{S(t)1}-\eqref{S(t)2}, \eqref{DF_prop}, and Corollary \ref{coro2} that
\begin{align*}
\|\widetilde{II}_2(t)\|_{L^p(\Omega)}&\leq C\delta^{1+\frac{\beta}{2}}\int_0^t\int_0^s(s-r)^{-\frac{\beta-1}{2}}\Big\|(1+\|X^{N,(1)}_{t_{m_r}}\|^2_E)\|F(X^{N,(1)}_{t_{m_r}})\|\chi_{\{\tau^{\delta}_{m_r}\ge\tau^{\delta}_{min}\}}\Big\|_{L^{2p}(\Omega)}\times\\
&\quad \|F(X^{N,(1)}_{t_{m_s}})\|_{L^{2p}(\Omega;\dot{H}^1)}\mathrm{d}r\mathrm{d}s\leq C\delta^{\beta}.
\end{align*} 
Hence, we get that for $\beta\in(0,1)\cup(1,2),$ 
$ \|K_{1,1,1}^2(t)\|_{L^p(\Omega)}\leq C\delta^{\beta},\;t\in[0,T].$

Similarly, the term $K_{1,1,2}$ can be split as
\begin{align*}
K_{1,1,2}(t)&=\int_0^t\Big\langle \int_0^sS^N(s-r)(S^N(r-t_{m_r})-\mathrm{Id})\frac{F^N(X^{N,(1)}_{t_{m_r}})}{1+\|F^N(X^{N,(1)}_{t_{m_r}})\|\tau^{\delta}_{min}}\chi_{\{\tau^{\delta}_{m_r}<\tau^{\delta}_{min}\}}\mathrm{d}r,\\
&\quad (S^N(s-t_{m_s})-\mathrm{Id})F^N(X^{N,(1)}_{t_{m_s}})\Big\rangle\chi_{\{\tau^{\delta}_{m_s}\ge \tau^{\delta}_{min}\}}\mathrm{d}s\\
&\quad +\int_0^t\Big\langle \int_0^sS^N(s-r)(F^N(X^{N,(1)}_{t_{m_r}})-F^N(X(r)))\chi_{\{\tau^{\delta}_{m_r}<\tau^{\delta}_{min}\}}\mathrm{d}r,\\
&\quad (S^N(s-t_{m_s})-\mathrm{Id})F^N(X^{N,(1)}_{t_{m_s}})\Big\rangle\chi_{\{\tau^{\delta}_{m_s}\ge \tau^{\delta}_{min}\}}\mathrm{d}s\\
&\quad -\int_0^t\Big\langle \int_0^sS^N(s-r)\frac{\|F^N(X^{N,(1)}_{t_{m_r}})\|\tau^{\delta}_{min}F^N(X^{N,(1)}_{t_{m_r}})}{1+\|F^N(X^{N,(1)}_{t_{m_r}})\|\tau^{\delta}_{min}}\chi_{\{\tau^{\delta}_{m_r}<\tau_{min}\}}\mathrm{d}r,\\
&\quad (S^N(s-t_{m_s})-\mathrm{Id})F^N(X^{N,(1)}_{t_{m_s}})\Big\rangle\chi_{\{\tau^{\delta}_{m_s}\ge \tau^{\delta}_{min}\}}\mathrm{d}s=:K^1_{1,1,2}(t)+K^2_{1,1,2}(t)+K^3_{1,1,2}(t).
\end{align*}
Terms $K^1_{1,1,2},K^2_{1,1,2}$ can be estimated similarly to terms $K^{1}_{1,1,1}$ and $K^2_{1,1,1}+K^3_{1,1,1},$ respectively, and we can get 
$\|K^{1}_{1,1,2}(t)\|_{L^p(\Omega)}+\|K^{2}_{1,1,2}(t)\|_{L^p(\Omega)}\leq C\int_0^t\|X^{N,(1)}_s-Y^N_s\|^2_{L^{2p}(\Omega;H)}\mathrm{d}s+C\lambda_N^{-\beta}+C\delta^{\beta}$ for $\beta\in(0,1)\cup(1,2).$ 
Then we show the estimate of the term $K^3_{1,1,2}.$ By \eqref{S(t)1}-\eqref{S(t)2} and the Young inequality, we arrive at for $\beta\in(0,2),$
\begin{align*}
K^3_{1,1,2}(t)&=-\int_0^t\Big\langle \int_0^sA^{\frac{\beta}{2}}S^N(s-r)\frac{\|F^N(X^{N,(1)}_{t_{m_r}})\|\tau^{\delta}_{min}F^N(X^{N,(1)}_{t_{m_r}})}{1+\|F^N(X^{N,(1)}_{t_{m_r}})\|\tau^{\delta}_{min}}\chi_{\{\tau^{\delta}_{m_r}<\tau^{\delta}_{min}\}}\mathrm{d}r,\\
&\quad \;\,A^{-\frac{\beta}{2}}(S^N(s-t_{m_s})-\mathrm{Id})F^N(X^{N,(1)}_{t_{m_s}})\Big\rangle\chi_{\{\tau^{\delta}_{m_s}\ge \tau^{\delta}_{min}\}}\mathrm{d}s\\
&\leq C\int_0^t\int_0^s(s-r)^{-\frac{\beta}{2}}\Big[\|F(X^{N,(1)}_{t_{m_r}})\|^4(\tau^{\delta}_{min})^2\chi_{\{\tau^{\delta}_{m_r}<\tau^{\delta}_{min}\}}\\
&\quad +(s-t_{m_s})^{\beta}\|F(X^{N,(1)}_{t_{m_s}})\|^2\chi_{\{\tau^{\delta}_{m_s}\ge \tau^{\delta}_{min}\}}\Big]\mathrm{d}r\mathrm{d}s,
\end{align*}
which gives $\|K^3_{1,1,2}(t)\|_{L^p(\Omega)}\leq C\delta^{\beta},\;t\in[0,T].$

Hence, we derive $\|K_{1,1}(t)\|_{L^p(\Omega)}\leq C\int_0^t\|X^{N,(1)}_r-Y^N_r\|^2_{L^{2p}(\Omega;H)}\mathrm{d}r+C\lambda_N^{-\beta}+C\delta^{\beta},\;t\in[0,T]$ for $\beta\in(0,2].$

The estimate of $K_{1,2}$ is similar to that of $J_{1,2},$ and we can obtain $\|K_{1,2}(t)\|_{L^p(\Omega)}\leq C\int_0^t\|X^{N,(1)}_s-Y^N_s\|^2_{L^{2p}(\Omega;H)}\mathrm{d}s+C\lambda_N^{-\beta}+C\delta^{\delta}$ for $\beta\in(0,2]$ and $t\in[0,T].$ The proof is omitted. 

Based on Corollary \ref{coro2},  similarly to the estimate of the term $K_{1,1,2}$ and \cite[Theorem 4.12]{WXJ20}, the term $K_2$ can be estimated as 
$\|K_2(t)\|_{L^p(\Omega)}\leq C\int_0^t\|X^{N,(1)}_s-Y^N_s\|^2_{L^{2p}(\Omega;H)}\mathrm{d}s+C\lambda_N^{-\beta}+C\delta^{\delta}$ for $\beta\in(0,2]$ and $t\in[0,T].$
The proof  is omitted. 

Therefore, combining estimates of terms $K_1$ and $K_2$ yields 
$$
\big\|X^{N,(1)}_t-Y^N_t\big\|^2_{L^{2p}(\Omega;H)}\leq C\Big(\int_0^t\big\|X^{N,(1)}_t-Y^N_s\big\|^2_{L^{2p}(\Omega;H)}\mathrm{d}s+\lambda_N^{-\beta}+\delta^{\beta}\Big),
$$
where $p\ge 1$ and the constant $C>0$ depends on $p,T,\tau_{min},L_0,L_1,L_2,X_0,$ and $\|A^{\frac{\beta-1}{2}}Q^{\frac{1}{2}}\|_{\mathcal{L}_2(H)}$.
Recalling that $\lambda_N\sim N^{-\frac 2d}$, and applying the Gr\"onwall inequality, we obtain
$$
\sup_{0\leq t\leq T}\|X^{N,(1)}_t-Y^N_t\|_{L^{2p}{(\Omega;H)}}\leq C(N^{-\frac{\beta}{d}}+\delta^{\frac{\beta}{2}}),\;p\ge 1,
$$
which combining \eqref{X-P^NX}-\eqref{eq4} 
finishes the proof of \eqref{ATEU}.
\end{proof}

%\section*{Acknowledgments}
%This work is funded by the National key R$\&$D Program of China under Grant (No. 2020YFA0713701), National Natural Science Foundation of China (No. 12022118, No. 12031020, No. 11971470,
%No. 11871068), and by Youth Innovation Promotion Association CAS.

\bibliographystyle{plain}
\bibliography{adaptive.bib}

\end{document}